\documentclass[a4paper,11pt]{article}
\usepackage{xypic}
\usepackage{mathtools}
\usepackage{amsmath}
\usepackage{mathrsfs}
\usepackage{amssymb}
\usepackage{amsthm}
\usepackage{dsfont}
\usepackage{graphicx}
\usepackage{bmpsize}
\usepackage{times}
\usepackage{ textcomp }
\usepackage{color}
\usepackage{mathtools}
\usepackage{latexsym, empheq, fancybox}
\usepackage{mathrsfs}
\usepackage{exscale}
\usepackage{indentfirst}
\usepackage{booktabs, array}
\usepackage[english]{babel}
\newtheorem{theorem}{Theorem}[section]
\newtheorem{lemma}[theorem]{Lemma}
\newtheorem{definition}[theorem]{Definition}

\newtheorem{corollary}[theorem]{Corollary}

\newtheorem{remark}[theorem]{Remark}
\newtheorem{proposition}[theorem]{Proposition}
\newtheorem{question}[theorem]{Question}
\usepackage{authblk}
\title{\Large\textbf{Weighted Topological Entropy of Random Dynamical Systems}}

\author{
		Kexiang Yang$^1$,
		Ercai Chen$^1$,
		Zijie Lin$^{1,2}$,
		Xiaoyao Zhou$^1$
		\thanks{ 1. School of Mathematical Sciences and Institute of Mathematics, Nanjing Normal University, Nanjing 210023, PR China\\
			2. School of Mathematics, University of Science and Technology of China, Hefei, Anhui, 230026, PR China
			(E-mail:kxyangs@163.com, ecchen@njnu.edu.cn, zjlin137@126.com, zhouxiaoyaodeyouxian@126.com)}
}

\begin{document}
\date{}

\maketitle
\begin{abstract}
Let $f_{i},i=1,2$ be continuous bundle random dynamical systems over an ergodic compact metric system $(\Omega,\mathcal{F},\mathbb{P},\vartheta)$. Assume that ${\bf a}=(a_{1},a_{2})\in\mathbb{R}^{2}$ with $a_{1}>0$ and $a_{2}\geq0$, $f_{2}$ is a factor of $f_{1}$ with a factor map $\Pi:\Omega\times X_{1}\rightarrow\Omega\times X_{2}$. We define the ${\bf a}$-weighted Bowen topological entropy of $h^{{\bf a}}(\omega,f_{1},X_{1})$ of $f_{1}$ with respect to $\omega\in \Omega$. It is shown that the quality $h^{{\bf a}}(\omega,f_{1},X_{1})$ is measurable in $\Omega$, and denoted that $h^{{\bf a}}(f_{1},\Omega\times X_{1})$ is the integration of $h^{{\bf a}}(\omega,f_{1},X_{1})$ against $\mathbb{P}$. We prove the following variational principle:
\begin{align*}
h^{{\bf a}}(f_{1},\Omega\times X_{1})=\sup\left\{a_{1}h_{\mu}^{(r)}(f_{1})+a_{2}h_{\mu\circ\Pi^{-1}}^{(r)}(f_{2})\right\},
\end{align*}
where the supremum is taken over the set of all $\mu\in\mathcal{M}_{\mathbb{P}}^{1}(\Omega\times X_{1},f_{1})$. In the case of random dynamical systems with an ergodic and compact driving system, this gives an affirmative answer to the question posed by Feng and Huang [Variational principle for weighted topological pressure, J. Math. Pures Appl. {\bf 106} (2016), 411-452]. It also generalizes the relativized variational principle for fiber
topological entropy, and provides a topological extension of Hausdorff dimension of invariant sets and random measures on the $2$-torus $\mathbb{T}^{2}$. In addition, the Shannon-McMillan-Breiman theorem, Brin-Katok local entropy formula and Katok entropy formula of weighted measure-theoretic entropy for random dynamical systems are also established in this paper.
\end{abstract}


\section{Introduction}

The framework of random transformations were studied by Ulam and von Neumann in \cite{Ulam1945Random} and by Kakutani in \cite{Kakutani1951Random}. With the development of random transformations, the ergodic theory of random dynamical systems was motivated in \cite{Kifer1986Ergodic}. Entropy is a key notion in ergodic theory and dynamical systems. For the deterministic case, in 1958 Kolmogorov \cite{Kolmogorov1958A} and Sinai \cite{Sinai1959On} introduced measure-theoretic entropy for measurable dynamical systems. Later, Adler, Konheim and McAndrew \cite{Adler1965Topological} introduced the notion of topological entropy, which is an invariant of topological conjugacy for topological dynamical systems. The
variational principle between topological entropy and
measure-theoretic entropy was proved by Goodwyn \cite{Goodwyn1969Topological}, Dinaburg \cite{Dinaburg1970The} and Goodman \cite{Goodman1971Relating}. The relativized variational principle plays a vital role in the study of equilibrium states for random dynamical systems. Later this topic was discussed in the work of the relativized ergodic theory in \cite{Ledrappier1977A} and then it was extended in \cite{Bogenschutz1992Entropy} to random transformations acting on one space. After that Kifer \cite{Kifer2001On} gave the relativized variational principle in random dynamical systems for the invariant measurable subset which the fibers are compact. For the other works, we refer to \cite{Ledrappier1988Dimension,Kifer1996Fractal,Feng2009Multifractal,Huang2017Entropy}.

The concept of random dynamical systems, introduced into the research of dynamics, has constituted an essential component in the complexity description of dynamical systems (see \cite{Arnold1998Random,Bogenschutz1992Entropy,Crauel2002Random}).
Let $(\Omega,\mathcal{F},\mathbb{P},\vartheta)$ be a Polish probability space together with an invertible ergodic $\mathbb{P}$-preserving transformation $\vartheta$. Let $(X,\mathcal{B})$ be a compact metric space together with the distance function $d$ and the Borel $\sigma$-algebra $\mathcal{B}$. A continuous bundle random dynamical system (RDS) $f$ over $(\Omega,\mathcal{F},\mathbb{P},\vartheta)$ is generated by mappings $f_{\omega}:X\rightarrow X$, i.e.,
\begin{align*}
f_{\omega}^{n}=\left\{
\begin{aligned}
f_{\vartheta^{n-1}\omega}\circ f_{\vartheta^{n-2}\omega}\circ\cdots \circ f_{\vartheta\omega}\circ f_{\omega} \ \ \ \ \ \ \ \ \ \ \ &\text{for} \ n>0,\\
Id \ \ \ \ \ \ \ \ \ \ \ \ \ \ \ \ \ \ \ \ \ \ \ \ \ \ \ \ \ \ \ \ \ \ \ \ \ \ \ \ \ \ \ \ \ \ \ \ \ \ \ \ \ \ \ \ \ \ \ &\text{for} \ n=0,
\end{aligned}
\right.
\end{align*}
so that the map $(\omega,x)\mapsto f_{\omega}x$ is measurable and the map $x\mapsto f_{\omega}x$ is continuous for $\mathbb{P}$-a.e. $\omega\in \Omega$. The map
$\Theta:\Omega\times X\rightarrow\Omega\times X$, $\Theta(\omega,x)=(\vartheta\omega,f_{\omega}x)$ is called the skew product transformation. For each $i=1,2$, let $X_{i}$ be a compact metric space, and $f_{i}$ the corresponding continuous bundle RDS over $(\Omega,\mathcal F,\mathbb P,\vartheta)$. By a factor map from $f_{1}$ to $f_{2}$ we mean a map $\Pi:\Omega\times X_{1}\rightarrow\Omega\times X_{2}$ which satisfies that there exists a continuous surjection map $\pi:X_{1}\rightarrow X_{2}$ such that for $\mathbb{P}$-a.e. $\omega\in \Omega$, the following diagram is commutative:
\begin{align*}
\xymatrix{
  X_{1} \ar[d]_{f_{1,\omega}} \ar[r]^{\pi}
                & X_{2} \ar[d]^{f_{2,\omega}}  \\
  X_{1} \ar[r]_{\pi}
                & X_{2}             }
\end{align*}
and $\Pi:\Omega\times X_{1}\rightarrow \Omega\times X_{2},(\omega,x)\mapsto (\omega,\pi x)$ constitutes a measurable map. Let ${\bf a}=(a_{1},a_{2})\in\mathbb{R}^{2}$ with $a_{1}>0$ and $a_{2}\geq0$. The main purpose of this paper is to consider the following:
\begin{question}[Feng and Huang \cite{Feng2008Variational}]
How can one define a meaningful term $h^{(a_{1},a_{2})}(f_{1},\Omega\times X_{1})$ such that the following variational principle holds?
\begin{align}\label{question}
h^{(a_{1},a_{2})}(f_{1},\Omega\times X_{1})=\sup\left\{a_{1}h_{\mu}^{(r)}(f_{1})+a_{2}h_{\mu\circ\Pi^{-1}}^{(r)}(f_{2})\right\},
\end{align}
where the supremum is taken over the set of all $\mu\in\mathcal{M}_{\mathbb{P}}^{1}(\Omega\times X_{1},f_{1})$, and $h_{\mu}^{(r)}(f_{1})$, $h_{\mu\circ\Pi^{-1}}^{(r)}(f_{2})$ are defined as the fiber or relativized entropy of $\mu$ and $\mu\circ\Pi^{-1}$ with respect to $f_{1}$ and $f_{2}$, respectively (see \eqref{fiber or relativized measure-theoretic entropy} in Section \ref{Preliminaries}).
\end{question}

By the relativized variational principle of RDS, the left-hand side of \eqref{question} is equal to $a_{1}h_{top}^{(r)}(f_{1})$ when $a_{2}=0$, where $h_{top}^{(r)}(f_{1})$ is defined as the fiber topological entropy of RDS (see \cite{Kifer2001On}). We mainly consider the case of $a_{2}\neq 0$. Let us turn back to the background of the Question. The interest is to study the dimensions of invariant sets and measures on the $2$-dimensional torus together with diagonal affine expanding maps (see Bedford \cite{Bedford1984Markov} and McMullen \cite{McMullen1984The}). Let $\mathbb{T}^{2}=\mathbb{R}^{2}/\mathbb{Z}^{2}$ be the $2$-dimensional torus and $A=\text{diag}(m_{1},m_{2})$ be an integral diagonal matrix, $2\leq m_{1}<m_{2}$. Let $T$ be the endomorphism, which is defined by $Tu=Au (\text{mod} \ 1)$ for any $u\in\mathbb{T}^{2}$.  Kenyon and Peres \cite{Kenyon1996Measures} generalized McMullen's result to any compact invariant set and obtained the variational principle of the Hausdorff dimension:
\begin{align}\label{Kenyon and Peres}
\text{dim}_{H}(K)=\sup\left\{\frac{1}{\log m_{2}}h_{\eta}(T)+(\frac{1}{\log m_{1}}-\frac{1}{\log m_{2}})h_{\eta\circ\pi^{-1}}(S)\right\}
\end{align}
where the supremum is taken over the set of all $T$-invariant Borel probability measures $\eta$ supported on $K$, $\pi:\mathbb{T}^{2}\rightarrow\mathbb{T}^{1}$, $(x,y)\mapsto x$, and $S:\mathbb{T}^{1}\rightarrow\mathbb{T}^{1}$, $x\mapsto m_{1}x(\text{mod} \ 1)$. When it comes to the weighted topological entropy and pressure, Barral and Feng gave the definition of weighted topological pressure for the case of subshifts over finite alphabets and its equilibrium measures (see \cite{Barral2012Weighted,Feng2011Equilibrium}). After these works, Feng and Huang \cite{Feng2008Variational} defined the weighted topological entropy and pressure, and established a variational principle for factors maps between general topological dynamical systems. In \cite{Feng2008Variational}, they asked a question: does this variational principle admit a randomized version? A question arises naturally whether there is a weighted Bowen topological entropy $h^{(a_{1},a_{2})}(f_{1},\Omega\times X_{1})$ such that the variational principle \eqref{question} holds. The other results of weighted entropy were referred to \cite{Zhao2018Weighted,Wang2019Weighted,Shen2020Weighted}. Recently, Stadlbauer, Suzuki and Varandas \cite{Stadlbauer2020Thermodynamic} studied the relative pressure of RDS by introducing a Carath\'{e}odory structure. For a RDS, the research with Bowen topological entropy has drawn our attention. Combining with the question from Feng and Huang \cite{Feng2008Variational}, the weighted topological entropy and its variational principle of RDSs are good
objects to start with.

In this paper, we give a definition of $h^{{\bf a}}(f_{1},\Omega\times X_{1})$ by following the dimension theory of affine invariant subsets of tori and that of the definitions of Bowen \cite{Bowen1973Topological}, Feng and Huang \cite{Feng2008Variational} in a way resembling Hausdorff dimension. Firstly, we give the definition of $h^{{\bf a}}(\omega,f_{1},Z)$, which is called weighted Bowen topological entropy of $Z$ of $f_{1}$ with respect to $\omega\in \Omega$ in the subsection \ref{Weighted topological entropy}. And then in the Section \ref{Measurability of weighted topological entropy}, we prove that the function $\omega\mapsto h^{{\bf a}}(\omega,f_{1},X_{1})$ is measurable when the space is $X_{1}$. The weighted Bowen topological entropy of $X_{1}$ of $f_{1}$ is given by
\begin{align*}
h^{{\bf a}}(f_{1},\Omega\times X_{1})= \int_{\Omega} h^{{\bf a}}(\omega,f_{1},X_{1})d\mathbb{P}(\omega).
\end{align*}

Inspired by the question in Feng and Huang \cite{Feng2008Variational}, we prove the variational principle \eqref{question} in the randomized version when the driving system is ergodic compact metric system.  Let $(\Omega,\mathcal{F},\mathbb{P},\vartheta)$ be an ergodic compact metric system. Namely, $\Omega$ is a compact metric
space, $\mathcal{F}$ is the Borel $\sigma$-algebra $\mathcal{B}_{\Omega}$ of $\Omega$, $\mathbb{P}$ is an ergodic Borel probability measure on $\Omega$, and $\vartheta:\Omega\rightarrow\Omega$ is an invertible transformation preserving the measure $\mathbb{P}$. The main result of this paper is the following variational principle for weighted Bowen topological entropy of RDSs.

\begin{theorem}\label{the first result}
Let $f_{i},i=1,2$ be continuous bundle RDSs over an ergodic compact metric system $(\Omega,\mathcal{F},\mathbb{P},\vartheta)$. Assume that ${\bf a}=(a_{1},a_{2})\in\mathbb{R}^{2}$ with $a_{1}>0$ and $a_{2}\geq0$, $f_{2}$ is a factor of $f_{1}$ with a factor map $\Pi:\Omega\times X_{1}\rightarrow\Omega\times X_{2}$. Then
\begin{align*}
h^{{\bf a}}(f_{1},\Omega\times X_{1})=\sup\left\{a_{1}h_{\mu}^{(r)}(f_{1})+a_{2}h_{\mu\circ\Pi^{-1}}^{(r)}(f_{2}):\mu\in \mathcal{M}_{\mathbb{P}}^{1}(\Omega\times X_{1},f_{1})\right\}.
\end{align*}
\end{theorem}

Theorem \ref{the first result} provides a weighted version of relativized variational principle for fiber topological entropy in RDSs. In the Section \ref{The proof of Theorem the first result lower bound}, we give the lower bound of Theorem \ref{the first result}. The other interest is the research in the dimension of $\mu_{\omega}$, $\omega\in\Omega$. The works on the theory of dimension formula of RDSs dates back to \cite{Ledrappier1988Dimension}.  Ledrappier and Young \cite{Ledrappier1988Dimension} showed that the dimension of sample measures equals Lyapunov dimension almost everywhere. By using Theorem \ref{the first result}, we can extend Kenyon-Peres' variational principle \eqref{Kenyon and Peres} to special RDSs consisting of skew product expanding maps on the $2$-torus $\mathbb{T}^{2}$. Let $2\leq m_{1}\leq m_{2}$ be integers. Let $\phi$ be a $C^{1}$ real-valued function on $\mathbb{T}^{1}$. Define $T_{1}:\mathbb{T}^{2}\rightarrow \mathbb{T}^{2}$ by
\begin{align*}
T_{1}(x_{1},x_{2})=(m_{1}x_{1},m_{2}x_{2}+\phi(x_{1})).
\end{align*}
Let $K\subset \mathbb{T}^{2}$ be a $T_{1}$-invariant compact set. Let $\pi$ be the canonical projection from $\mathbb{T}^{2}$ to $\mathbb{T}^{1}$. i.e.,
$
\pi(x_{1},x_{2})=x_{1}.
$
Set $X_{1}=K$ and $X_{2}=\pi(K)$. Define $T_{2}:X_{2}\rightarrow X_{2}$, $x_{1}\mapsto m_{1}x_{1}$. Our target systems are the product spaces $(\Omega\times X_{1},\mathcal{F}\times \mathcal{B}_{X_{1}})$ and $(\Omega\times X_{2},\mathcal{F}\times \mathcal{B}_{X_{2}})$, with the measurable transformation $\Theta_{1}:=\vartheta\times T_{1}$ and $\Theta_{2}:=\vartheta\times T_{2}$, respectively, which can be viewed as special RDSs $f_{1}$ and $f_{2}$ over an ergodic compact metric system $(\Omega,\mathcal{F},\mathbb{P},\vartheta)$. Then $f_{2}$ is the factor of $f_{1}$ with a factor map $\Pi:\Omega\times X_{1}\rightarrow \Omega\times X_{2},(\omega,x)\mapsto(\omega,\pi x)$. Define ${\bf a}=(a_{1},a_{2})$ with
\begin{align*}
a_{1}=\frac{1}{\log m_{2}}, \ a_{2}=\frac{1}{\log m_{1}}-\frac{1}{\log m_{2}}.
\end{align*}
Therefore, by the definition of $h^{\bf a}(f_{1}, \Omega\times K)$, we obtain that $h^{\bf a}(\omega,f_{1}, \Omega\times K)$ is a constant for $\mathbb{P}$-a.e. $\omega\in \Omega$ and the constant is equal to $\text{dim}_{H}(K)$. By Theorem \ref{the first result}, we have
\begin{align*}
\text{dim}_{H}(K)=h^{\bf a}(f_{1}, \Omega\times K)=\sup_{\mu\in \mathcal{M}_{\mathbb{P}}^{1}(\Omega\times K,f_{1})}h_{\mu}^{\bf a}(f_{1}),
\end{align*}
where the supremum is attainable at some ergodic $\mu\in \mathcal{M}_{\mathbb{P}}^{1}(\Omega\times K,f_{1})$. By Theorem \ref{Brin-Katok local entropy formula}, we have for $\mathbb{P}$-a.e. $\omega\in \Omega$, $\text{dim}_{H}(\mu_{\omega})=h_{\mu}^{\bf a}(f_{1})$ for each ergodic $\mu\in \mathcal{M}_{\mathbb{P}}^{1}(\Omega\times K,f_{1})$. Then there exists an ergodic $\mu\in \mathcal{M}_{\mathbb{P}}^{1}(\Omega\times K,f_{1})$ with the disintegration of full Hausdorff dimension, i.e.,
\begin{align*}
\text{dim}_{H}(\mu_{\omega})=\text{dim}_{H}(K)
\end{align*}
for $\mathbb{P}$-a.e. $\omega\in \Omega$. This extends the work of Kenyon and Peres \cite{Kenyon1996Measures} in a random case.

In the Section \ref{The proof of Theorem the first result upper bound}, we prove the upper bound of Theorem \ref{the first result}. To prove Theorem \ref{the first result}, we obtain some results with weighted measure-theoretic entropy $h_{\mu}^{\bf a}(f_{1})$ of RDSs, where $h_{\mu}^{\bf a}(f_{1})$ is defined by
\begin{align*}
h_{\mu}^{\bf a}(f_{1})=a_{1}h_{\mu}^{(r)}(f_{1})+a_{2}h_{\mu\circ\Pi^{-1}}^{(r)}(f_{2}).
\end{align*}
The theorem of Shannon-McMillan-Breiman provides a complementary view of entropy, more detailed and more local in nature.
The theorem was conceived by Shannon \cite{shannon1948A} and proved in increasing strength by McMillan \cite{mcmillan1953The} ($L^{1}$-convergence),  Breiman \cite{Breiman1957The}  (almost everywhere convergence) and Chung \cite{chung1961A} (for countable partitions). Bogenschutz \cite{Bogenschutz1992Entropy} established the Shannon-McMillan-Breiman theorem  of RDS. In 2011, Downarowicz gave the full conditional version of the Shannon-McMillan-Breiman Theorem valid for endomorphisms and subinvariant conditional sigma-algebras in the appendix B of the book \cite{Downarowicz2011Entropy}.
In order to establish the variational principle of weighted topological entropy and pressure, Feng and Huang \cite{Feng2008Variational} provided a weighted version of the  Shannon-McMillan-Breiman Theorem. In the Section \ref{Appendix A: weighted Shannon-McMillan-Breiman Theorem} of this paper, we prove the weighted Shannon-McMillan-Breiman Theorem of RDSs (see Theorem \ref{weighted SMB theorem of RDS}).

As an application of the Shannon-McMillan-Breiman Theorem, Brin and Katok \cite{Brin1983On} gave a topological version of the Shannon-McMillan-Breiman Theorem, that is, Brin-Katok entropy formula, which answered a question posed by Young and Ledrappier.
Later, Zhu \cite{Zhu2008On,Zhu2009Two} proved the random cases of Brin-Katok entropy formula and Katok entropy formula.
Motivated by the Kenyon-Peres's \cite{Kenyon1996Measures} proof of  Ledrappier-Young dimension formula and Brin-Katok's \cite{Brin1983On} argument,  Feng and Huang \cite{Feng2008Variational} also proved the weighted version of Brin-Katok entropy formula, which is also a topological extension of the Ledrappier-Young dimension formula.  In the Section \ref{Appendix B: a weighted version of the Brin-Katok theorem} of the paper, we establish the weighted version of Brin-Katok entropy formula of RDSs (see Theorem \ref{Brin-Katok local entropy formula}) and extend Feng and Huang's result to the random case, i.e., for any $\mu\in \mathcal{E}_{\mathbb{P}}^{1}(\Omega\times X_{1},f_{1})$, then for $\mu$-a.e. $(\omega,x)\in\Omega\times X_{1}$, we have
\begin{align*}
\lim\limits_{\epsilon\to0}\liminf\limits_{n\to\infty}-\dfrac{1}{n}\log\mu_\omega(B^{\bf a}_\omega(x,n,\epsilon))
&=\lim\limits_{\epsilon\to0}\limsup\limits_{n\to\infty}-\dfrac{1}{n}\log\mu_\omega(B^{\bf a}_\omega(x,n,\epsilon))\\&=a_{1}h_{\mu}^{(r)}(f_{1})+a_{2}h_{\mu\circ\Pi^{-1}}^{(r)}(f_{2})\\
&=h^{{\bf a}}_{\mu}(f_{1}).
\end{align*}

For an ergodic Borel probability measure, Katok \cite{Katok1980Lyapunov} defined a topological version of measure-theoretic entropy of a continuous map on a compact metric space. Later, Wang and Huang \cite{Wang2019Weighted} gave the Katok entropy formula with Weighted measure-theoretic entropy for topological dynamical systems. In
the Section \ref{Appendix C: a weighted Katok entropy formula} of the paper, by using cover sets we define the weighted measure-theoretic entropy of RDSs (see Definition \ref{an alternative definition of weighted measure-theoretic entropy}) and establish its Katok entropy formula (see Theorem \ref{weighted Katok entropy formula}), i.e.,
for any $\mu\in \mathcal{E}_{\mathbb{P}}^{1}(\Omega\times X_{1},f_{1})$ and $0<\delta<1$, then for $\mathbb{P}$-a.e. $\omega\in\Omega$, we have
\begin{align*}
\lim_{\epsilon\rightarrow0}\liminf_{n\rightarrow\infty}\frac{1}{n}\log N_{f_{1}}(\omega,{\bf a},n,\epsilon,\delta)&=\lim_{\epsilon\rightarrow0}\limsup_{n\rightarrow\infty}\frac{1}{n}\log N_{f_{1}}(\omega,{\bf a},n,\epsilon,\delta)\\&={\rm ent}_{\mu}^{\bf a}(\omega,f_{1})\\&=h_\mu^{{\bf a}}(f_{1}).
\end{align*}

The paper is organized as follows. In Section \ref{Preliminaries}, we introduce the definitions of fiber measure-theoretic entropy, weighted topological entropy and average weighted topological entropy of RDSs. We prove the upper semi-continuity of some entropy functions, the equivalence of $h^{{\bf a}}$ and $h^{{\bf a}}_{W}$ and a dynamical Frostman lemma of RDSs. In Section \ref{Measurability of weighted topological entropy}, we prove the measurability of $h^{{\bf a}}(\omega,f_{1},X_{1})$ in $\Omega$. In the Section \ref{The proof of Theorem the first result lower bound}, we give the lower bound of Theorem \ref{the first result}. In the Section \ref{The proof of Theorem the first result upper bound}, we prove the upper bound of Theorem \ref{the first result}.  In the Section \ref{Appendix A: weighted Shannon-McMillan-Breiman Theorem}, we prove the weighted Shannon-McMillan-Breiman Theorem of RDSs. In the Section \ref{Appendix B: a weighted version of the Brin-Katok theorem}, we establish the weighted version of Brin-Katok entropy formula of RDSs. In
the Section \ref{Appendix C: a weighted Katok entropy formula}, we prove the weighted version of Katok entropy formula of RDSs.

\section{Preliminaries}\label{Preliminaries}

Denote by $\mathcal{P}_{\mathbb{P}}(\Omega\times X)$ the space of probability measures on $\Omega\times X$ having the marginal $\mathbb{P}$ on $\Omega$. For each measurable in $(\omega,x)$ and continuous in $x\in X$ function $\varphi$ on $\Omega\times X$, let us set
\begin{align*}
\|\varphi\|_{1}=\int\|\varphi(\omega)\|_{\infty}d\mathbb{P}(\omega), \ \text{where} \ \|\varphi(\omega)\|_{\infty}=\sup_{x\in X}|\varphi(\omega,x)|.
\end{align*}
Let $L^{1}(\Omega,C(X))$ the space of such functions $\varphi$ with $\|\varphi\|_{1}<\infty$. For $\mu,\mu^{(n)}\in\mathcal{P}_{\mathbb{P}}(\Omega\times X),n=1,2,\ldots$, write $\mu^{(n)}\Rightarrow \mu$ if $\int \varphi d\mu^{(n)}\rightarrow \int \varphi d\mu$ as $n\rightarrow\infty$ for $\varphi\in L^{1}(\Omega,C(X))$ that introduces a weak* topology in $\mathcal{P}_{\mathbb{P}}(\Omega\times X)$. By [\cite{Kifer2001On}, Lemma 2.1 (i)], the space $\mathcal{P}_{\mathbb{P}}(\Omega\times X)$ is compact in this weak* topology.

We denote by $\pi_{\Omega}$ and $\pi_{X}$ the projections from $\Omega\times X$
onto $\Omega$ and $X$, respectively. Any $\mu\in\mathcal{P}_{\mathbb{P}}(\Omega\times X)$ on $\Omega\times X$ disintegrates $d\mu(\omega,x)=d\mu_{\omega}(x)d\mathbb{P}(\omega)$, where $\mu_{\omega}$ are regular conditional probabilities with respect to the $\sigma$-algebra $\pi_{\Omega}^{-1}\mathcal{F}$. Denote by $\mathcal{M}_{\mathbb{P}}^{1}(\Omega\times X,f)$ and $\mathcal{E}_{\mathbb{P}}^{1}(\Omega\times X,f)$ respectively the sets of $\Theta$-invariant measures $\mu\in \mathcal{P}_{\mathbb{P}}(\Omega\times X)$ and $\Theta$-invariant ergodic measures $\mu\in \mathcal{P}_{\mathbb{P}}(\Omega\times X)$. By Bogensch\"{u}tz \cite{Bogenschutz1993Equilibrium}, $\mu$ is $\Theta$-invariant if and only if $f_{\omega}\mu_{\omega}=\mu_{\vartheta\omega}$ for $\mathbb{P}$-a.e. $\omega\in \Omega$. Let $\mu\in M_{\mathbb{P}}^{1}(\Omega\times X,f)$. Let $\mathcal{Q}=\{\mathcal{Q}_{i}\}$ be a finite measurable partition of $\Omega\times X$, and denote $\mathcal{Q}_{\omega}=\{(\mathcal{Q}_{i})_{\omega}\}$, where $(\mathcal{Q}_{i})_{\omega}=\{x\in X:(\omega,x)\in \mathcal{Q}_{i}\}$ is a partition of $X$. The {\em conditional entropy of $\mathcal{Q}$ given the $\sigma$-algebra $\pi_{\Omega}^{-1}\mathcal{F}$} is defined by
\begin{align*}
H_{\mu}(\mathcal{Q} \ | \ \pi_{\Omega}^{-1}\mathcal{F})=-\int\sum_{i}\mu(\mathcal{Q}_{i} \ | \ \pi_{\Omega}^{-1}\mathcal{F})\log \mu(\mathcal{Q}_{i} \ | \ \pi_{\Omega}^{-1}\mathcal{F})d\mu=\int H_{\mu_{\omega}}(\mathcal{Q}_{\omega})d\mathbb{P}(\omega),
\end{align*}
The {\em fiber or relativized measure-theoretic entropy with respect to $\mu\in \mathcal{M}_{\mathbb{P}}^{1}(\Omega\times X,f)$} is defined by the formula
\begin{align}\label{fiber or relativized measure-theoretic entropy}
h_{\mu}^{(r)}(f)=\sup_{\mathcal{Q}}h_{\mu}^{(r)}(f,\mathcal{Q}), \ h_{\mu}^{(r)}(f,\mathcal{Q})=\lim_{n\rightarrow\infty}\frac{1}{n}H_{\mu}(\bigvee_{i=0}^{n-1}(\Theta^{i})^{-1}\mathcal{Q} \ | \ \pi_{\Omega}^{-1}\mathcal{F})
\end{align}
where the supremum is taken over all finite measurable partitions $\mathcal{Q}$ of $\Omega\times X$ into sets $\mathcal{Q}_{i}$ of the form $\mathcal{Q}_{i}=\Omega\times P_{i}$, where $\mathcal{P}=\{P_{i}\}$ is a finite partition of $X$ into measurable sets, and the limit in \eqref{fiber or relativized measure-theoretic entropy} exists in view of subadditivity of conditional entropy (see [\cite{Kifer1986Ergodic}, Section 2.1]). By the subadditive ergodic theorem \cite{Kifer1986Ergodic} and the ergodicity of $\mathbb{P}$, then the limit
\begin{align}\label{full measure-theoretic entropy}
h_{\mu}^{(r)}(f,\mathcal{Q})=\lim_{n\rightarrow \infty}\frac{1}{n} H_{\mu_{\omega}}(\bigvee_{i=0}^{n-1}(f_{\omega}^{i})^{-1}\mathcal{Q}_{\vartheta^{i}\omega})
\end{align}
exists and is constant $\mathbb{P}$-a.e. $\omega\in \Omega$.

\subsection{Weighted topological entropy}\label{Weighted topological entropy}
For $x\in X_{1}, n\in\mathbb N$ and $\epsilon>0,$ denote the ${\bf a}$-weighted Bowen ball $B^{\bf a}_{\omega}(x, n,\epsilon)$ of RDSs by the set of $y\in X_{1}$ satisfying the following two
conditions:
\begin{align*}
d_{1}(f_{1,\omega}^{j}x,f_{1,\omega}^{j}y)&<\epsilon, \ 0\leq j\leq \lceil a_{1}n\rceil-1,\\
d_{2}(f_{2,\omega}^j\pi x,f_{2,\omega}^j\pi y)&<\epsilon, \ 0\leq j\leq \lceil(a_{1}+a_{2})n\rceil-1.
\end{align*}
Here $\lceil u\rceil$ denotes the least integer not less than $u$. For any subset $Z\subseteq X_{1}$, $s\geq 0$ and $N\in\mathbb{N}$, define
\begin{align*}
\Lambda^{{\bf a},s}_{N,\epsilon}(\omega,Z)=\inf\sum_{j}\exp(-sn_{j})
\end{align*}
where the infimum is taken over all countable collections $\Gamma=\{B^{\bf a}_\omega(x_{j},n_{j},\epsilon)\}_{j}$ with $n_{j}\geq N$, and $\bigcup_{j}B^{\bf a}_\omega(x_{j},n_{j},\epsilon)\supseteq Z$. The quantity $\Lambda^{{\bf a},s}_{N,\epsilon}(\omega,Z)$ does not decrease with $N$, hence the following limit exists:
\begin{align*}
\Lambda^{{\bf a},s}_{\epsilon}(\omega,Z)=\lim_{N\rightarrow\infty}\Lambda^{{\bf a},s}_{N,\epsilon}(\omega,Z).
\end{align*}
There exists a critical value of the parameter $s$, which we will denote by $h^{{\bf a}}(\omega,f_{1},Z,\epsilon)$, where $\Lambda^{{\bf a},s}_{\epsilon}(\omega,Z)$ jumps from $\infty$ to $0$, i.e.
\begin{align*}
\Lambda^{{\bf a},s}_{\epsilon}(\omega,Z)=\left\{
\begin{aligned}
0, \ \ \ \ \ \ \ \ \ \ &\text{if} \ s>h^{{\bf a}}(\omega,f_{1},Z,\epsilon),\\
\infty, \ \ \ \ \ \ \ \ &\text{if} \ s<h^{{\bf a}}(\omega,f_{1},Z,\epsilon).
\end{aligned}
\right.
\end{align*}
Since the quantity $h^{{\bf a}}(\omega,f_{1},Z,\epsilon)$ is monotone with respect to $\epsilon$, we let $h^{{\bf a}}(\omega,f_{1},Z)=\lim_{\epsilon\rightarrow0}h^{{\bf a}}(\omega,f_{1},Z,\epsilon)$. We call $h^{{\bf a}}(\omega,f_{1},Z)$ the ${\bf a}$-weighted Bowen topological entropy of $Z$ of $f_{1}$ with respect to $\omega\in \Omega$.

\subsection{Average weighted topological entropy}

For any function $g:X_{1}\rightarrow [0,\infty)$, $s\geq 0$ and $N\in\mathbb{N}$, define
\begin{align*}
\mathcal{W}^{{\bf a},s}_{N,\epsilon}(\omega,g)=\inf\sum_{j}c_{j}\exp(-sn_{j})
\end{align*}
where the infimum is taken over all countable collections $\Gamma=\{B^{\bf a}_\omega(x_{j},n_{j},\epsilon),c_{j}\}_{j}$ with $n_{j}\geq N$, $0<c_{j}<\infty$ and $\sum_{j}c_{j}\chi_{B^{\bf a}_\omega(x_{j},n_{j},\epsilon)}\geq g$ where $\chi_{A}$ denotes the characteristic function of $A$, i.e., $\chi_{A}(x)=1$ if $x\in A$ and $0$ if $x\in X_{1}\setminus A$. For $Z\subseteq X_{1}$, we set $\mathcal{W}^{{\bf a},s}_{N,\epsilon}(\omega,Z)=\mathcal{W}^{{\bf a},s}_{N,\epsilon}(\omega,\chi_{Z})$. The quantity $\mathcal{W}^{{\bf a},s}_{N,\epsilon}(\omega,Z)$ does not decrease with $N$, hence the following limit exists:
\begin{align*}
\mathcal{W}^{{\bf a},s}_{\epsilon}(\omega,Z)=\lim_{N\rightarrow\infty}\mathcal{W}^{{\bf a},s}_{N,\epsilon}(\omega,Z).
\end{align*}
There exists a critical value of the parameter $s$, which we will denote by $h^{{\bf a}}_{W}(\omega,f_{1},Z,\epsilon)$, where $\mathcal{W}^{{\bf a},s}_{\epsilon}(\omega,Z)$ jumps from $\infty$ to $0$, i.e.
\begin{align*}
\mathcal{W}^{{\bf a},s}_{\epsilon}(\omega,Z)=\left\{
\begin{aligned}
0, \ \ \ \ \ \ \ \ \ \ &\text{if} \ s>h^{{\bf a}}_{W}(\omega,f_{1},Z,\epsilon),\\
\infty, \ \ \ \ \ \ \ \ &\text{if} \ s<h^{{\bf a}}_{W}(\omega,f_{1},Z,\epsilon).
\end{aligned}
\right.
\end{align*}
Since the quantity $h^{{\bf a}}_{W}(\omega,f_{1},Z,\epsilon)$ is monotone with respect to $\epsilon$, we let $h^{{\bf a}}_{W}(\omega,f_{1},Z)=\lim_{\epsilon\rightarrow0}h^{{\bf a}}_{W}(\omega,f_{1},Z,\epsilon)$. We call $h^{{\bf a}}_{W}(\omega,f_{1},Z)$ the average ${\bf a}$-weighted Bowen topological entropy of $Z$ of $f_{1}$ with respect to $\omega\in \Omega$.

\begin{proposition}\label{three properties of Bowen topological entropy}
The following results hold.
\begin{itemize}
  \item[(i)] For $Z\subseteq Z'$, then
  \begin{equation*}
  h^{\bf a}(\omega,f_{1},Z)\leq h^{\bf a}(\omega,f_{1}, Z') \ \ \text{and} \ \ h_{W}^{\bf a}(\omega,f_{1}, Z)\leq h_{W}^{\bf a}(\omega,f_{1}, Z').
  \end{equation*}
  \item[(ii)] For $Z\subseteq \bigcup_{i=1}^{\infty}Z_{i}$, $\epsilon>0$, $s \geq 0$ and $N \geq 1$, then
  \begin{equation*}
  \Lambda_{N,\epsilon}^{{\bf a},s}(\omega, Z)\leq \sum_{i=1}^{\infty}\Lambda_{N,\epsilon}^{{\bf a},s}(\omega, Z_{i}).
  \end{equation*}
\end{itemize}
\end{proposition}
\begin{proof}
Those properties follow directly from the definitions. The proof details can be seen in \cite{Pesin1997Dimension}.
\end{proof}

\subsection{Upper semi-continuity of some entropy functions}

In this subsection, we obtain the upper semi-continuity of some entropy functions by using the definition of conditional entropy (see Lemma \ref{upper semi-continuous function}). For any $\epsilon>0$ and $M\in\mathbb{N}$, we define
\begin{align}\label{definition of finite Borel partition}
\mathcal{P}_{X}(\epsilon,M)=\{\alpha:\alpha \ \text{is a finite Borel partition of} \ X \ \text{with} \ \text{diam}(\alpha)<\epsilon,\ \text{card}(\alpha)\leq M\},
\end{align}
where $\text{diam}(\alpha):=\max_{A\in \alpha}\text{diam}(A)$, and $\text{card}(\alpha)$ stands for the cardinality of $\alpha$.
Then we define
\begin{align*}
\mathcal{P}_{X}(\epsilon)=\{\alpha:\alpha \ \text{is a finite Borel partition of} \ X \ \text{with} \ \text{diam}(\alpha)<\epsilon\}.
\end{align*}
\begin{remark}
It is clear that for any $\epsilon>0$, there exists $N:=N(\epsilon)\in\mathbb{N}$ such that $\mathcal{P}_{X}(\epsilon,M)\neq \emptyset$ for any $M\geq N$.
\end{remark}

\begin{lemma}\label{upper semi-continuous function}
Let $f$ be a continuous bundle RDS over $(\Omega,\mathcal{F},\mathbb{P},\vartheta)$. Then
\begin{itemize}
  \item[(1)] If $M\in\mathbb{N}$ with $\mathcal{P}_{X}(\epsilon,M)\neq\emptyset$, then the map
  \begin{align*}
  \mu\in\mathcal{P}_{\mathbb{P}}(\Omega\times X)\mapsto H_{\mu}(\epsilon,M;l):=\inf_{\alpha\in \mathcal{P}_{X}(\epsilon,M)}\frac{1}{l}H_{\mu}(\bigvee_{i=0}^{l-1}\Theta^{-i}\hat{\alpha} \ | \ \pi_{\Omega}^{-1}\mathcal{F})
  \end{align*}
  where $\hat{\alpha}=\{\Omega\times A:A\in \alpha\}$, is upper semi-continuous from $\mathcal{P}_{\mathbb{P}}(\Omega\times X)$ to $[0,\log M]$ for each $l\in\mathbb{N}$.
  \item[(2)] The map
  \begin{align*}
  \mu\in\mathcal{P}_{\mathbb{P}}(\Omega\times X)\mapsto H_{\mu}(\epsilon;l):=\inf_{\alpha\in \mathcal{P}_{X}(\epsilon)}\frac{1}{l}H_{\mu}(\bigvee_{i=0}^{l-1}\Theta^{-i}\hat{\alpha} \ | \ \pi_{\Omega}^{-1}\mathcal{F})
  \end{align*}
  where $\hat{\alpha}=\{\Omega\times A:A\in \alpha\}$, is a bound upper semi-continuous non-negative function for each $l\in\mathbb{N}$.
  \item[(3)] The map
  \begin{align*}
  \mu\in\mathcal{M}_{\mathbb{P}}^{1}(\Omega\times X,f)\mapsto h_{\mu}^{(r)}(f,\epsilon):=\inf_{\alpha\in \mathcal{P}_{X}(\epsilon)}h^{(r)}_{\mu}(f,\hat{\alpha})
  \end{align*}
  where $\hat{\alpha}=\{\Omega\times A:A\in \alpha\}$, is a bound upper semi-continuous non-negative function.
\end{itemize}
\end{lemma}
\begin{proof}
Firstly, we prove (1). Let $M\in\mathbb{N}$ with $\mathcal{P}_{X}(\epsilon,M)\neq\emptyset$, and $l\in\mathbb{N}$. Clearly, the map $H_{\bullet}(\epsilon,M;l)$ is defined from $\mathcal{P}_{\mathbb{P}}(\Omega\times X)$ to $[0,\log M]$. Let $\mu_{0}\in\mathcal{P}_{\mathbb{P}}(\Omega\times X)$. It is sufficient to show that the map $H_{\bullet}(\epsilon,M;l)$  is upper semi-continuous at $\mu_{0}$. Let $\delta>0$. By the definition of $H_{\mu_{0}}(\epsilon,M;l)$, there exists $\alpha\in \mathcal{P}_{X}(\epsilon,M)$ such that
\begin{align}\label{equation 1 upper semi-continuous function}
\frac{1}{l}H_{\mu_{0}}(\bigvee_{i=0}^{l-1}\Theta^{-i}\hat{\alpha} \ | \ \pi_{\Omega}^{-1}\mathcal{F})\leq H_{\mu_{0}}(\epsilon,M;l)+\delta.
\end{align}
Let $\alpha=\{A_{1},\ldots,A_{u}\}$. Then $u\leq M$ and $\text{diam}(A_{i})<\epsilon$ for $i=1,2,\ldots,u$. By [\cite{Walters1982An}, Lemma 4.15], there exists $\delta_{1}=\delta_{1}(u,\delta)>0$ such that whenever $\gamma_{1}=\{E_{1},\ldots,E_{u}\}$, $\gamma_{2}=\{F_{1},\ldots,F_{u}\}$ are two Borel partitions of $X$ with $\sum_{j=1}^{u}\sum_{i=0}^{l-1}\mu_{0}\circ\Theta^{-i}(\Omega\times (E_{j}\triangle F_{j}))<\delta_{1}$, then
\begin{align}\label{equation 2 upper semi-continuous function}
\begin{split}
&\frac{1}{l}\left|H_{\mu_{0}}(\bigvee_{i=0}^{l-1}\Theta^{-i}\mathcal{Q}_{1} \ | \ \pi_{\Omega}^{-1}\mathcal{F})-H_{\mu_{0}}(\bigvee_{i=0}^{l-1}\Theta^{-i}\mathcal{Q}_{2} \ | \ \pi_{\Omega}^{-1}\mathcal{F})\right|\\ \leq &\frac{1}{l}\sum_{i=0}^{l-1}\left|H_{\mu_{0}\circ\Theta^{-i}}(\mathcal{Q}_{1} \ | \ \mathcal{Q}_{2}\vee\pi_{\Omega}^{-1}\mathcal{F})+H_{\mu_{0}\circ\Theta^{-i}}(\mathcal{Q}_{2} \ | \ \mathcal{Q}_{1}\vee\pi_{\Omega}^{-1}\mathcal{F})\right|<\delta,
\end{split}
\end{align}
where $\mathcal{Q}_{1}=\{\Omega\times A:A\in \gamma_{1}\}$ and $\mathcal{Q}_{2}=\{\Omega\times A:A\in \gamma_{2}\}$.
Write $\eta=\sum_{i=0}^{l-1}\mu_{0}\circ\Theta^{-i}$. Next, we are going to construct a Borel partition $\beta=\{B_{1},\ldots,B_{u}\}$ of $X$ so that $\text{diam}(\beta)<\epsilon$, $\sum_{j=1}^{u}\eta(\Omega\times (A_{j}\triangle B_{j}))<\delta_{1}$ and $\eta(\Omega\times \partial\beta)=0$.

In fact, note that $\eta(\Omega\times X)=l<\infty$, hence $\eta(\Omega\times \cdot)$ is regular on $X$. Then there exist open subsets $V_{j}$ of $X$ such that $A_{j}\subseteq V_{j}$, $\text{diam}(V_{j})<\epsilon$ and $\eta(\Omega\times (V_{j}\setminus A_{j}))<\frac{\delta_{1}}{u^{2}}$ for $j=1,\ldots,u$. Clearly, $\mathcal{V}=\{V_{1},\ldots,V_{u}\}$ is an open cover. Let $t>0$ be a Lebesgue number of $\mathcal{V}$. For any $x\in X$, there exists $0<t_{x}\leq \frac{t}{3}$ such that $\eta(\Omega\times \partial B(x,t_{x}))=0$. Then $\{B(x,t_{x}):x\in X\}$ forms an open cover of $X$. Take its finite subcover $\{B(x_{i},t_{x_{i}})\}_{i=1}^{r}$, that is, $\bigcup_{i=1}^{r}B(x_{i},t_{x_{i}})=X$. Obviously, each $B(x_{i},t_{x_{i}})$ is a subset of some $V_{j(i)}$, $j(i)\in\{1,\ldots,u\}$ since $t_{x_{i}}\leq\frac{t}{3}$.

Let $I_{j}=\{i\in\{1,\ldots,r\}:B(x_{i},t_{x_{i}})\subset V_{j}\}$ for $j=1,\ldots,u$. Then $\bigcup_{j=1}^{u}I_{j}=\{1,\ldots,r\}$. Put $B_{1}=\bigcup_{i\in I_{1}}B(x_{i},t_{x_{i}})$ and $B_{j}=(\bigcup_{i\in I_{j}}B(x_{i},t_{x_{i}}))\setminus \bigcup_{m=1}^{j-1}B_{m}$ inductively for $j=2,\ldots,u$. It is clear that $\beta=\{B_{1},\ldots,B_{u}\}$ is a Borel partition of $X$ with $B_{j}\subseteq V_{j}$ and $\eta(\Omega\times \partial B_{j})=0$ for $j=1,\ldots,u$. For any $j=1,\ldots,u$,
\begin{align*}
A_{j}\triangle B_{j}&=(B_{j}\setminus A_{j})\cup(A_{j}\cap(X\setminus B_{j}))\subseteq (V_{j}\setminus A_{j})\cup \bigcup_{k\neq j}(A_{j}\cap B_{k})\\
&\subseteq (V_{j}\setminus A_{j})\cup\bigcup_{k\neq j}(A_{j}\cap V_{k})\subseteq (V_{j}\setminus A_{j})\cup\bigcup_{k\neq j}(A_{j}\cap (V_{k}\setminus A_{k}))\\
&\subseteq \bigcup_{k=1}^{u}(V_{k}\setminus A_{k}).
\end{align*}
Then $\sum_{j=1}^{u}\eta(\Omega\times (A_{j}\triangle B_{j}))\leq u\sum_{k=1}^{u}\eta(\Omega\times (V_{k}\setminus A_{k}))<\delta_{1}$.

By the above process, we have constructed a Borel partition $\beta=\{B_{1},\ldots,B_{u}\}\in \mathcal{P}_{X}(\epsilon,M)$ so that $\sum_{j=1}^{u}\eta(\Omega\times (A_{j}\triangle B_{j}))<\delta_{1}$ and $\eta(\Omega\times \partial\beta)=0$. By \eqref{equation 1 upper semi-continuous function} and \eqref{equation 2 upper semi-continuous function}, we have
\begin{align*}
\frac{1}{l}H_{\mu_{0}}(\bigvee_{i=0}^{l-1}\Theta^{-i}\hat{\beta} \ | \ \pi_{\Omega}^{-1}\mathcal{F})&\leq \frac{1}{l}H_{\mu_{0}}(\bigvee_{i=0}^{l-1}\Theta^{-i}\hat{\alpha} \ | \ \pi_{\Omega}^{-1}\mathcal{F})+\delta\\&\leq H_{\mu_{0}}(\epsilon,M;l)+2\delta
\end{align*}
where $\hat{\beta}=\{\Omega\times A:A\in\beta\}$. Since $\eta(\Omega\times \partial\beta)=0$, we have $\mu_{0}\circ\Theta^{-i}(\Omega\times \partial\beta)=0$ for $i=0,1\ldots,l-1$. Then for $\mathbb{P}$-a.e. $\omega\in \Omega$, we have $(\mu_{0})_{\omega}((f_{\omega}^{i})^{-1}(\partial\beta))=0$. As $\partial (f_{\omega}^{i})^{-1} A\subseteq (f_{\omega}^{i})^{-1} \partial A$ for any $A\subseteq X$, one has $(\mu_{0})_{\omega}( \partial ((f_{\omega}^{i})^{-1}\beta))=0$ for $i=0,1\ldots,l-1$. Hence $(\mu_{0})_{\omega}( \partial (\bigvee_{i=0}^{l-1}(f_{\omega}^{i})^{-1}\beta))=0$ for $\mathbb{P}$-a.e. $\omega\in \Omega$. By Lemma [\cite{Kifer2001On}, Lemma 2.1], we have
\begin{align*}
\limsup_{\mu\rightarrow\mu_{0}}H_{\mu}(\epsilon,M;l)&\leq \limsup_{\mu\rightarrow\mu_{0}}\frac{1}{l}H_{\mu}(\bigvee_{i=0}^{l-1}\Theta^{-i}\hat{\beta} \ | \ \pi_{\Omega}^{-1}\mathcal{F})\\
&\leq \frac{1}{l}H_{\mu_{0}}(\bigvee_{i=0}^{l-1}\Theta^{-i}\hat{\beta} \ | \ \pi_{\Omega}^{-1}\mathcal{F})\\
&\leq H_{\mu_{0}}(\epsilon,M;l)+2\delta.
\end{align*}
Finally letting $\delta\searrow 0$, we obtain that the map $H_{\bullet}(\epsilon,M;l)$  is upper semi-continuous at $\theta_{0}$. This completes the proof of (1).

Now we turn to the proof of (2). Let $l\in\mathbb{N}$. Since $\mathcal{P}_{X}(\epsilon)=\bigcup_{M\in\mathbb{N},\mathcal{P}_{X}(\epsilon,M)\neq\emptyset}\mathcal{P}_{X}(\epsilon,M)$, we have
\begin{align*}
H_{\mu}(\epsilon;l)=\inf_{M\in\mathbb{N},\mathcal{P}_{X}(\epsilon,M)\neq\emptyset}H_{\mu}(\epsilon,M;l)
\end{align*}
for $\mu\in\mathcal{P}_{\mathbb{P}}(\Omega\times X)$. Moreover, by (1) and the fact that the infimum of any family of any family of upper semi-continuous function is again an upper semi-continuous one, we know that the map
\begin{align*}
  \mu\in\mathcal{P}_{\mathbb{P}}(\Omega\times X)\mapsto H_{\mu}(\epsilon;l):=\inf_{\alpha\in \mathcal{P}_{X}(\epsilon)}\frac{1}{l}H_{\mu}(\bigvee_{i=0}^{l-1}\Theta^{-i}\hat{\alpha} \ | \ \pi_{\Omega}^{-1}\mathcal{F})
  \end{align*}
where $\hat{\alpha}=\{\Omega\times A:A\in \alpha\}$, is a bound upper semi-continuous non-negative function for each $l\in\mathbb{N}$. This proves (2).

In the end we prove (3). Note that
\begin{align*}
h_{\mu}^{(r)}(f,\epsilon)&=\inf_{\alpha\in \mathcal{P}_{X}(\epsilon)}h^{(r)}_{\mu}(f,\hat{\alpha})=\inf_{\alpha\in \mathcal{P}_{X}(\epsilon)}\inf_{l\geq 1}\frac{1}{l}H_{\mu}(\bigvee_{i=0}^{l-1}\Theta^{-i}\hat{\alpha} \ | \ \pi_{\Omega}^{-1}\mathcal{F})\\
&=\inf_{l\geq 1}\inf_{\alpha\in \mathcal{P}_{X}(\epsilon)}\frac{1}{l}H_{\mu}(\bigvee_{i=0}^{l-1}\Theta^{-i}\hat{\alpha} \ | \ \pi_{\Omega}^{-1}\mathcal{F})=\inf_{l\geq 1}H_{\mu}(\epsilon;l)
\end{align*}
where $\hat{\alpha}=\{\Omega\times A:A\in \alpha\}$ and $\mu\in\mathcal{M}_{\mathbb{P}}^{1}(\Omega\times X,f)$. Using (2) and the fact that the infimum of any family of any
family of upper semi-continuous function is again an upper semi-continuous one, we know that the map
\begin{align*}
\mu\in\mathcal{M}_{\mathbb{P}}^{1}(\Omega\times X,f)\mapsto h_{\mu}^{(r)}(f,\epsilon)
\end{align*}
is a bound upper semi-continuous non-negative function. This completes the proof.

\end{proof}

\subsection{Equivalence of $h^{{\bf a}}$ and $h^{{\bf a}}_{W}$}

In this subsection, we prove that the ${\bf a}$-weighted Bowen topological entropy is equal to the average ${\bf a}$-weighted Bowen topological entropy of $\Omega\times X_{1}$ of $f_{1}$ with respect to $\omega\in \Omega$.. The main result of this subsection is the following:
\begin{proposition}\label{Equivalence of weighted and average topological entropy}
Let $Z\subseteq X_{1}$. Then for any $\omega\in \Omega$, $s\geq 0$ and $\epsilon,\delta>0$, we have
\begin{align*}
\Lambda^{{\bf a},s+\delta}_{N,6\epsilon}(\omega,Z)\leq \mathcal{W}^{{\bf a},s}_{N,\epsilon}(\omega,Z)\leq \Lambda^{{\bf a},s}_{N,\epsilon}(\omega,Z),
\end{align*}
where $N$ is large enough. As a consequence, $h^{{\bf a}}(\omega,f_{1},X_{1})=h^{{\bf a}}_{W}(\omega,f_{1},X_{1})$.
\end{proposition}

\begin{remark}
By the above Proposition and Theorem \ref{the measurability of weighted topological entropy}, we obtain that the function $h^{{\bf a}}_{W}(\omega,f_{1},X_{1})$ is measurable in $\Omega$, and $h^{{\bf a}}(f_{1},\Omega\times X_{1})=\int_{\Omega} h^{{\bf a}}_{W}(\omega,f_{1},X_{1})d\mathbb{P}(\omega)$.
\end{remark}

To prove Proposition \ref{Equivalence of weighted and average topological entropy}, we need the covering lemma in the following. The details can be seen in [\cite{Mattila1995Geometry}, Lemma 2.1].

\begin{lemma}\label{covering lemma of RDS}
Let $ \epsilon>0 $ and $ \mathcal{B}_{\omega}(\epsilon)=\{B_{\omega}^{\bf a}(x,n,\epsilon):x\in
X_{1}, n\in\mathbb{N}^{+}\} $. For any family $ \mathcal{B}_{\omega}'(\epsilon)\subseteq\mathcal{B}_{\omega}(\epsilon) $, there exists a (not necessarily countable) subfamily $ \mathcal{G}_{\omega}\subseteq\mathcal{B}_{\omega}'(\epsilon) $ consisting of disjoint balls such that
\begin{align}\label{5 balls covering}
\bigcup_{B\in\mathcal{B}_{\omega}'(\epsilon)}B\subseteq\bigcup_{B_{\omega}^{\bf a}(x,n,\epsilon)\in\mathcal{G}_{\omega}}B_{\omega}^{\bf a}(x,n,5\epsilon).
\end{align}
\end{lemma}

Next, we give the proof of Proposition \ref{Equivalence of weighted and average topological entropy}.

\begin{proof}[Proof of proposition \ref{Equivalence of weighted and average topological entropy}]
Let $Z\subseteq X_{1}$, $s\geq 0$ and $\epsilon,\delta>0$. Take $g=\chi_{Z}$ and $c_{i}=1$, we have
\begin{align*}
\mathcal{W}^{{\bf a},s}_{N,\epsilon}(\omega,Z)\leq \Lambda^{{\bf a},s}_{N,\epsilon}(\omega,Z)
\end{align*}
for any $N\in\mathbb{N}$. In the following, we prove that $
\Lambda^{{\bf a},s+\delta}_{N,6\epsilon}(\omega,Z)\leq \mathcal{W}^{{\bf a},s}_{N,\epsilon}(\omega,Z)$ when $N$ is large enough.

Assume that $ N\geq2 $ is such that $ n^{2}e^{-n\delta}\leq1 $ for $ n\geq
N $. Let $\big\{(B_{\omega}^{\bf a}(x_{i},n_{i},\epsilon),c_{i})\big\}_{i\in
\mathcal{I}}$ be a family so that $ \mathcal{I}\subseteq
\mathbb{N} $, $ x_{i}\in
X_{1}$, $ 0<c_{i}<\infty $, $ n_{i}\geq
N $ and
\begin{align*}
\sum_{i}c_{i}\chi_{B_{i}}\geq
\chi_{Z}
\end{align*}
where $B_{i}:=B_{\omega}^{\bf a}(x_{i},n_{i},\epsilon)$.

In order to prove $\Lambda^{{\bf a},s+\delta}_{N,6\epsilon}(\omega,Z)\leq \mathcal{W}^{{\bf a},s}_{N,\epsilon}(\omega,Z)$, it suffices to prove that
\begin{align*}
\Lambda^{{\bf a},s+\delta}_{N,6\epsilon}(\omega,Z)\leq
\sum_{i\in
\mathcal{I}}c_{i}\exp(-sn_{i}).
\end{align*}
Denote $ \mathcal{I}_{n}:=\{i\in
\mathcal{I}:n_{i}=n\} $ and $\mathcal{I}_{n,k}=\{i\in
\mathcal{I}_{n}:i\leq
k\}$ for $n\geq
N$ and $k\in\mathbb{N}$. For the simplicity of notations, set
\begin{equation*}
B_{i}:=B_{\omega}^{\bf a}(x_{i},n_{i},\epsilon) \ \text{and} \ 5B_{i}:=B_{\omega}^{\bf a}(x_{i},n_{i},5\epsilon)
\end{equation*}
for $i\in
\mathcal{I}$. Without loss of generality, assume that $B_{i}\neq
B_{j}$ for $i\neq
j$. For $t>0$, set
\begin{align*}
Z_{n,t}=\left\{x\in
Z:\sum_{i\in
\mathcal{I}_{n}}c_{i}\chi_{B_{i}}(x)>t\right\},
\end{align*}
and
\begin{align*}
 Z_{n,k,t}=\left\{x\in
Z:\sum_{i\in
\mathcal{I}_{n,k}}c_{i}\chi_{B_{i}}(x)>t\right\}.
\end{align*}
In the following, we divide the proof into three steps.

{\bf Step 1.} For any $ n\geq
N $, $ k\in\mathbb{N} $ and $ t>0 $, there exists a finite set $ \mathcal{J}_{n,k,t}\subseteq
\mathcal{I}_{n,k} $ such that the balls $ B_{i} $ $ (i\in
\mathcal{J}_{n,k,t}) $ are pairwise disjoint, $ Z_{n,k,t}\subseteq
\bigcup_{i\in
\mathcal{J}_{n,k,t}}5B_{i} $ and
\begin{align*}
\text{card}(
	\mathcal{J}_{n,k,t})\exp(-sn)\leq
\frac{1}{t}\sum_{i\in
	\mathcal{I}_{n,k}}c_{i}\exp(-sn).
\end{align*}

In order to prove the previous inequality, using Federer's method \cite{Federer1969Geometric} (see also Mattila \cite{Mattila1995Geometry}), we may assume that each $ c_{i} $ is a positive integer (see \cite{Feng2012Variational}). Let $m $ be the smallest integer satisfying $ m\geq
t $. Let $ \mathscr{B}=\{B_{i}:i\in\mathcal{I}_{n,k}\} $, and $ u:\mathscr{B}\rightarrow\mathbb{Z} $ by $ u(B_{i})=c_{i} $. By induction, we can define integer-valued functions $ v_{0},v_{1},\ldots
,v_{m}$ on $\mathscr{B} $ and subfamilies $ \mathscr{B}_{1},\ldots,\mathscr{B}_{m} $ of $ \mathscr{B} $ with $ v_{0}=u $. Using Lemma \ref{covering lemma of RDS}, there exists a pairwise disjoint subfamily $ \mathscr{B}_{1} $ of $ \mathscr{B} $ such that
$ \bigcup_{B\in
	\mathscr{B}}B\subseteq
\bigcup_{B\in
	\mathscr{B}_{1}}5B $, and then $ Z_{n,k,t}\subseteq
\bigcup_{B\in
	\mathscr{B}_{1}}5B $. Set
\begin{align*}
v_{1}(B)=\left\{
\begin{aligned}
v_{0}(B)-1 \ &\text{for} \ B\in\mathscr{B}_{1},\\
v_{0}(B)\ \ \ \ \ \ \ \ &\text{for} \ B\in\mathscr{B}\setminus
\mathscr{B}_{1},
\end{aligned}
\right.
\end{align*}
and $\mathscr{C}_{1}=\{B\in \mathscr{B}:v_{1}(B)\geq 1\}$. Note that $ \mathscr{B}_{1} $ is pairwise disjoint, we have $ Z_{n,k,t}\subseteq\{x:\sum_{B\in\mathscr{B}:x\in
	B}v_{1}(B)\geq
m-1\} $, where any $ x\in
Z_{n,k,t} $ belongs to some ball $ B\in\mathscr{B} $ with $ v_{1}(B)\geq1 $. Then $ Z_{n,k,t}\subseteq \bigcup_{B\in \mathscr{C}_{1}}B$. By using Lemma \ref{covering lemma of RDS}, there exists a pairwise disjoint subfamily $ \mathscr{B}_{2} $ of $ \mathscr{C}_{1} $ such that $ \bigcup_{B\in
	\mathscr{C}_{1}}B\subseteq
\bigcup_{B\in
	\mathscr{B}_{2}}5B $, and then $ Z_{n,k,t}\subseteq
\bigcup_{B\in
	\mathscr{B}_{2}}5B $. Set
\begin{align*}
v_{2}(B)=\left\{
\begin{aligned}
v_{1}(B)-1 \ &\text{for} \ B\in\mathscr{B}_{2},\\
v_{1}(B)\ \ \ \ \ \ \ \ &\text{for} \ B\in\mathscr{B}\setminus
\mathscr{B}_{2},
\end{aligned}
\right.
\end{align*}
and $\mathscr{C}_{2}=\{B\in \mathscr{B}:v_{2}(B)\geq 1\}$. Note that $ \mathscr{B}_{2} $ is pairwise disjoint, we have $ Z_{n,k,t}\subseteq\{x:\sum_{B\in\mathscr{B}:x\in
	B}v_{2}(B)\geq
m-2\} $, where any $ x\in
Z_{n,k,t} $ belongs to some ball $ B\in\mathscr{B} $ with $ v_{2}(B)\geq1 $. Then $ Z_{n,k,t}\subseteq \bigcup_{B\in \mathscr{C}_{2}}B$. By using Lemma \ref{covering lemma of RDS}, there exists a pairwise disjoint subfamily $ \mathscr{B}_{3} $ of $ \mathscr{C}_{2} $ such that $ \bigcup_{B\in
	\mathscr{C}_{2}}B\subseteq
\bigcup_{B\in
	\mathscr{B}_{3}}5B $, and then $ Z_{n,k,t}\subseteq
\bigcup_{B\in
	\mathscr{B}_{3}}5B $.
Using Lemma \ref{covering lemma of RDS} repeatedly, for $ j=1,\ldots,m $, we can obtain disjoint subfamilies $ \mathscr{B}_{j} $ of $ \mathscr{B}$ such that
\begin{align*}
\mathscr{B}_{j}\subseteq\{B\in
\mathscr{B}:v_{j-1}(B)\geq1\}, \ \ Z_{n,k,t}\subseteq
\bigcup_{B\in
	\mathscr{B}_{j}}5B
\end{align*}
and the function $ v_{j} $ such that
\begin{align*}
v_{j}(B)=\left\{
\begin{aligned}
v_{j-1}(B)-1 \ &\text{for} \ B\in\mathscr{B}_{j},\\
v_{j-1}(B)\ \ \ \ \ \ \ \ &\text{for} \ B\in\mathscr{B}\setminus
\mathscr{B}_{j}.
\end{aligned}
\right.
\end{align*}
This is possible since $ Z_{n,k,t}\subseteq\{x:\sum_{B\in\mathscr{B}:x\in
	B}v_{j}(B)\geq
m-j\} $ for $ j<m $, where any $ x\in
Z_{n,k,t} $ belongs to some ball $ B\in\mathscr{B} $ with $ v_{j}(B)\geq1 $. Hence,
\begin{align*}
\sum_{j=1}^{m}\text{card}(
	\mathscr{B}_{j})\exp(-sn)
&=\sum_{j=1}^{m}\sum_{B\in
	\mathscr{B}_{j}}(v_{j-1}(B)-v_{j}(B))\exp(-sn)\\
&\leq\sum_{B\in
	\mathscr{B}}\sum_{j=1}^{m}(v_{j-1}(B)-v_{j}(B))\exp(-sn)\\
&\leq\sum_{B\in
	\mathscr{B}}u(B)\exp(-sn)\\
&=\sum_{i\in
	\mathcal{I}_{n,k}}c_{i}\exp(-sn).
\end{align*}

Select $ j_{0}\in\{1,2,\ldots,l\} $ such that $ \text{card}(
	\mathscr{B}_{j_{0}}) $ is the smallest. Then we have
\begin{align*}
\text{card}(
	\mathscr{B}_{j_{0}})\exp(-sn)
&\leq
\frac{1}{m}\sum_{i\in
	\mathcal{I}_{n,k}}c_{i}\exp(-sn)\\
&\leq\frac{1}{t}\sum_{i\in
	\mathcal{I}_{n,k}}c_{i}\exp(-sn).
\end{align*}
Hence, $\mathcal{J}_{n,k,t}=\{i\in\mathcal{I}_{n,k}:B_{i}\in\mathcal{B}_{j_{0}}\}$ is as desired.

{\bf Step 2.} For any $ n\geq
N $ and $ t>0 $, then
\begin{align*}
\Lambda_{N,6\epsilon}^{s+\delta}(\omega,Z_{n,t})
\leq\frac{1}{n^{2}t}\sum_{i\in
	\mathcal{I}_{n}}c_{i}\exp(-sn).
\end{align*}

Without loss of generality, suppose that $ Z_{n,t}\neq\emptyset $; otherwise there is nothing to prove. Since $ Z_{n,k,t}\rightarrow
Z_{n,t} $ as $ k\rightarrow\infty $, $ Z_{n,k,t}\neq\emptyset $ for all sufficiently large $ k $. Let $ \mathcal{J}_{n,k,t} $ be the set constructed in step 1, then $ \mathcal{J}_{n,k,t}\neq\emptyset $ for sufficiently large $ k $. Define $ E_{n,k,t}=\{x_{i}:i\in
\mathcal{J}_{n,k,t}\} $. Since the space of all non-empty compact subsets of $ X_{1} $ is compact in terms of the Hausdorff distance (cf. \cite{Federer1969Geometric}), there exists a subsequence $ \{k_{j}\}_{j\geq1} $ of positive integers and a non-empty compact set $ E_{n,t}\subseteq
X_{1} $ such that $ E_{n,k_{j},t} $ converges to $ E_{n,t}$ under the Hausdorff distance as $ j\rightarrow\infty $. Since the distance of any two points in $ E_{n,k,t} $ is not less than $ \epsilon $, so do the points in $ E_{n,t} $. Then $ E_{n,t} $ is a finite set and $ \text{card}(E_{n,k_{j},t})=\text{card}(E_{n,t})$ for sufficiently large $ j $.
Hence
\begin{align*}
\bigcup_{x\in
	E_{n,t}}B_{\omega}^{\bf a}(x,n,5.5\epsilon)\supseteq
\bigcup_{x\in
	E_{n,k_{j},t}}B_{\omega}^{\bf a}(x,n,5\epsilon)=\bigcup_{i\in\mathcal{J}_{n,k_{j},t}}5B_{i}\supseteq
Z_{n,k_{j},t}
\end{align*}
for any sufficiently large $ j $, and then $ \bigcup_{x\in
	E_{n,t}}B_{\omega}^{\bf a}(x,n,6\epsilon)\supseteq
Z_{n,t}$.

Since $ \text{card}(E_{n,k_{j},t})=\text{card}(E_{n,t})$ for sufficiently large $ j $, using the conclusion in Step 1, hence
\begin{align*}
\text{card}(E_{n,t})\exp(-sn)
&=\text{card}(E_{n,k_{j},t})\exp(-sn)\\
&\leq\frac{1}{t}\sum_{i\in
	\mathcal{I}_{n}}c_{i}\exp(-sn).
\end{align*}
Then
\begin{align*}
\Lambda_{N,6\epsilon}^{s+\delta}(\omega,Z_{n,t})
&\leq\text{card}(
	E_{n,t})\exp(-n(s+\delta))\\
&\leq\frac{1}{\exp(n\delta)t}\sum_{i\in
	\mathcal{I}_{n}}c_{i}\exp(-sn)\\
&\leq\frac{1}{n^{2}t}\sum_{i\in
	\mathcal{I}_{n}}c_{i}\exp(-sn).
\end{align*}

{\bf Step 3.} For any $ t\in(0,1) $, then
\begin{align*}
\Lambda_{N,6\epsilon}^{s+\delta}(\omega,Z)\leq\frac{1}{t}\sum_{i\in
	\mathcal{I}}c_{i}\exp(-sn_{i}).
\end{align*}

Fix $ t\in(0,1) $. Observe that $ \sum_{n=N}^{\infty}n^{-2}<1 $ and $ Z\subseteq
\bigcup_{n=N}^{\infty}Z_{n,n^{-2}t} $. By (ii) of Proposition \ref{three properties of Bowen topological entropy} and Step 2, then
\begin{align*}
\Lambda_{N,6\epsilon}^{s+\delta}(\omega,Z)
&\leq\sum_{n=N}^{\infty}\Lambda_{N,6\epsilon}^{s+\delta}(\omega,Z_{n,n^{-2}t})\\
&\leq\sum_{n=N}^{\infty}\frac{1}{t}\sum_{i\in
	\mathcal{I}_{n}}c_{i}\exp(-sn)\\
&\leq\frac{1}{t}\sum_{i\in
	\mathcal{I}}c_{i}\exp(-sn_{i}).
\end{align*}
which finishes the proof.
\end{proof}

\subsection{A dynamical Frostman lemma of RDSs}

In our proof of Theorem \ref{the first result}, we need the following dynamical Frostman lemma of RDSs. The case of $a_{1}=1,a_{2}=0$ was first proved in [\cite{Feng2012Variational}, Lemma 3.4], and then the weighted version of $a_{1}>0,a_{2}\geq0$ was given by [\cite{Feng2008Variational}, Lemma 3.3] for topological dynamical systems.

\begin{lemma}\label{a dynamical Frostman lemma of RDS}
Suppose that $h^{{\bf a}}(\omega_{0},f_{1},X_{1})>0$. Then for any $0<s<h^{{\bf a}}(\omega_{0},f_{1},X_{1})$, there exist a Borel probability measure $\nu$ in $X_{1}$ and $\epsilon>0$, $N\in\mathbb{N}$ such that for any $x\in X_{1}$ and $n\geq N$ we have
\begin{align*}
\nu(B^{\bf a}_{\omega_{0}}(x,n,\epsilon))\leq \exp(-sn).
\end{align*}
\end{lemma}

The proof of Lemma \ref{a dynamical Frostman lemma of RDS} follows Proposition \ref{Equivalence of weighted and average topological entropy} and Lemma \ref{Frostman's Lemma of RDS}.

\begin{lemma}\label{Frostman's Lemma of RDS}
Let $ s\geq 0 $, $ N\in\mathbb{N} $ and $ \epsilon>0 $. Suppose that $c:=\mathcal{W}^{{\bf a},s}_{N,\epsilon}(\omega_{0},X_{1})>0$. Then there is a Borel probability measure $ \nu$ on $ X_{1}$ such that for any $n\geq N$, $x\in X_{1}$, and any compact $K\subset B^{\bf a}_{\omega_{0}}(x,n,\epsilon)$,
\begin{align*}
\nu(K)\leq\frac{1}{c}\exp(-sn).
\end{align*}
\end{lemma}
\begin{proof}
It is clear that $ c<\infty $. Define a function $ p $ on the space $ C(X_{1}) $ of continuous real-valued functions on $X_{1}$ by
\begin{align*}
p(\phi)=\frac{1}{c}\mathcal{W}^{{\bf a},s}_{N,\epsilon}(\omega_{0},\phi).
\end{align*}

Let $ {\bf 1}\in
C(X_{1}) $ denote the constant function $ {\bf 1}(x)\equiv1 $. In the following, it is easy to confirm that
\begin{itemize}
	\item[(1)] $ p(\phi_{1}+\phi_{2})\leq
	p(\phi_{1})+p(\phi_{2}) $ for all $ \phi_{1},\phi_{2}\in
	C(X_{1}) $;
	\item[(2)] $ p(t\phi)=tp(\phi) $ for all $ t\geq0 $ and $ \phi\in
	C(X_{1}) $;
	\item[(3)] $ p({\bf 1})=1 $, $ 0\leq
	p(\phi)\leq\|\phi\|_{\infty} $ for any $ \phi\in
	C(X_{1}) $, and $ p(\phi)=0 $ for $ \phi\in
	C(X_{1}) $ with $ \phi\leq0 $.
\end{itemize}

By the Hahn-Banach theorem, we can extend the linear functional $t\rightarrow t p({\bf 1}), t\in\mathbb{R}$, from
the subspace of the constant functions to a linear functional $ L:C(X_{1})\rightarrow\mathbb{R} $ satisfying
\begin{align*}
L({\bf 1})=p({\bf 1})=1 \ \text{and} \ -p(-\phi)\leq
L(\phi)\leq
p(\phi) \ \text{for any} \ \phi\in
C(X_{1}).
\end{align*}
If $ \phi\in
C(X_{1}) $ with $ \phi\geq0 $, then $ p(-\phi)=0 $ and so $ L(\phi)\geq0 $. Then combining $ L({\bf 1})=1 $, by using the Riesz representation theorem, we can find a Borel probability measure $ \nu $ on $ X_{1} $ such that $ L(\phi)=\int
\phi d\nu $ for $ \phi\in
C(X_{1}) $.

Let $x\in X_{1}$ and $n\geq N$. Suppose that $K$ is a compact subset of $B_{\omega_{0}}^{\bf a}(x,n,\epsilon)$. By the Uryson lemma, there exists $\phi\in C(X_{1})$ such that $ 0\leq
\phi\leq1 $, $ \phi(y)=1 $ for $ y\in
K $ and $ \phi(y)=0 $ for $ y\in
X_{1}\setminus
B_{\omega_{0}}^{\bf a}(x,n,\epsilon) $. Hence $\nu(K)\leq L(\phi)\leq p(\phi)$. Since $\phi\leq
\chi_{B_{\omega_{0}}^{\bf a}(x,n,\epsilon)} $, then $ \mathcal{W}_{N,\epsilon}^{{\bf a},s}(\omega_{0},
\phi)\leq
\exp(-sn) $ and hence $ p(\phi)\leq
\frac{1}{c}\exp(-sn) $. Then $ \nu(K)\leq
\frac{1}{c}\exp(-sn) $. This implies the proof.
\end{proof}

\section{Measurability of weighted topological entropy}\label{Measurability of weighted topological entropy}

In this section, we mainly consider the measurability of $h^{{\bf a}}(\omega,f_{1},X_{1})$ in $\Omega$. We obtain the following theorem which states that the function $h^{{\bf a}}(\omega,f_{1},X_{1})$ is measurable in $\Omega$. Motivated by the measurability methods of Kifer in \cite{Kifer2001On}, we give the proof.

\begin{theorem}\label{the measurability of weighted topological entropy}
Let $f_{i},i=1,2$ be continuous bundle RDSs over $(\Omega,\mathcal{F},\mathbb{P},\vartheta)$. Assume that ${\bf a}=(a_{1},a_{2})\in\mathbb{R}^{2}$ with $a_{1}>0$ and $a_{2}\geq0$, $f_{2}$ is a factor of $f_{1}$ with a factor map $\Pi:\Omega\times X_{1}\to \Omega\times X_{2}$. The function $h^{{\bf a}}(\omega,f_{1},X_{1})$ is measurable in $\Omega$.
\end{theorem}

In the following, we give an equivalent definition of weighted fiber topological entropy by using open
covers. In the following, for any subset $A$ of $X$, set $\text{diam}(A)=\sup_{x,y\in A}d(x, y)$. For a
family of open covers $\{\mathcal{U}_{i}\}_{i=1}^{2}$ where $\mathcal{U}_{i}$ is an open cover of $X_{i}$, define $\text{diam}(\{\mathcal{U}_{i}\}_{i=1}^{2}):=\max_{1\leq i\leq 2}
\text{diam}(\mathcal{U}_{i})$, where $\text{diam}(\mathcal{U}_{i})=\max_{U\in \mathcal{U}_{i}}\text{diam}(U)$. For $i=1, 2$, fix open covers $\{\mathcal{U}_{i}\}_{i=1}^{2}$ where $\mathcal{U}_{i}$ is a finite open cover of $X_{i}$. For ${\bf a }=(a_1,a_2)\in\mathbb R^2$ with $a_1>0$ and $a_2\geq0$, define the weighted string
\begin{align*}
\mathbf{U}^{{\bf a}}:=U^{1}_{1}U^{2}_{1} U^{3}_{1}\cdots U^{\lceil a_{1}n\rceil}_{1}U^{1}_{2} U^{2}_{2} U^{3}_{2}\cdots U^{\lceil(a_{1}+a_{2})n\rceil}_{2},
\end{align*}
where $U^{j}_{1}\in \mathcal{U}_{1}$ for $1\leq j\leq\lceil a_1n\rceil$ and $U^{j}_{2}\in \mathcal{U}_{2}$ for $1\leq j\leq\lceil(a_1+a_2)n\rceil$. Let
\begin{align*}
S_{n}(\{\mathcal{U}_{i}\}_{i=1}^{2}):=\{\mathbf{U}^{{\bf a}}:&U^{j}_{1}\in \mathcal{U}_{1} \ \text{for} \ 1\leq j\leq\lceil a_1n\rceil \ \text{and} \\& U^{j}_{2}\in \mathcal{U}_{2} \ \text{for} \ 1\leq j\leq\lceil(a_1+a_2)n\rceil\}
\end{align*}
and
\begin{align*}
S(\{\mathcal{U}_{i}\}_{i=1}^{2}):=\bigcup_{n>0}S_{n}(\{\mathcal{U}_{i}\}_{i=1}^{2}).
\end{align*}
For $\mathbf{U}^{{\bf a}}\in S_{n}(\{\mathcal{U}_{i}\}_{i=1}^{2})$, we call the integer $m(\mathbf{U}^{{\bf a}})=n$ the length of $\mathbf{U}^{{\bf a}}$ and define
\begin{align*}
X_{1}^{\omega}(\mathbf{U}^{{\bf a}}):=U^{1}_{1}&\cap (f_{1,\omega})^{-1}U^{2}_{1}\cap(f_{1,\omega}^{2})^{-1}U^{3}_{1}\cap\cdots\cap (f_{1,\omega}^{\lceil a_{1}n\rceil-1})^{-1}U^{\lceil a_{1}n\rceil}_{1}\\
&\cap\pi^{-1}U^{1}_{2}\cap\pi^{-1}(f_{2,\omega})^{-1}U^{2}_{2}\cap\pi^{-1}(f_{2,\omega}^{2})^{-1}U^{3}_{2}\cap\cdots\\
&\cap\pi^{-1}(f_{2,\omega}^{\lceil(a_{1}+a_{2})n\rceil-1})^{-1}U^{\lceil(a_{1}+a_{2})n\rceil}_{2}.
\end{align*}

For any subset $Z\subseteq X_{1}$, $s\geq 0$ and $N\in\mathbb{N}$, define
\begin{align*}
\mathcal{M}^{{\bf a},s}_{N,\{\mathcal{U}_{i}\}_{i=1}^{2}}(\omega,Z)=\inf\sum_{\mathbf{U}^{{\bf a}}\in \Gamma}\exp(-sm(\mathbf{U}^{{\bf a}}))
\end{align*}
where the infimum is taken over all countable collections $\Gamma\subset\bigcup_{n\geq N}S_{n}(\{\mathcal{U}_{i}\}_{i=1}^{2})$ such that $\bigcup_{\mathbf{U}^{{\bf a}}\in \Gamma}X_{1}^{\omega}(\mathbf{U}^{{\bf a}})\supseteq Z$. The quantity $\mathcal{M}^{{\bf a},s}_{N,\{\mathcal{U}_{i}\}_{i=1}^{2}}(\omega,Z)$ does not decrease
with $N$, hence the following limit exists:
\begin{align*}
\mathcal{M}^{{\bf a},s}_{\{\mathcal{U}_{i}\}_{i=1}^{2}}(\omega,Z)=\lim_{N\rightarrow\infty}\mathcal{M}^{{\bf a},s}_{N,\{\mathcal{U}_{i}\}_{i=1}^{2}}(\omega,Z).
\end{align*}
There exists a critical value of the parameters,which is defined by $h^{\bf a}(\omega,f_{1},Z,\{\mathcal{U}_{i}\}_{i=1}^{2})$, where $\mathcal{M}^{{\bf a},s}_{\{\mathcal{U}_{i}\}_{i=1}^{2}}(\omega,Z)$ jumps from $\infty$ to $0$, i.e.
\begin{align*}
\mathcal{M}^{{\bf a},s}_{\{\mathcal{U}_{i}\}_{i=1}^{2}}(\omega,Z)=\left\{
\begin{aligned}
0, \ \ \ \ \ \ \ \ \ \ &\text{if} \ s>h^{\bf a}(\omega,f_{1},Z,\{\mathcal{U}_{i}\}_{i=1}^{2}),\\
\infty, \ \ \ \ \ \ \ \ &\text{if} \ s<h^{\bf a}(\omega,f_{1},Z,\{\mathcal{U}_{i}\}_{i=1}^{2}).
\end{aligned}
\right.
\end{align*}

\begin{proposition}\label{the existence of limit}
The limit
\begin{equation*}
\lim_{\text{diam}(\{\mathcal{U}_{i}\}_{i=1}^{2})\rightarrow 0}h^{\bf a}(\omega,f_{1},Z,\{\mathcal{U}_{i}\}_{i=1}^{2})=\sup_{\{\mathcal{U}_{i}\}_{i=1}^{2}}h^{\bf a}(\omega,f_{1},Z,\{\mathcal{U}_{i}\}_{i=1}^{2})
\end{equation*}
exists (but may be $+\infty$).
\end{proposition}
\begin{proof}
Let $\mathcal{V}=\{\mathcal{V}_{1},\mathcal{V}_{2}\}$, where $\mathcal{V}_{i}$ is an open cover of $X_{i}$ and $\mathcal{U}=\{\mathcal{U}_{1},\mathcal{U}_{2}\}$ where $\mathcal{U}_{i}$ is an open cover of $X_{i}$ for $1\leq i\leq k$. Assume that $\mathcal{V}_{i}$ refines $\mathcal{U}_{i}$, this implies that
\begin{align*}
\mathcal{M}^{{\bf a},s}_{\{\mathcal{U}_{i}\}_{i=1}^{2}}(\omega,Z)\leq \mathcal{M}^{{\bf a},s}_{\{\mathcal{V}_{i}\}_{i=1}^{2}}(\omega,Z),
\end{align*}
and then
\begin{align*}
h^{\bf a}(\omega,f_{1},Z,\{\mathcal{U}_{i}\}_{i=1}^{2})\leq h^{\bf a}(\omega,f_{1},Z,\{\mathcal{V}_{i}\}_{i=1}^{2}).
\end{align*}
It follows that
\begin{align*}
h^{\bf a}(\omega,f_{1},Z,\{\mathcal{U}_{i}\}_{i=1}^{2})\leq \liminf_{\text{diam}(\{\mathcal{V}_{i}\}_{i=1}^{2})\rightarrow 0}h^{\bf a}(\omega,f_{1},Z,\{\mathcal{V}_{i}\}_{i=1}^{2}),
\end{align*}
which implies that
\begin{align*}
\limsup_{\text{diam}(\{\mathcal{U}_{i}\}_{i=1}^{2})\rightarrow 0}h^{\bf a}(\omega,f_{1},Z,\{\mathcal{U}_{i}\}_{i=1}^{2})\leq \liminf_{\text{diam}(\{\mathcal{V}_{i}\}_{i=1}^{2})\rightarrow 0}h^{\bf a}(\omega,f_{1},Z,\{\mathcal{V}_{i}\}_{i=1}^{2}),
\end{align*}
and
\begin{align*}
\sup_{\{\mathcal{U}_{i}\}_{i=1}^{2}}h^{\bf a}(\omega,f_{1},Z,\{\mathcal{U}_{i}\}_{i=1}^{2})\leq \liminf_{\text{diam}(\{\mathcal{V}_{i}\}_{i=1}^{2})\rightarrow 0}h^{\bf a}(\omega,f_{1},Z,\{\mathcal{V}_{i}\}_{i=1}^{2}).
\end{align*}
Combining this with the fact that
\begin{align*}
\sup_{\{\mathcal{U}_{i}\}_{i=1}^{2}}h^{\bf a}(\omega,f_{1},Z,\{\mathcal{U}_{i}\}_{i=1}^{2})\geq \liminf_{\text{diam}(\{\mathcal{V}_{i}\}_{i=1}^{2})\rightarrow 0}h^{\bf a}(\omega,f_{1},Z,\{\mathcal{V}_{i}\}_{i=1}^{2}).
\end{align*}
This ends the proof.
\end{proof}

\begin{proposition}\label{equality for weighted topological entropy}
Let $Z$ be a subset of $X_{1}$. Then
\begin{align*}
h^{{\bf a}}(\omega,f_{1},Z)=\sup_{\{\mathcal{U}_{i}\}_{i=1}^{2}}h^{\bf a}(\omega,f_{1},Z,\{\mathcal{U}_{i}\}_{i=1}^{2}).
\end{align*}
\end{proposition}

\begin{proof}
Fix $t>h^{{\bf a}}(\omega,f_{1},Z)$ and a family of open covers $\{\mathcal{U}_{i}\}_{i=1}^{2}$. For any open cover $\mathcal{U}$, let $\gamma(\mathcal{U})$ be the Lebesgue number of $\mathcal{U}$ and define $\gamma(\{\mathcal{U}_{i}\}_{i=1}^{2}):=\min_{1\leq i\leq 2}\{\gamma(\mathcal{U}_{i})\}$. Pick $0<\epsilon<\frac{\gamma(\{\mathcal{U}_{i}\}_{i=1}^{2})}{2}$ and $N$ large enough such that $\Lambda^{{\bf a},t}_{N,\epsilon}(\omega,Z)\leq\frac{1}{2}$. Then by the definition of $\Lambda^{{\bf a},t}_{N,\epsilon}(\omega,Z)$,
there exists $\Gamma=\{B^{\bf a}_\omega(x_{j},n_{j},\epsilon)\}_{j}$ with $n_{j}\geq N$, and $\bigcup_{j}B^{\bf a}_\omega(x_{j},n_{j},\epsilon)\supseteq Z$ and $\sum_{j}\exp(-t n_{j})\leq 1$. Since $\epsilon<\frac{\gamma(\{\mathcal{U}_{i}\}_{i=1}^{2})}{2}$, there exists $\mathbf{U}^{{\bf a}}$ such that $ B^{\bf a}_{\omega}(x, n_{j},\epsilon)\subseteq X_{1}^{\omega}(\mathbf{U}^{{\bf a}})$. Then we have $\mathcal{M}^{{\bf a},t}_{N,\{\mathcal{U}_{i}\}_{i=1}^{2}}(\omega,Z)\leq\sum_{j}\exp(-t n_{j})\leq 1$. Letting $N\rightarrow\infty$, $\mathcal{M}^{{\bf a},t}_{\{\mathcal{U}_{i}\}_{i=1}^{2}}(\omega,Z)\leq 1$. This implies that $h^{\bf a}(\omega,f_{1},Z,\{\mathcal{U}_{i}\}_{i=1}^{2})\leq t$. By Proposition \ref{the existence of limit}, we have $h^{{\bf a}}(\omega,f_{1},Z)\geq\sup_{\{\mathcal{U}_{i}\}_{i=1}^{2}}h^{\bf a}(\omega,f_{1},Z,\{\mathcal{U}_{i}\}_{i=1}^{2})$.

For the other inequality, fix any $t>\sup_{\{\mathcal{U}_{i}\}_{i=1}^{2}}h^{\bf a}(\omega,f_{1},Z,\{\mathcal{U}_{i}\}_{i=1}^{2})$. For any $\delta>0$, we can choose a family of open covers $\{\mathcal{U}_{i}\}_{i=1}^{2}$ with $\text{diam}(\{\mathcal{U}_{i}\}_{i=1}^{2})<\delta$ and $N$ large enough such that $\mathcal{M}^{{\bf a},t}_{N,\{\mathcal{U}_{i}\}_{i=1}^{2}}(\omega,Z)<1$. Pick each $X_{1}^{\omega}(\mathbf{U}^{{\bf a}})\cap Z\neq\emptyset$ for $\{\mathcal{U}_{i}\}_{i=1}^{2}$, we choose
$x_{\mathbf{U}^{{\bf a}}}\in X_{1}^{\omega}(\mathbf{U}^{{\bf a}})\cap Z$ such that $X_{1}^{\omega}(\mathbf{U}^{{\bf a}})\subseteq B^{\bf a}_{\omega}(x_{\mathbf{U}^{{\bf a}}}, n_{j},\delta)$. Then we have $1\geq \mathcal{M}^{{\bf a},t}_{N,\{\mathcal{U}_{i}\}_{i=1}^{2}}(\omega,Z)\geq \Lambda^{{\bf a},t}_{N,\delta}(\omega,Z)$. Letting $N\rightarrow\infty$, we have $1\geq \Lambda^{{\bf a},t}_{\delta}(\omega,Z)$. This implies $h^{{\bf a}}(\omega,f_{1},Z,\delta)\leq t$. Letting $\delta\rightarrow0$, this finishes the proof.

\end{proof}

Next, we prove the Theorem \ref{the measurability of weighted topological entropy}.

\begin{proposition}
Let $f_{i},i=1,2$ be continuous bundle RDSs over $(\Omega,\mathcal{F},\mathbb{P},\vartheta)$. Assume that ${\bf a}=(a_{1},a_{2})\in\mathbb{R}^{2}$ with $a_{1}>0$ and $a_{2}\geq0$, $f_{2}$ is a factor of $f_{1}$ with a factor map $\Pi:\Omega\times X_{1}\to \Omega\times X_{2}$. Suppose that $\{\mathcal{U}_{i}\}_{i=1}^{2}$ are open covers where $\mathcal{U}_{i}$ is a finite open cover of $X_{i}$. Fix $s\in\mathbb{R}$ and $N\geq 1$. Then the function
\begin{align*}
(\omega,s)\mapsto\mathcal{M}^{{\bf a},s}_{N,\{\mathcal{U}_{i}\}_{i=1}^{2}}(\omega,X_{1})
\end{align*}
is measurable.
\end{proposition}

\begin{proof}
Modifying Kifer's arguments in \cite{Kifer2001On}, we give the proof process. Fix $s\in\mathbb{R}$ and $N\geq 1$. For $i=1,2$, let $\mathcal{U}_{i}=\{V^{(1)}_{i},\ldots,V^{(l_{i})}_{i}\}$.
Consider the elements of $\bigcup_{n\geq N}S_{n}(\{\mathcal{U}_{i}\}_{i=1}^{2})$, $\bigcup_{n\geq N}S_{n}(\{\mathcal{U}_{i}\}_{i=1}^{2})$ is a set consisting of strings
\begin{align*}
\mathbf{U}^{{\bf a}}((a_{0}^{1},\ldots,a_{\lceil(a_{1}+a_{2})n_{j}\rceil-1}^{2}))=&V_{1}^{(a_{0}^{1})} V_{1}^{(a_{1}^{1})} V_{1}^{(a_{2}^{1})}\cdots V_{1}^{(a_{\lceil a_{1}n_{j}\rceil-1}^{1})}\\
& V_{2}^{(a_{0}^{2})} V_{2}^{(a_{1}^{2})} V_{2}^{(a_{2}^{2})}\cdots V_{2}^{(a_{\lceil(a_{1}+a_{2})n_{j}\rceil-1}^{2})},
\end{align*}
where $V_{1}^{(a^{1}_{s})}\in \mathcal{U}_{1}$ for $0\leq s\leq\lceil a_1n_{j}\rceil-1$, $V_{2}^{(a^{2}_{s})}\in \mathcal{U}_{2}$ for $0\leq s\leq\lceil(a_1+a_{2})n_{j}\rceil-1$, and $m(\mathbf{U}^{{\bf a}}((a_{0}^{1},\ldots,a_{\lceil(a_{1}+a_{2})n_{j}\rceil-1}^{2})))=n_{j}\geq N$.

Fix an $n_{j}$-string $\mathbf{U}^{{\bf a}}((a_{0}^{1},\ldots,a_{\lceil(a_{1}+a_{2})n_{j}\rceil-1}^{2}))\in \bigcup_{n\geq N}S_{n}(\{\mathcal{U}_{i}\}_{i=1}^{2})$, define the set $A_{\Omega\times X_{1}}^{(a_{0}^{1},\ldots,a_{\lceil(a_{1}+a_{2})n_{j}\rceil-1}^{2})}$ as
\begin{align*}
\left\{(\omega,x):\omega\in \Omega,x\in X_{1}^{\omega}(\mathbf{U}^{{\bf a}}((a_{0}^{1},\ldots,a_{\lceil(a_{1}+a_{2})n_{j}\rceil-1}^{2})))\right\}.
\end{align*}
Since $\mathrm{U}_{i}^{a^{i}_{s}}=\Omega\times V_{i}^{(a^{i}_{s})}\in \mathcal{F}_{\Omega}\times \mathcal{B}_{X_{i}}$, where $V^{(a_{s}^{i})}_{i}\in \mathcal{U}_{i}$ for $i=1,2$,
\begin{align*}
A_{\Omega\times X_{1}}^{(a_{0}^{1},\ldots,a_{\lceil(a_{1}+a_{2})n_{j}\rceil-1}^{2})}&=\mathop{\bigcap}_{s=0}^{\lceil a_1n_{j}\rceil-1}(\Theta_{1}^{s})^{-1}\mathrm{U}_{1}^{a^{1}_{s}}\cap \mathop{\bigcap}_{s=0}^{\lceil (a_1+a_{2})n_{j}\rceil-1}(\Theta_{1}^{s})^{-1}\Pi^{-1}\mathrm{U}_{2}^{a^{2}_{s}}\\&\in \mathcal{F}_{\Omega}\times \mathcal{B}_{X_{1}}.
\end{align*}

{\bf Claim 1.} For any $n_{j}\in \mathbb{N}$. Define the function $\varphi_{(a_{0}^{1},\ldots,a_{\lceil(a_{1}+a_{2})n_{j}\rceil-1}^{2})}$ on $\Omega\times \mathbb{R}$ by
\begin{align*}
\varphi_{(a_{0}^{1},\ldots,a_{\lceil(a_{1}+a_{2})n_{j}\rceil-1}^{2})}=\exp\left(-s n_{j}\right).
\end{align*}
Then the function
\begin{align*}
(\omega,s)\mapsto \varphi_{(a_{0}^{1},\ldots,a_{\lceil(a_{1}+a_{2})n_{j}\rceil-1}^{2})}(\omega,s)
\end{align*}
is measurable.

\noindent{\bf Proof of Claim 1.} For each fixed element $\omega$, we have $s\mapsto\exp\left(-s n_{j}\right)$
is continuous. Fix $s\in\mathbb{R}$, we have $\omega\mapsto\exp\left(-s n_{j}\right)$ is measurable. By the [\cite{Crauel2002Random}, Lemma 1.1], this implies Claim 1.

Set $\overline{b}^{(j)}=(a_{0}^{1},\ldots,a_{\lceil(a_{1}+a_{2})n_{j}\rceil-1}^{2})$, where $j=1,\ldots,d$. Since $\Omega\times X_{1}\in \mathcal{F}_{\Omega}\times \mathcal{B}_{X_{1}}$, it follows that by [\cite{Castaing1977Convex}, Theorem II.30], we have
\begin{align*}
\Omega_{\{\overline{b}^{(1)},\overline{b}^{(2)},\ldots,\overline{b}^{(d)}\}}&=\left\{\omega:\mathop{\bigcup}_{j=1}^{d}X_{1}^{\omega}(\mathbf{U}^{{\bf a}}(\overline{b}^{(j)}))\supset X_{1}\right\}\\
&=\Omega\setminus \text{Pr}_{\Omega}\left(\left\{(\omega,x):(\omega,x)\in\left(\Omega\times X_{1}\setminus \mathop{\bigcup}_{j=1}^{d}A_{\Omega\times X_{1}}^{\overline{b}^{(j)}}\right)\right\}\right)\in \mathcal{F}_{\Omega}.
\end{align*}
Therefore for any $\beta>0$, we have
\begin{align*}
&\left\{(\omega,s):\mathcal{M}^{{\bf a},s}_{N,\{\mathcal{U}_{i}\}_{i=1}^{2}}(\omega,X_{1})<\beta\right\}\\=&\left\{(\omega,s):\inf_{\Gamma\subset\bigcup_{n\geq N}S_{n}(\{\mathcal{U}_{i}\}_{i=1}^{2}), \ \bigcup_{\mathbf{U}^{{\bf a}}\in \Gamma}X_{1}^{\omega}(\mathbf{U}^{{\bf a}})\supseteq X_{1}}\sum_{\mathbf{U}^{{\bf a}}\in \Gamma}\exp(-sm(\mathbf{U}^{{\bf a}}))<\beta\right\}\\
=&\bigcup_{\overline{b}^{(1)},\ldots,\overline{b}^{(d)}, \ d\in\mathbb{N}}\left\{(\omega,s): \sum_{i=1}^{d} \varphi_{\overline{b}^{(i)}}(\omega,s)<\beta\right\}\cap\left\{(\omega,s):\omega\in\Omega_{\{\overline{b}^{(1)},\ldots,\overline{b}^{(d)}\}},s\in\mathbb{R}\right\}.
\end{align*}
Notice that the number of interactions is at most countable in the last item. Then
\begin{align*}
\left\{(\omega,s):\mathcal{M}^{{\bf a},s}_{N,\{\mathcal{U}_{i}\}_{i=1}^{2}}(\omega,X_{1})<\beta\right\}\in \mathcal{F}_{\Omega}\times\mathcal{B}_{\mathbb{R}}.
\end{align*}
This implies this Proposition.
\end{proof}

By the above Proposition, the function $(\omega,s)\mapsto\mathcal{M}^{{\bf a},s}_{\{\mathcal{U}_{i}\}_{i=1}^{2}}(\omega,X_{1})$ is measurable in $\Omega$.

\begin{proof}[Proof of Theorem \ref{the measurability of weighted topological entropy}]
In order to prove $h^{{\bf a}}(\omega,f_{1},X_{1})$ is measurable in $\Omega$, it suffices to prove that
$$h^{\bf a}(\omega,f_{1},X_{1},\{\mathcal{U}_{i}\}_{i=1}^{2})=\inf \left\{s: \mathcal{M}^{{\bf a},s}_{\{\mathcal{U}_{i}\}_{i=1}^{2}}(\omega,X_{1})=0\right\}$$
is measurable in $\Omega$ by Proposition \ref{the existence of limit} and Proposition \ref{equality for weighted topological entropy}.

For each $\beta>0$, on the one hand, suppose that $h^{\bf a}(\omega,f_{1},X_{1},\{\mathcal{U}_{i}\}_{i=1}^{2})<\beta$, then by the definition of $h^{\bf a}(\omega,f_{1},X_{1},\{\mathcal{U}_{i}\}_{i=1}^{2})$, there exists a number $\alpha$ satisfying $$\mathcal{M}^{{\bf a},\alpha}_{\{\mathcal{U}_{i}\}_{i=1}^{2}}(\omega,X_{1})=0, $$ such that $\alpha<\beta$. On the another hand, if there exists a number $\alpha$ satisfying $$\mathcal{M}^{{\bf a},\alpha}_{\{\mathcal{U}_{i}\}_{i=1}^{2}}(\omega,X_{1})=0,$$ such that $\alpha<\beta$, we have $h^{\bf a}(\omega,f_{1},X_{1},\{\mathcal{U}_{i}\}_{i=1}^{2})<\beta$.

Consider the map $\text{\^{Pr}}_{\Omega}:\Omega\times\mathbb{R}\rightarrow \Omega$, $(\omega,\alpha)\mapsto \omega$, then by \cite{Crauel2002Random}, $\text{\^{Pr}}_{\Omega}$ is a measurable map. The set $\{(\omega,\alpha)\in\Omega\times\mathbb{R}:\alpha<\beta\}$ is measurable. By the above Proposition, the set $\{(\omega,\alpha)\in \Omega\times\mathbb{R}:\mathcal{M}^{{\bf a},\alpha}_{\{\mathcal{U}_{i}\}_{i=1}^{2}}(\omega,X_{1})=0\}$ is measurable.
We observe that
\begin{align*}
\{\omega&\in \Omega:h^{\bf a}(\omega,f_{1},X_{1},\{\mathcal{U}_{i}\}_{i=1}^{2})<\beta\}\\&=\{\omega: \ \exists \ \alpha\in \mathbb{R} \ \text{satisfies} \ \mathcal{M}^{{\bf a},\alpha}_{\{\mathcal{U}_{i}\}_{i=1}^{2}}(\omega,X_{1})=0,\alpha<\beta\}\\&=\text{\^{Pr}}_{\Omega}(\{(\omega,\alpha)\in\Omega\times \mathbb{R} :\alpha<\beta\}\cap\{(\omega,\alpha)\in\Omega\times\mathbb{R}:\mathcal{M}^{{\bf a},\alpha}_{\{\mathcal{U}_{i}\}_{i=1}^{2}}(\omega,X_{1})=0\}).
\end{align*}
It follows that $\{\omega\in \Omega:h^{\bf a}(\omega,f_{1},X_{1},\{\mathcal{U}_{i}\}_{i=1}^{2})<\beta\}$ is measurable in $\Omega$.
\end{proof}

By Theorem \ref{the measurability of weighted topological entropy}, we give the following definition.

\begin{definition}
The ${\bf a}$-weighted Bowen topological entropy of $X_{1}$ of $f_{1}$ is defined by
\begin{align*}
h^{{\bf a}}(f_{1},\Omega\times X_{1})= \int_{\Omega} h^{{\bf a}}(\omega,f_{1},X_{1})d\mathbb{P}(\omega).
\end{align*}
\end{definition}

\section{The proof of Theorem \ref{the first result}: lower bound}\label{The proof of Theorem the first result lower bound}

In this section, we prove the lower bound of Theorem \ref{the first result}. The following result is a weighted version of Brin-Katok local entropy formula:

\begin{theorem}\label{Brin-Katok local entropy formula}
Let $f_{i},i=1,2$ be continuous bundle RDSs over $(\Omega,\mathcal{F},\mathbb{P},\vartheta)$. Assume that ${\bf a}=(a_{1},a_{2})\in\mathbb{R}^{2}$ with $a_{1}>0$ and $a_{2}\geq0$, $f_{2}$ is a factor of $f_{1}$ with a factor map $\Pi:\Omega\times X_{1}\to \Omega\times X_{2}$. Let $\mu\in \mathcal{E}_{\mathbb{P}}^{1}(\Omega\times X_{1},f_{1})$, then for $\mu$-a.e. $(\omega,x)\in\Omega\times X_{1}$, we have
\begin{align*}
\lim\limits_{\epsilon\to0}\liminf\limits_{n\to\infty}-\dfrac{1}{n}\log\mu_\omega(B^{\bf a}_\omega(x,n,\epsilon))
&=\lim\limits_{\epsilon\to0}\limsup\limits_{n\to\infty}-\dfrac{1}{n}\log\mu_\omega(B^{\bf a}_\omega(x,n,\epsilon))\\&=a_{1}h_{\mu}^{(r)}(f_{1})+a_{2}h_{\mu\circ\Pi^{-1}}^{(r)}(f_{2})\\
&=h^{{\bf a}}_{\mu}(f_{1}).
\end{align*}
\end{theorem}

Next, the lower bound of Theorem \ref{the first result} can be obtained in the following theorem.

\begin{theorem}
Let $f_{i},i=1,2$ be continuous bundle RDSs over $(\Omega,\mathcal{F},\mathbb{P},\vartheta)$. Assume that ${\bf a}=(a_{1},a_{2})\in\mathbb{R}^{2}$ with $a_{1}>0$ and $a_{2}\geq0$, $f_{2}$ is a factor of $f_{1}$ with a factor map $\Pi:\Omega\times X_{1}\to \Omega\times X_{2}$. For any $\mu\in\mathcal{M}_{\mathbb{P}}^{1}(\Omega\times X_{1},f_{1})$, then we have
\begin{align*}
 h^{{\bf a}}(f_{1},\Omega\times X_{1})\geq a_{1}h_{\mu}^{(r)}(f_{1})+a_{2}h_{\mu\circ\Pi^{-1}}^{(r)}(f_{2}).
\end{align*}
\end{theorem}

\begin{proof}
We consider the case of ergodic measure $\mu\in \mathcal{E}_{\mathbb{P}}^{1}(\Omega\times X_{1},f_{1})$. In fact, by the ergodic decomposition thereom [\cite{Shu2009The}, Proposition 2.1], let $Q_{\mu}$ be its ergodic decomposition, $\mu=\int_{\mathcal{E}_{\mathbb{P}}^{1}(\Omega\times X_{1},f_{1})}\nu d Q_{\mu}(\nu)$, then we have
$\mu\circ\Pi^{-1}=\int_{\mathcal{E}_{\mathbb{P}}^{1}(\Omega\times X_{2},f_{2})}\nu d Q_{\mu}\circ\Pi^{-1}(\nu)$. Then we have
\begin{align*}
h^{{\bf a}}_{\mu}(f_{1})=&a_{1}h_{\mu}^{(r)}(f_{1})+a_{2}h_{\mu\circ\Pi^{-1}}^{(r)}(f_{2})
\\=&a_1\int_{\mathcal{E}_{\mathbb{P}}^{1}(\Omega\times X_{1},f_{1})}h_{\nu}^{(r)}(f_{1})d Q_{\mu}(\nu)+a_2\int_{\mathcal{E}_{\mathbb{P}}^{1}(\Omega\times X_{2},f_{2})}h_{\nu}^{(r)}(f_{2})d Q_{\mu}\circ\Pi^{-1}(\nu)\\
=&\int_{\mathcal{E}_{\mathbb{P}}^{1}(\Omega\times X_{1},f_{1})}a_{1}h_{\nu}^{(r)}(f_{1})+a_{2}h_{\nu\circ\Pi^{-1}}^{(r)}(f_{2})d Q_{\mu}(\nu)\\
=&\int_{\mathcal{E}_{\mathbb{P}}^{1}(\Omega\times X_{1},f_{1})}h^{{\bf a}}_{\nu}(f_{1})d Q_{\mu}(\nu).
\end{align*}
Hence it suffices to show that
\begin{align*}
h^{{\bf a}}(f_{1},\Omega\times X_{1})\geq \min\left\{\frac{1}{\delta},h^{{\bf a}}_{\mu}(f_{1})-\delta\right\}-\delta
\end{align*}
for any $\delta>0$.

Fix $\delta>0$ and an ergodic measure $\mu\in \mathcal{E}_{\mathbb{P}}^{1}(\Omega\times X_{1},f_{1})$. Let
\begin{align*}
H:=\min\left\{\frac{1}{\delta},h^{{\bf a}}_{\mu}(f_{1})-\delta\right\}.
\end{align*}
Let $G$ be a $\mu$ full measure set of elements satisfying
\begin{align*}
\lim\limits_{\epsilon\to0}\liminf\limits_{n\to\infty}-\dfrac{1}{n}\log\mu_\omega(B^{\bf a}_\omega(x,n,\epsilon))=h^{{\bf a}}_{\mu}(f_{1})
\end{align*}
for any $(\omega,x)\in G$ in Theorem \ref{Brin-Katok local entropy formula}. Fix $\omega\in \pi_{\Omega}(G)$. By Theorem \ref{Brin-Katok local entropy formula}, we can choose $\epsilon>0$ such that
\begin{align*}
\liminf\limits_{n\to\infty}-\dfrac{1}{n}\log\mu_\omega(B^{\bf a}_\omega(x,n,\epsilon))>H \ \text{for} \ \ \mu_{\omega}\text{-a.e.} \  \ x\in X_{1}.
\end{align*}
Then there exists a large $N\in\mathbb{N}$ and a Borel set $E_{N}\subset X_{1}$ with $\mu_{\omega}(E_{N})>\frac{1}{2}$ such that for any $x\in E_{N}$ and $n\geq N$,
\begin{align*}
\mu_{\omega}(B^{\bf a}_\omega(x,n,\epsilon))<\exp(-n H).
\end{align*}
Assume that $\Gamma=\{B^{\bf a}_\omega(x_{j},n_{j},\frac{\epsilon}{2})\}_{j}$ with $n_{j}\geq N$, and $\bigcup_{j}B^{\bf a}_\omega(x_{j},n_{j},\frac{\epsilon}{2})=X_{1}$. Set
\begin{align*}
\mathcal{I}:=\left\{j:B^{\bf a}_\omega(x_{j},n_{j},\frac{\epsilon}{2})\cap E_{N}\neq \emptyset\right\}.
\end{align*}
For $j\in\mathcal{I}$, pick $y_{j}\in B^{\bf a}_\omega(x_{j},n_{j},\frac{\epsilon}{2})\cap E_{N}$, then we have $B^{\bf a}_\omega(x_{j},n_{j},\frac{\epsilon}{2})\subseteq B^{\bf a}_\omega(y_{j},n_{j},\epsilon)$ and
\begin{align*}
\mu_{\omega}(B^{\bf a}_\omega(y_{j},n_{j},\epsilon))\leq \exp(-n_{j} H).
\end{align*}
Set $s:=H-\delta$. Then for any $j\in \mathcal{I}$, we have
\begin{align*}
\exp(-sn_{j})\geq \mu_{\omega}(B^{\bf a}_\omega(y_{j},n_{j},\epsilon))\exp(n_{j}(-s+H-\delta))=\mu_{\omega}(B^{\bf a}_\omega(y_{j},n_{j},\epsilon)).
\end{align*}
Therefore,
\begin{align*}
\sum_{j\in\mathcal{I}}\exp(-sn_{j})\geq \sum_{j\in\mathcal{I}}\mu_{\omega}(B^{\bf a}_\omega(y_{j},n_{j},\epsilon))\geq \mu_{\omega}\left(\bigcup_{j\in\mathcal{I}}B^{\bf a}_\omega(y_{j},n_{j},\epsilon)\right)\geq \mu_{\omega}(E_{N})\geq \frac{1}{2}.
\end{align*}
It follows that $\Lambda^{{\bf a},s}_{\epsilon}(\omega,X_{1})\geq \Lambda^{{\bf a},s}_{N,\epsilon}(\omega,X_{1})\geq \frac{1}{2}$. Then for any $\omega\in\pi_{\Omega}(G)$,
\begin{align}\label{inequation for lower bound}
h^{{\bf a}}(\omega,f_{1},X_{1})\geq \min\left\{\frac{1}{\delta},h^{{\bf a}}_{\mu}(f_{1})-\delta\right\}-\delta.
\end{align}
Since $\omega$ was chosen from the $\mathbb{P}$ full measure set $\pi_{\Omega}(G)$, integrating in \eqref{inequation for lower bound} against $\mathbb{P}$, this obtains
\begin{align*}
h^{{\bf a}}(f_{1},\Omega\times X_{1})\geq \min\left\{\frac{1}{\delta},h^{{\bf a}}_{\mu}(f_{1})-\delta\right\}-\delta,
\end{align*}
as desired.

\end{proof}

\section{The proof of Theorem \ref{the first result}: upper bound}\label{The proof of Theorem the first result upper bound}

In this section, we prove the upper bound of Theorem \ref{the first result}, this is,
\begin{align}\label{upper bound for weighted topological entropy}
 h^{{\bf a}}(f_{1},\Omega\times X_{1})\leq \sup\left\{a_{1}h_{\mu}^{(r)}(f_{1})+a_{2}h_{\mu\circ\Pi^{-1}}^{(r)}(f_{2}):\mu\in \mathcal{M}_{\mathbb{P}}^{1}(\Omega\times X_{1},f_{1})\right\}.
\end{align}
Before proving the above result, we first give some lemmas. We follows the definition of conditional entropy and that of Feng and Huang's proof \cite{Feng2008Variational}.

\begin{lemma}\label{the first estimation for measurable entropy function}
Let $f$ be a continuous bundle RDS over $(\Omega,\mathcal{F},\mathbb{P},\vartheta)$ and $\mu$ a probability measure on $\Omega\times X$. Let
$\mathcal{Q}=\{Q_{1},Q_{2},\ldots,Q_{M}\}$ be a Borel partition of $\Omega\times X$ with cardinality $M$. Write for brevity
\begin{align*}
h(n):=H_{\frac{1}{n}\sum_{i=0}^{n-1}\mu\circ \Theta^{-i}}(\mathcal{Q} \ | \ \pi_{\Omega}^{-1}\mathcal{F}), \  h(n,m):=H_{\frac{1}{m}\sum_{i=n}^{n+m-1}\mu\circ \Theta^{-i}}(\mathcal{Q} \ | \ \pi_{\Omega}^{-1}\mathcal{F})
\end{align*}
for $n,m\in\mathbb{N}$. Then we have
\begin{itemize}
  \item[(i)] $h(n)\leq \log M$ and $h(n,m)\leq \log M$ for any $n,m\in\mathbb{N}$.
  \item[(ii)] $|h(n+1)-h(n)|\leq \frac{1}{n+1}\log(3M^{2}(n+1))$ for any $n\in\mathbb{N}$.
  \item[(iii)] $|h(n,m)-\frac{n}{n+m}h(n)-\frac{m}{n+m}h(n,m)|\leq \log 2$ for any $n,m\in\mathbb{N}$.
\end{itemize}
\end{lemma}

\begin{proof}
(i) is obvious. Next we prove (ii). Let $\mathcal{W}$ be a finite Borel partition of $\Omega\times X$, $\nu_{1},\nu_{2}$ be $\Theta$-invariant probability measures on $\Omega\times X$ and $p\in[0,1]$, we have
\begin{align}\label{an inequation for conditional entropy}
\begin{split}
0&\leq H_{p\nu_{1}+(1-p)\nu_{2}}(\mathcal{W} \ | \ \pi_{\Omega}^{-1}\mathcal{F})-pH_{\nu_{1}}(\mathcal{W} \ | \ \pi_{\Omega}^{-1}\mathcal{F})-(1-p)H_{\nu_{2}}(\mathcal{W} \ | \ \pi_{\Omega}^{-1}\mathcal{F})\\
&\leq -(p\log p+(1-p)\log(1-p))\\
&\leq \log 2.
\end{split}
\end{align}
Let $n\in\mathbb{N}$. Applying \eqref{an inequation for conditional entropy} and (i), we have
\begin{align*}
&|h(n+1)-h(n)|\\=&|h(n+1)-\frac{n}{n+1}h(n)-\frac{1}{n+1}H_{\mu\circ\Theta^{-n}}(\mathcal{Q} \ | \ \pi_{\Omega}^{-1}\mathcal{F})-\frac{1}{n+1}h(n)\\&+\frac{1}{n+1}H_{\mu\circ\Theta^{-n}}(\mathcal{Q} \ | \ \pi_{\Omega}^{-1}\mathcal{F})|\\
\leq&|h(n+1)-\frac{n}{n+1}h(n)-\frac{1}{n+1}H_{\mu\circ\Theta^{-n}}(\mathcal{Q} \ | \ \pi_{\Omega}^{-1}\mathcal{F})|+\frac{2}{n+1}\log M\\
\leq&-\frac{n}{n+1}\log \frac{n}{n+1}-\frac{1}{n+1}\log \frac{1}{n+1}+\frac{2}{n+1}\log M\\
\leq&\frac{1}{n+1}\log (3M^{2}(n+1)),
\end{align*}
where we use the fact $(1+\frac{1}{n})^{n}<e<3$ in the last inequality.
This implies (ii). Since
\begin{align*}
\frac{1}{n+m}\sum_{i=0}^{n+m-1}\mu\circ \Theta^{-i}=\frac{n}{n+m}\left(\frac{1}{n}\sum_{i=0}^{n-1}\mu\circ \Theta^{-i}\right)+\frac{m}{n+m}\left(\frac{1}{m}\sum_{i=n}^{n+m-1}\mu\circ \Theta^{-i}\right)
\end{align*}
for $n,m\in\mathbb{N}$, (iii) follows from \eqref{an inequation for conditional entropy}.
\end{proof}

\begin{lemma}\label{the third inequation for conditional measure entropy}
Let $f$ be a continuous bundle RDS over $(\Omega,\mathcal{F},\mathbb{P},\vartheta)$ and $\mu$ a probability measure on $\Omega\times X$. For $\epsilon>0$ and $l,M\in\mathbb{N}$. Then the following statements hold.
\begin{itemize}
  \item[(1)] For any $n\in\mathbb{N}$,
  \begin{align*}
  &|H_{\frac{1}{n}\sum_{i=0}^{n-1}\mu\circ \Theta^{-i}}(\epsilon,M;l)-H_{\frac{1}{n+1}\sum_{i=0}^{n}\mu\circ \Theta^{-i}}(\epsilon,M;l)|\\
  \leq&\frac{1}{l(n+1)}\log(3M^{2l}(n+1)).
  \end{align*}
  \item[(2)] For any $n,m\in\mathbb{N}$,
  \begin{align*}
  &\frac{n}{n+m}H_{\frac{1}{n}\sum_{i=0}^{n-1}\mu\circ \Theta^{-i}}(\epsilon,M;l)+\frac{m}{n+m}H_{\frac{1}{n}\sum_{i=n}^{n+m-1}\mu\circ \Theta^{-i}}(\epsilon,M;l)\\
\leq& H_{\frac{1}{n+m}\sum_{i=0}^{n+m-1}\mu\circ \Theta^{-i}}(\epsilon,M;l)+\frac{\log 2}{l}.
  \end{align*}
\end{itemize}
\end{lemma}

\begin{proof}
The statements follow from the definition of $H_{\bullet}(\epsilon,M;l)$ and Lemma \ref{the first estimation for measurable entropy function}.
\end{proof}

The following result mainly follows the definition of conditional entropy and [\cite{Cao2008The}, Lemma 2.4].

\begin{lemma}\label{the second inequation for conditional measure entropy}
Let $\mu$ be a probability measure on $\Omega\times X$ and $M\in\mathbb{N}$. Suppose $\xi=\{A_{1},\ldots,A_{j}\}$ is a Borel partition of $X$ with $j\leq M$. Then for any position integers $n,l$ with $n\geq 2l$, we have
\begin{align*}
\frac{1}{n}H_{\nu}(\bigvee_{i=0}^{n-1}\Theta^{-i}\mathcal{Q} \ | \ \pi_{\Omega}^{-1}\mathcal{F})\leq \frac{1}{l}H_{\nu_{n}}(\bigvee_{i=0}^{l-1}\Theta^{-i}\mathcal{Q} \ | \ \pi_{\Omega}^{-1}\mathcal{F})+\frac{2l}{n}\log M,
\end{align*}
where $\nu_{n}=\frac{1}{n}\sum_{i=0}^{n-1}\nu\circ \Theta^{-i}$ and $\mathcal{Q}=\{\Omega\times A:A\in\xi\}$.
\end{lemma}
\begin{proof}
For $0\leq s\leq l-1$, let $t_{s}$ denote the largest integer $t$ so that $tl+s\leq n$. Then we
have
\begin{align*}
\bigvee_{i=0}^{n-1}\Theta^{-i}\mathcal{Q}=\bigvee_{r=0}^{t_{s}-1}\Theta^{-rl-s}\bigvee_{i=0}^{l-1}\Theta^{-i}\mathcal{Q}\vee\bigvee_{m\in C_{s}}\Theta^{-m}\mathcal{Q}
\end{align*}
for $0\leq s\leq l-1$, where $C_{s}$ is a subset of $\{0,\ldots,n-1\}$ with cardinality at most $2l$. Then
\begin{align*}
H_{\nu}(\bigvee_{i=0}^{n-1}\Theta^{-i}\mathcal{Q} \ | \ \pi_{\Omega}^{-1}\mathcal{F})\leq\sum_{r=0}^{t_{s}-1}H_{\nu}(\Theta^{-rl-s}\bigvee_{i=0}^{l-1}\Theta^{-i}\mathcal{Q} \ | \ \pi_{\Omega}^{-1}\mathcal{F})+2l\log M.
\end{align*}
Summing this over $s$ from $0$ to $l-1$ gives
\begin{align*}
lH_{\nu}(\bigvee_{i=0}^{n-1}\Theta^{-i}\mathcal{Q} \ | \ \pi_{\Omega}^{-1}\mathcal{F})&\leq\sum_{s=0}^{l-1}\sum_{r=0}^{t_{s}-1}H_{\nu}(\Theta^{-rl-s}\bigvee_{i=0}^{l-1}\Theta^{-i}\mathcal{Q} \ | \ \pi_{\Omega}^{-1}\mathcal{F})+2l^{2}\log M\\
&=\sum_{p=0}^{n-l}H_{\nu}(\Theta^{-p}\bigvee_{i=0}^{l-1}\Theta^{-i}\mathcal{Q} \ | \ \pi_{\Omega}^{-1}\mathcal{F})+2l^{2}\log M\\
&\leq\sum_{p=0}^{n-1}H_{\nu}(\Theta^{-p}\bigvee_{i=0}^{l-1}\Theta^{-i}\mathcal{Q} \ | \ \pi_{\Omega}^{-1}\mathcal{F})+2l^{2}\log M\\
&=\sum_{p=0}^{n-1}H_{\nu\circ\Theta^{-p}}(\bigvee_{i=0}^{l-1}\Theta^{-i}\mathcal{Q} \ | \ \pi_{\Omega}^{-1}\mathcal{F})+2l^{2}\log M\\
&=nH_{\nu_{n}}(\bigvee_{i=0}^{l-1}\Theta^{-i}\mathcal{Q} \ | \ \pi_{\Omega}^{-1}\mathcal{F})+2l^{2}\log M.
\end{align*}
This implies the desired inequality.
\end{proof}

Next, we give the proof of upper bound of Theorem \ref{the first result}:

\begin{proof}[Proof of \eqref{upper bound for weighted topological entropy}]
We divide the proof into two cases:

{\bf Case 1.} If $h^{{\bf a}}(f_{1},\Omega\times X_{1})=0$, this directly obtains \eqref{upper bound for weighted topological entropy} since ${\bf a}=(a_{1},a_{2})\in\mathbb{R}^{2}$ with $a_{1}>0$ and $a_{2}\geq0$.

{\bf Case 2.} If $h^{{\bf a}}(f_{1},\Omega\times X_{1})>0$, set
\begin{align*}
\Omega=\left\{\omega:h^{{\bf a}}(\omega,f_{1},X_{1})>0\right\}\cup\left\{\omega:h^{{\bf a}}(\omega,f_{1},X_{1})=0\right\}
\end{align*}
where $F=\{\omega:h^{{\bf a}}(\omega,f_{1},X_{1})>0\}$, by Proposition \ref{the measurability of weighted topological entropy}, $F\subseteq \Omega$ is a positive measurable subset such that for any $\omega\in F$, $h^{{\bf a}}(\omega,f_{1},X_{1})>0$. Then
\begin{align*}
h^{{\bf a}}(f_{1},\Omega\times X_{1})&= \int_{F} h^{{\bf a}}(\omega,f_{1},X_{1})d\mathbb{P}(\omega)+\int_{\{\omega:h^{{\bf a}}(\omega,f_{1},X_{1})=0\}} h^{{\bf a}}(\omega,f_{1},X_{1})d\mathbb{P}(\omega)\\
&=\int_{F} h^{{\bf a}}(\omega,f_{1},X_{1})d\mathbb{P}(\omega).
\end{align*}

Let $G_{\mathbb{P}}$ be the set of all generic points of $(\Omega,\mathcal{F},\mathbb{P},\vartheta)$, this is, $\omega\in G_{\mathbb{P}}$ if and
only if $\mathbb{P}=\lim_{m\rightarrow\infty}\frac{1}{m}\sum_{j=0}^{m-1}\delta_{\vartheta^{j}\omega}$ in the weak* topology. By the Birkhoff
pointwise ergodic theorem, we have $\mathbb{P}(G_{\mathbb{P}})=1$ as $(\Omega,\mathcal{F},\mathbb{P},\vartheta)$ is ergodic. Then we have $\mathbb{P}(F\cap G_{\mathbb{P}})>0$, and
\begin{align*}
h^{{\bf a}}(f_{1},\Omega\times X_{1})=\int_{F\cap G_{\mathbb{P}}} h^{{\bf a}}(\omega,f_{1},X_{1})d\mathbb{P}(\omega).
\end{align*}
In order to prove \eqref{upper bound for weighted topological entropy}, it suffices to prove that for any $\omega\in F\cap G_{\mathbb{P}}$,
\begin{align*}
 h^{{\bf a}}(\omega,f_{1},\Omega\times X_{1})\leq \sup\left\{a_{1}h_{\mu}^{(r)}(f_{1})+a_{2}h_{\mu\circ\Pi^{-1}}^{(r)}(f_{2}):\mu\in \mathcal{M}_{\mathbb{P}}^{1}(\Omega\times X_{1},f_{1})\right\},
\end{align*}
that is, there exists $\mu\in \mathcal{M}_{\mathbb{P}}^{1}(\Omega\times X_{1},f_{1})$ such that
\begin{align*}
h^{{\bf a}}(\omega,f_{1},\Omega\times X_{1})\leq a_{1}h_{\mu}^{(r)}(f_{1})+a_{2}h_{\mu\circ\Pi^{-1}}^{(r)}(f_{2}).
\end{align*}

Fix $\omega_{0}\in F\cap G_{\mathbb{P}}$. Then $h^{{\bf a}}(\omega_{0},f_{1},X_{1})>0$. Fix $0<s<h^{{\bf a}}(\omega_{0},f_{1},X_{1})$. By Lemma \ref{a dynamical Frostman lemma of RDS}, there exists a probability measure $\nu_{\omega_{0}}$ of $X_{1}$ with $\nu_{\omega_{0}}(X_{1})=1$, $\epsilon>0$, and $N\in\mathbb{N}$ such that
\begin{align}\label{the inequality for lemma 2.7}
\nu_{\omega_{0}}(B^{\bf a}_{\omega_{0}}(x,n,\epsilon))\leq \exp(-sn)
\end{align}
for any $n\geq N$ and $x\in X_{1}$. By the continuity, there exists $\tau\in(0,\epsilon)$ such that if $x,y\in X_{1}$ satisfy $d_{1}(x,y)<\tau$, then
\begin{align*}
d_{2}(\pi(x),\pi(y))<\epsilon.
\end{align*}
Take $M_{0}\in\mathbb{N}$ with $\mathcal{P}_{X_{i}}(\tau,M_{0})\neq\emptyset$ for $i=1,2$, where $\mathcal{P}_{X_{i}}(\tau,M_{0})$ is defined as in \eqref{definition of finite Borel partition}. Now fix $M\in\mathbb{N}$ with $M\geq M_{0}$. Let $\alpha_{i}\in\mathcal{P}_{X_{i}}(\tau,M)$ for $i=1,2$. Set $\beta_{1}=\alpha_{1}$, $\beta_{2}=\pi^{-1}\alpha_{2}$. Then for any $n\in\mathbb{N}$ and $x\in X_{1}$, we have
\begin{align*}
\bigvee_{j=0}^{\lceil a_{1}n\rceil-1}(f_{1,\omega_{0}}^{j})^{-1}\beta_{1}\vee\bigvee_{j=\lceil a_{1}n\rceil}^{\lceil (a_{1}+a_{2})n\rceil-1}(f_{1,\omega_{0}}^{j})^{-1}\beta_{2}(x)\subseteq B^{\bf a}_{\omega_{0}}(x,n,\epsilon).
\end{align*}
Assume that $n\geq N$. By \eqref{the inequality for lemma 2.7}, we have
\begin{align*}
\nu_{\omega_{0}}(\bigvee_{j=0}^{\lceil a_{1}n\rceil-1}(f_{1,\omega_{0}}^{j})^{-1}\beta_{1}\vee\bigvee_{j=\lceil a_{1}n\rceil}^{\lceil (a_{1}+a_{2})n\rceil-1}(f_{1,\omega_{0}}^{j})^{-1}\beta_{2}(x))\leq\exp(-sn)
\end{align*}
for any $x\in X_{1}$. Then we construct a measurable $\nu$ with support $\{\omega_{0}\}\times X_{1}$ such that
\begin{align*}
d\nu(\omega,x)=d\nu_{\omega_{0}}(x)d\delta_{\omega_{0}}(\omega).
\end{align*}
Then we have
\begin{align*}
&H_{\nu}(\bigvee_{j=0}^{\lceil a_{1}n\rceil-1}(\Theta_{1}^{j})^{-1}\mathcal{Q}_{1}\vee\bigvee_{j=\lceil a_{1}n\rceil}^{\lceil (a_{1}+a_{2})n\rceil-1}(\Theta_{1}^{j})^{-1}\mathcal{Q}_{2} \ | \ \pi_{\Omega}^{-1}\mathcal{F})\\=&-\int\log\nu_{\omega_{0}}(\bigvee_{j=0}^{\lceil a_{1}n\rceil-1}(f_{1,\omega_{0}}^{j})^{-1}\beta_{1}\vee\bigvee_{j=\lceil a_{1}n\rceil}^{\lceil (a_{1}+a_{2})n\rceil-1}(f_{1,\omega_{0}}^{j})^{-1}\beta_{2}(x))d\nu_{\omega_{0}}(x)\\
\geq &sn,
\end{align*}
where $\mathcal{Q}_{1}=\{\Omega\times A:A\in \beta_{1}\}$ and $\mathcal{Q}_{2}=\{\Omega\times A:A\in \beta_{2}\}$. Hence
\begin{align}\label{upper estmation 1}
H_{\nu}(\bigvee_{j=0}^{\lceil a_{1}n\rceil-1}(\Theta_{1}^{j})^{-1}\mathcal{Q}_{1} \ | \ \pi_{\Omega}^{-1}\mathcal{F})+H_{\nu}(\bigvee_{j=\lceil a_{1}n\rceil}^{\lceil (a_{1}+a_{2})n\rceil-1}(\Theta_{1}^{j})^{-1}\mathcal{Q}_{2} \ | \ \pi_{\Omega}^{-1}\mathcal{F})\geq sn.
\end{align}
Now fix $l\in\mathbb{N}$. By Lemma \ref{the second inequation for conditional measure entropy}, for sufficiently large $n$, the left-hand side of \eqref{upper estmation 1} is bounded from above by
\begin{align*}
&\frac{\lceil a_{1}n\rceil}{l}H_{\varpi_{1,n}}(\bigvee_{j=0}^{l-1}(\Theta_{1}^{j})^{-1}\mathcal{Q}_{1} \ | \ \pi_{\Omega}^{-1}\mathcal{F})\\&+\frac{\lceil (a_{1}
+a_{2})n\rceil-\lceil a_{1}n\rceil}{l}H_{\varpi_{2,n}}(\bigvee_{j=0}^{l-1}(\Theta_{1}^{j})^{-1}\mathcal{Q}_{2} \ | \ \pi_{\Omega}^{-1}\mathcal{F})+4l\log M,
\end{align*}
where
\begin{align*}
\varpi_{1,n}=\frac{\sum_{j=0}^{\lceil a_{1}n\rceil-1}\nu\circ(\Theta_{1})^{-j}}{\lceil a_{1}n\rceil} \ \ \ \text{and} \ \ \ \varpi_{2,n}=\frac{\sum_{j=\lceil a_{1}n\rceil}^{\lceil (a_{1}+a_{2})n\rceil-1}\nu\circ(\Theta_{1})^{-j}}{\lceil (a_{1}+a_{2})n\rceil-\lceil a_{1}n\rceil}.
\end{align*}
Hence by \eqref{upper estmation 1} and the definition of $H_{\bullet}(\tau,M;l)$, we have
\begin{align}\label{upper estmation 2}
\lceil a_{1}n\rceil H_{\varpi_{1,n}}(\tau,M;l)+(\lceil (a_{1}
+a_{2})n\rceil-\lceil a_{1}n\rceil)H_{\varpi_{2,n}\circ\Pi^{-1}}(\tau,M;l)
\geq sn -4l\log M.
\end{align}
Define $\nu_{m}=\frac{\sum_{j=0}^{m-1}\nu\circ(\Theta_{1})^{-j}}{m}$ for $m\in\mathbb{N}$. Then
\begin{align*}
\nu_{m}\circ\Pi^{-1}=\frac{\sum_{j=0}^{m-1}(\nu\circ\Pi^{-1})\circ(\Theta_{2})^{-j}}{m},
\end{align*}
\begin{align*}
\varpi_{2,n}\circ\Pi^{-1}=\frac{\sum_{j=\lceil a_{1}n\rceil}^{\lceil (a_{1}+a_{2})n\rceil-1}(\nu\circ\Pi^{-1})\circ(\Theta_{2})^{-j}}{\lceil (a_{1}+a_{2})n\rceil-\lceil a_{1}n\rceil}, \ \ \nu_{\lceil a_{1}n\rceil}=\varpi_{1,n},
\end{align*}
and
\begin{align*}
\nu_{\lceil (a_{1}+a_{2})n\rceil}\circ\Pi^{-1}=&\frac{\lceil a_{1}n\rceil}{\lceil (a_{1}+a_{2})n\rceil}\nu_{\lceil a_{1}n\rceil}\circ\Pi^{-1}+\frac{\lceil (a_{1}+a_{2})n\rceil-\lceil a_{1}n\rceil}{\lceil (a_{1}+ a_{2})n\rceil}\varpi_{2,n}\circ\Pi^{-1}.
\end{align*}
By Lemma \ref{the third inequation for conditional measure entropy} (2), we have
\begin{align*}
&\frac{\lceil a_{1}n\rceil}{\lceil (a_{1}+a_{2})n\rceil}H_{\nu_{\lceil a_{1}n\rceil}\circ\Pi^{-1}}(\tau,M;l)+\frac{\lceil (a_{1}+a_{2})n\rceil-\lceil a_{1}n\rceil}{\lceil (a_{1}+ a_{2})n\rceil}H_{\varpi_{2,n}\circ\Pi^{-1}}(\tau,M;l)\\
\leq&H_{\nu_{\lceil (a_{1}+a_{2})n\rceil}\circ\Pi^{-1}}(\tau,M;l)+\frac{\log 2}{l}.
\end{align*}
That is,
\begin{align*}
&\lceil (a_{1}+a_{2})n\rceil H_{\nu_{\lceil (a_{1}+a_{2})n\rceil}\circ\Pi^{-1}}(\tau,M;l)-\lceil a_{1}n\rceil H_{\nu_{\lceil a_{1}n\rceil}\circ\Pi^{-1}}(\tau,M;l)\\
\geq& (\lceil (a_{1}+a_{2})n\rceil-\lceil a_{1}n\rceil)H_{\varpi_{2,n}\circ\Pi^{-1}}(\tau,M;l)-\frac{\lceil (a_{1}+a_{2})n\rceil}{l}\log 2.
\end{align*}
Combining the above inequality with \eqref{upper estmation 2}, we have
\begin{align}\label{upper estmation 3}
\begin{split}
&\lceil a_{1}n\rceil H_{\nu_{\lceil a_{1}n\rceil}}(\tau,M;l)+\lceil (a_{1}+a_{2})n\rceil H_{\nu_{\lceil (a_{1}+a_{2})n\rceil}\circ\Pi^{-1}}(\tau,M;l)\\&-\lceil a_{1}n\rceil H_{\nu_{\lceil a_{1}n\rceil}\circ\Pi^{-1}}(\tau,M;l)\\
\geq&\lceil a_{1}n\rceil H_{\varpi_{1,n}}(\tau,M;l)-\frac{\lceil a_{1}n\rceil}{l}+(\lceil (a_{1}+a_{2})n\rceil-\lceil a_{1}n\rceil)H_{\varpi_{2,n}\circ\Pi^{-1}}(\tau,M;l)\\&-\frac{\lceil (a_{1}+a_{2})n\rceil}{l}\log 2\\
\geq&sn -4l\log M-\frac{2\lceil (a_{1}+a_{2})n\rceil}{l}\log 2.
\end{split}
\end{align}
By Lemma \ref{the third inequation for conditional measure entropy} (1), we have
\begin{align}\label{inequation for estmation}
\begin{split}
|H_{\nu_{n}\circ\Pi^{-1}}(\tau,M;l)-H_{\nu_{n+1}\circ\Pi^{-1}}(\tau,M;l)|&\leq\frac{1}{l(n+1)}\log(3M^{2l}(n+1)).
\end{split}
\end{align}
Set
\begin{align*}
\gamma(n)=&\lceil (a_{1}+a_{2})n\rceil (H_{\nu_{\lceil (a_{1}+a_{2})n\rceil}\circ\Pi^{-1}}(\tau,M;l)- H_{\nu_{\lceil a_{1}n\rceil}\circ\Pi^{-1}}(\tau,M;l)).
\end{align*}
Then we have
\begin{align*}
&\lceil a_{1}n\rceil H_{\nu_{\lceil a_{1}n\rceil}}(\tau,M;l)+\lceil (a_{1}+a_{2})n\rceil H_{\nu_{\lceil (a_{1}+a_{2})n\rceil}\circ\Pi^{-1}}(\tau,M;l)\\&-\lceil a_{1}n\rceil H_{\nu_{\lceil a_{1}n\rceil}\circ\Pi^{-1}}(\tau,M;l)\\
=&\gamma(n)+\lceil a_{1}n\rceil H_{\nu_{\lceil a_{1}n\rceil}}(\tau,M;l)+(\lceil (a_{1}+a_{2})n\rceil-\lceil a_{1}n\rceil)H_{\nu_{\lceil a_{1}n\rceil}\circ\Pi^{-1}}(\tau,M;l).
\end{align*}
By \eqref{upper estmation 3}, we have
\begin{align}\label{upper estmation 4}
\begin{split}
&\frac{\lceil a_{1}n\rceil}{n} H_{\nu_{\lceil a_{1}n\rceil}}(\tau,M;l)+\frac{\lceil (a_{1}+a_{2})n\rceil-\lceil a_{1}n\rceil}{n}H_{\nu_{\lceil a_{1}n\rceil}\circ\Pi^{-1}}(\tau,M;l)\\
\geq& -\frac{\gamma(n)}{n}+s -\frac{4l\log M}{n}-\frac{2\lceil (a_{1}+a_{2})n\rceil}{nl}\log 2.
\end{split}
\end{align}
Define
\begin{align*}
\varpi(n)=-(a_{1}+a_{2})(H_{\nu_{\lceil (a_{1}+a_{2})n\rceil}\circ\Pi^{-1}}(\tau,M;l)- H_{\nu_{\lceil a_{1}n\rceil}\circ\Pi^{-1}}(\tau,M;l)).
\end{align*}
Then we have $\limsup_{n\rightarrow\infty}\varpi(n)\geq 0$ by applying [\cite{Feng2008Variational}, Lemma 5.4], if which we take $p=2$,
\begin{align*}
u_{j}(n)=\left\{
\begin{aligned}
a_{1}H_{\nu_{n}\circ\Pi^{-1}}(\tau,M;l) \ \ \ \ \ \ \ \ \ \ \ \ \ \ \ \ &\text{for} \ j=1,\\
-(a_{1}+a_{2})H_{\nu_{n}\circ\Pi^{-1}}(\tau,M;l) \ \ &\text{for} \ j=2,
\end{aligned}
\right.
\end{align*}
and
\begin{align*}
c_{j}(n)=\left\{
\begin{aligned}
a_{1} \ \ \ \ \ \ \ \ \ \ \ \ &\text{for} \ j=1,\\
a_{1}+a_{2} \ \ \  &\text{for} \ j=2,
\end{aligned}
\right.
\end{align*}
and $r_{j}=a_{1}$ for any $j=1,2$. The condition $\lim_{n\rightarrow\infty}|u_{j}(n+1)-u_{j}(n)|=0$ fulfils from \eqref{inequation for estmation}. Then we have
\begin{align*}
\limsup_{n\rightarrow\infty}\frac{-\gamma(n)}{n}=\limsup_{n\rightarrow\infty}\varpi(n)\geq 0.
\end{align*}
Letting $n\rightarrow\infty$ in \eqref{upper estmation 4} and taking the upper limit, we obtain
\begin{align}\label{upper estmation 5}
\limsup_{n\rightarrow\infty}(a_{1}H_{\nu_{\lceil a_{1}n\rceil}}(\tau,M;l)+a_{2}H_{\nu_{\lceil a_{1}n\rceil}\circ\Pi^{-1}}(\tau,M;l)) \geq s-\frac{2(a_{1}+a_{2})}{l}\log 2.
\end{align}
Take a subsequence $(n_{j})$ of natural numbers so that the left-hand side of \eqref{upper estmation 5} equals
\begin{align*}
\limsup_{j\rightarrow\infty}(a_{1}H_{\nu_{\lceil a_{1}n_{j}\rceil}}(\tau,M;l)+a_{2}H_{\nu_{\lceil a_{1}n_{j}\rceil}\circ\Pi^{-1}}(\tau,M;l)),
\end{align*}
and moreover, $\nu_{\lceil a_{1}n_{j}\rceil}$ converges to an element $\lambda\in  \mathcal{M}_{\mathbb{P}}^{1}(\Omega\times X_{1},f_{1})$ in the weak* topology by [\cite{Kifer2001On}, Lemma 2.1]. Indeed, for any $\varphi\in C(\Omega)$, we have
\begin{align*}
\int_{\Omega}\varphi d\lambda\circ\pi_{\Omega}^{-1}&=\int_{\Omega\times X_{1}}\varphi\circ\pi_{\Omega}d\lambda\\
&=\lim_{j\rightarrow\infty}\frac{1}{\lceil a_{1}n_{j}\rceil}\sum_{s=0}^{\lceil a_{1}n_{j}\rceil-1}\int_{\Omega\times X_{1}}\varphi\circ\pi_{\Omega}d\nu\circ(\Theta_{1})^{-s}\\
&=\lim_{j\rightarrow\infty}\frac{1}{\lceil a_{1}n_{j}\rceil}\sum_{s=0}^{\lceil a_{1}n_{j}\rceil-1}\int_{\Omega\times X_{1}}\varphi\circ\pi_{\Omega}\circ(\Theta_{1})^{s}d\nu\\
&=\lim_{j\rightarrow\infty}\frac{1}{\lceil a_{1}n_{j}\rceil}\sum_{s=0}^{\lceil a_{1}n_{j}\rceil-1}\int_{\Omega}\varphi(\vartheta^{s}\omega)d\delta_{\omega_{0}}\\
&=\lim_{j\rightarrow\infty}\frac{1}{\lceil a_{1}n_{j}\rceil}\sum_{s=0}^{\lceil a_{1}n_{j}\rceil-1}\varphi(\vartheta^{s}\omega_{0})\\
&=\int_{\Omega}\varphi d\mathbb{P}(\omega).
\end{align*}
Since the map $H_{\bullet}(\tau,M;l)$ is upper semi-continuous on $\mathcal{P}_{\mathbb{P}}(\Omega\times X)$ (see Lemma \ref{upper semi-continuous function}), we have
\begin{align}\label{upper estmation 6}
a_{1}H_{\lambda}(\tau,M;l)+a_{2}H_{\lambda\circ\Pi^{-1}}(\tau,M;l)\geq s-\frac{2(a_{1}+a_{2})}{l}\log 2.
\end{align}
Define
\begin{align*}
\mathcal{E}:=\left\{(M,l,\delta):M,l\in\mathbb{N},\delta>0 \ \text{with} \ M \geq M_{0}, \ l\geq \frac{2(a_{1}+a_{2})}{l}\log 2\right\},
\end{align*}
and
\begin{align*}
\Omega_{M,l,\delta}:=\left\{\eta\in \mathcal{M}_{\mathbb{P}}^{1}(\Omega\times X_{1},f_{1}):a_{1}H_{\eta}(\tau,M;l)+a_{2}H_{\eta\circ\Pi^{-1}}(\tau,M;l)\geq s-\delta\right\}.
\end{align*}
By \eqref{upper estmation 6}, $\Omega_{M,l,\delta}$ is a non-empty compact set whenever $(M,l,\delta)\in\mathcal{E}$. However,
\begin{align*}
\Omega_{M_{1},l_{1},\delta_{1}}\cap\Omega_{M_{2},l_{2},\delta_{2}}\supseteq \Omega_{M_{1}+M_{1},l_{1}l_{2},\min\{\delta_{1},\delta_{2}\}}
\end{align*}
for any $(M_{1},l_{1},\delta_{1}), (M_{2},l_{2},\delta_{2})\in \mathcal{E}$. It follows that
\begin{align*}
\bigcap_{(M,l,\delta)\in\mathcal{E}}\Omega_{M,l,\delta}\neq\emptyset.
\end{align*}
Take $\mu_{s}\in \bigcap_{(M,l,\delta)\in\mathcal{E}}\Omega_{M,l,\delta}$. Then
\begin{align*}
a_{1}h_{\mu_{s}}^{(r)}(f_{1},\epsilon)+a_{2}h_{\mu_{s}\circ\Pi^{-1}}^{(r)}(f_{2},\epsilon)\geq s.
\end{align*}
Since the map $\theta\in\mathcal{M}_{\mathbb{P}}^{1}(\Omega\times X_{1},f_{1})\mapsto a_{1}h_{\theta}^{(r)}(f_{1},\epsilon)+a_{2}h_{\theta\circ\Pi^{-1}}^{(r)}(f_{2},\epsilon)$ is upper semi-continuous (see Lemma \ref{upper semi-continuous function}), we can find $\mu\in \mathcal{M}_{\mathbb{P}}^{1}(\Omega\times X,f_{1})$ such that
\begin{align*}
a_{1}h_{\mu}^{(r)}(f_{1},\epsilon)+a_{2}h_{\mu\circ\Pi^{-1}}^{(r)}(f_{2},\epsilon)\geq s.
\end{align*}
Since $a_{1}h_{\mu}^{(r)}(f_{1})+a_{2}h_{\mu\circ\Pi^{-1}}^{(r)}(f_{2})\geq a_{1}h_{\mu}^{(r)}(f_{1},\epsilon)+a_{2}h_{\mu\circ\Pi^{-1}}^{(r)}(f_{2},\epsilon)$, we have
\begin{align*}
a_{1}h_{\mu}^{(r)}(f_{1})+a_{2}h_{\mu\circ\Pi^{-1}}^{(r)}(f_{2})\geq s.
\end{align*}
Letting $s\nearrow h^{{\bf a}}(\omega_{0},f_{1},X_{1})$, this completes the proof the upper bound in Theorem \ref{the first result}.

\end{proof}

\section{Appendix A: weighted Shannon-McMillan-Breiman Theorem}\label{Appendix A: weighted Shannon-McMillan-Breiman Theorem}

In this Appendix A, we give the proof of weighted Shannon-McMillan-Breiman Theorem of
RDSs. First we need the following results:

\begin{lemma}\label{lemmaduo}
Let $f$ be a continuous bundle RDS over $(\Omega,\mathcal F,\mathbb P,\vartheta)$, $\mu\in\mathcal{M}_{\mathbb{P}}^{1}(\Omega\times X,f)$, $\alpha,\beta$ be two finite measurable partitions of $\Omega\times X$ and $\mathcal A$ be a sub-$\sigma$-algebra of $\mathcal B$. Let $I_{\mu_\omega}(\cdot|\cdot)$ be the conditional information of $\mu_\omega.$ Then for $\mu$-a.e. $(\omega,x)\in\Omega\times X$, we have

\begin{enumerate}
  \item[(1)] $
 I_{\mu_\omega}((f_\omega)^{-1}\alpha_{\vartheta\omega}|(f_{\omega})^{-1} \beta_{\vartheta\omega}) (x)=I_{\mu_{\vartheta\omega}}(\alpha_{\vartheta\omega}| \beta_{\vartheta\omega})(f_{\omega}(x));
$
  \item[(2)] $ I_{\mu_\omega}(\alpha_\omega\vee\beta_\omega|\mathcal A)(x)=I_{\mu_\omega}(\alpha_\omega |\mathcal A)(x)+ I_{\mu_\omega}(\beta_\omega|\alpha_\omega \vee\mathcal A)(x);$
  \item[(3)] if $\mathcal A_1\subset \mathcal A_2\subset \cdots$ is an increasing sequence of $\sigma$-algebra of $\mathcal B$ with $\mathcal A_n\to \mathcal A,$ then
\begin{align*}
\lim\limits_{n\to\infty} I_{\mu_\omega}(\alpha_\omega| \mathcal A_n)(x)=I_{\mu_\omega}(\alpha_\omega| \mathcal A )(x).
\end{align*}
\end{enumerate}
\end{lemma}
\begin{proof}
The proof details can be seen in \cite{Parry1981Topics}.
\end{proof}

The following Lemma was proved in [\cite{Feng2008Variational}, Lemma A.5].

\begin{lemma}\label{lemmaergodic}
Let $F_n\in L^1(\Omega\times X,\mathcal F\times \mathcal{B},\mu)$ be a sequence that converges  almost everywhere and in $L^1$ to  $F\in L^1(\Omega\times X,\mathcal F\times \mathcal{B},\mu)$ and $\int_{\Omega\times X}\sup_n|F_n(\omega,x)|d\mu(\omega,x)<+\infty.$ If $\varphi:\mathbb N\to\mathbb N$ with $\varphi(n)\geq n,n\in\mathbb N,$ then
\begin{align*}
\lim\limits_{n\to\infty}\dfrac{1}{n}\sum\limits_{j=0}^{n-1}F_{\varphi(n)-j}(\Theta^j(\omega,x))=\mathbb E_\mu(F|\mathcal I_\mu)(\omega,x),
\end{align*}
almost everywhere and in $L^{1}$, where $\mathcal{I}_\mu=\{B\in\mathcal{F}\times \mathcal{B}:\mu(B\triangle \Theta^{-1}B)=0\}$ and $\mathbb E_\mu(F|\mathcal I_\mu)$ is the conditional expectation of $F$ with respect to $\mathcal I_\mu.$
\end{lemma}

Next, we give the weighted Shannon-McMillan-Breiman Theorem of RDS.

\begin{theorem}\label{weighted SMB theorem of RDS}
Let $f$ be a continuous bundle RDS over $(\Omega,\mathcal F,\mathbb P,\vartheta)$. Let $\alpha^{(1)},\alpha^{(2)}$ be two finite measurable partitions of $\Omega\times X$ and ${\bf a}=(a_1,a_2)\in\mathbb R^2$ where $a_1>0$ and $a_2\geq0$. Then
\begin{align*}
&\lim\limits_{N\to\infty}\dfrac{1}{N}I_{\mu_\omega}(\bigvee\limits_{n=0}^{\lceil a_1N\rceil-1}(f^n_\omega)^{-1}\alpha^{(1)}_{\vartheta^n\omega}\vee \bigvee\limits_{n=0}^{\lceil(a_1 +a_2)N\rceil-1}(f^n_\omega)^{-1}\alpha^{(2)}_{\vartheta^n\omega})(x)\\=&a_1\mathbb E_\mu(F_1| \mathcal I_\mu)(\omega,x)+a_2\mathbb E_\mu(F_2| \mathcal I_\mu)(\omega,x)
\end{align*}
almost everywhere, where
$F_1(\omega,x):=I_{\mu_\omega}(\alpha^{(1)}_\omega\vee\alpha^{(2)}_\omega|\bigvee\limits_{n=1}^\infty(f^n_\omega)^{-1}
(\alpha^{(1)}_{\vartheta^{n}\omega}\vee\alpha^{(2)}_{\vartheta^{n}\omega}))(x)$, $F_2(\omega,x):=I_{\mu_\omega}(\alpha^{(2)}_\omega|\bigvee\limits_{n=1}^\infty(f^n_\omega)^{-1}
(\alpha^{(2)}_{\vartheta^{n}\omega}))(x)$ and $\mathcal{I}_\mu=\{B\in\mathcal{F}\times \mathcal{B}:\mu(B\triangle \Theta^{-1}B)=0\}$. In addition, if $\mu$ is ergodic, then
\begin{align*}
&\lim\limits_{N\to\infty}\dfrac{1}{N}I_{\mu_\omega}(\bigvee\limits_{n=0}^{\lceil a_1N\rceil-1}(f^n_\omega)^{-1}\alpha^{(1)}_{\vartheta^n\omega}\vee \bigvee\limits_{n=0}^{\lceil(a_1 +a_2)N\rceil-1}(f^n_\omega)^{-1}\alpha^{(2)}_{\vartheta^n\omega})(x)\\=& a_1 h_{\mu}^{(r)}(f,\alpha^{(1)}\vee \alpha^{(2)})+a_2 h_{\mu}^{(r)}(f,\alpha^{(2)})
\end{align*}
almost everywhere.
\end{theorem}

\begin{proof}[Proof of Theorem \ref{weighted SMB theorem of RDS}]
First we show that for $a>0,b\geq0$ and a finite measurable partition $\beta$ of $\Omega\times X$,
\begin{align}\label{equation in the first case for RDS}
\lim\limits_{N\to\infty}\dfrac{1}{N}I_{\mu_\omega}(\bigvee\limits_{n= \lceil aN\rceil}^{\lceil (a+b)N\rceil-1 }(f^n_\omega)^{-1}\beta_{\vartheta^n\omega})(x)=b \mathbb E_\mu(G|\mathcal I_\mu)(\omega,x),
\end{align}
almost everywhere, where $G(\omega,x):=I_{\mu_\omega}(\beta_\omega|\bigvee_{n=1}^\infty(f^n_\omega)^{-1}\beta_{\vartheta^n\omega})(x)$.

It is obvious for $b=0$. Given $b>0$. Sine $\vartheta$ is invertible, for any $\omega\in \Omega$ and $n\in\mathbb{N}$ we have
\begin{align}\label{invariant measure formula}
f_{\omega}^{n}\mu_{\omega}=\mu_{\vartheta^{n}\omega}.
\end{align}
Remark that $G_{k}(\omega,x)=I_{\mu_\omega}(\beta_\omega|\bigvee_{n=1}^{k-1}(f^n_\omega)^{-1}\beta_{\vartheta^n\omega})(x)$ for each $k\in\mathbb{N}$. By making repeated use of the formula \eqref{invariant measure formula} and the standard property of the information function, we have
\begin{align*}
&I_{\mu_\omega}(\bigvee\limits_{n=\lceil aN\rceil}^{\lceil(a+b)N\rceil-1}(f^n_\omega)^{-1}\beta_{\vartheta^n\omega})(x)\\=&I_{\mu_\omega}(\bigvee\limits_{n=0}^{\lceil (a+b)N\rceil-1}(f^n_\omega)^{-1}\beta_{\vartheta^n\omega})(x)-I_{\mu_\omega}(\bigvee\limits_{n=0}^{\lceil aN\rceil-1}(f^n_\omega)^{-1}\beta_{\vartheta^n\omega}| \bigvee\limits_{n=\lceil aN\rceil}^{\lceil(a+b)N\rceil-1}(f^n_\omega)^{-1}\beta_{\vartheta^n\omega})(x).
\end{align*}
By the Shannon-McMillan-Breiman theorem for RDS in \cite{Zhu2008On}, in order to prove \eqref{equation in the first case for RDS}, we only need to prove
\begin{align}\label{equation in the second case for RDS}
\lim\limits_{N\to\infty}\dfrac{1}{N}I_{\mu_\omega}(\bigvee\limits_{n=0}^{\lceil aN\rceil-1}(f^n_\omega)^{-1}\beta_{\vartheta^n\omega}| \bigvee\limits_{n=\lceil aN\rceil}^{\lceil(a+b)N\rceil-1}(f^n_\omega)^{-1}\beta_{\vartheta^n\omega})(x)=a\mathbb E_\mu(G|\mathcal I_\mu)(\omega,x)
\end{align}
almost everywhere.

Note that
\begin{align*}
&I_{\mu_\omega}(\bigvee\limits_{n=0}^{\lceil aN\rceil-1}(f^n_\omega)^{-1}\beta_{\vartheta^n\omega}| \bigvee\limits_{n=\lceil aN\rceil}^{\lceil(a+b)N\rceil-1}(f^n_\omega)^{-1}\beta_{\vartheta^n\omega})(x)\\
=&G_{\lceil (a+b)N\rceil}(\omega,x)+I_{\mu_\omega}( \bigvee\limits_{n=1}^{\lceil aN\rceil-1}(f^n_\omega)^{-1}\beta_{\vartheta^n\omega}| \bigvee\limits_{n=\lceil aN\rceil}^{\lceil(a+b)N\rceil-1}(f^n_\omega)^{-1}\beta_{\vartheta^n\omega})(x)\\
=&G_{\lceil (a+b)N\rceil}(\omega,x)+I_{\mu_\omega}((f_\omega)^{-1} \bigvee\limits_{n=0}^{\lceil aN\rceil-2}(f^{n-1}_{\vartheta\omega})^{-1}\beta_{\vartheta^{n}\omega}|(f_\omega)^{-1}  \bigvee\limits_{n=\lceil aN\rceil-1}^{\lceil(a+b)N\rceil-2}(f^{n-1}_{\vartheta\omega})^{-1}\beta_{\vartheta^{n}\omega})(x)\\
=&G_{\lceil (a+b)N\rceil}(\omega,x)+I_{\mu_{\vartheta\omega}}(\bigvee\limits_{n=0}^{\lceil aN\rceil-2}(f^{n-1}_{\vartheta\omega})^{-1}\beta_{\vartheta^{n}\omega}|   \bigvee\limits_{n=\lceil aN\rceil-1}^{\lceil(a+b)N\rceil-2}(f^{n-1}_{\vartheta\omega})^{-1}\beta_{\vartheta^{n}\omega})(f_\omega x) \\
=&G_{\lceil (a+b)N\rceil}(\omega,x)+G_{\lceil (a+b)N\rceil-1}\circ\Theta(\omega,x)+\cdots+G_{\lceil (a+b)N\rceil-\lceil aN\rceil+1}\circ\Theta(\omega,x)\\
=&\sum\limits_{j=0}^{\lceil aN\rceil-1}G_{\lceil(a+b)N\rceil-j}\circ\Theta^j(\omega,x),
\end{align*}
where
$G_k(\omega,x):=I_{\mu_\omega}(\beta_\omega|\bigvee\limits_{n=1}^{k-1}(f^n_\omega)^{-1}\beta_{\vartheta^n\omega})(x).$
Since $G_k\in L^1(\Omega\times X,\mathcal F\times\mathcal{B},\mu),$  we get that
\begin{align*}
\lim\limits_{k\to\infty}G_k(\omega,x)=I_{\mu_\omega}(\beta_\omega|\bigvee\limits_{n=1}^\infty(f^n_\omega)^{-1}\beta_{\vartheta^n\omega})(x)=G(\omega,x)
\end{align*}
for $\mu$-a.e. $(\omega,x)\in \Omega\times X$ from the Martingale Convergence Theorem.
By Chung's lemma (see \cite{chung1961A}), $\int\sup_{n\in\mathbb N}G_n(\omega,x)d\mu(\omega,x)\leq \int H_{\mu_\omega}(\beta_\omega)d\mathbb P(\omega)+1<\infty.$
It follows from the Dominated Convergence Theorem that $G\in L^1(\Omega\times X,\mathcal F\times\mathcal B,\mu). $
By Lemma \ref{lemmaergodic}, we have
\begin{align*}
&\lim\limits_{N\to\infty}\dfrac{1}{N}I_{\mu_\omega}(\bigvee\limits_{n=0}^{\lceil aN\rceil-1}(f^n_\omega)^{-1}\beta_{\vartheta^n\omega}| \bigvee\limits_{n=\lceil aN\rceil}^{\lceil(a+b)N\rceil-1}(f^n_\omega)^{-1}\beta_{\vartheta^n\omega})(x)\\
=&a\lim\limits_{N\to\infty}\dfrac{1}{\lceil aN\rceil}\sum\limits_{j=0}^{\lceil aN\rceil-1}G_{\lceil(a+b)N\rceil-j}\circ\Theta^j(\omega,x)\\
=&a \mathbb E_\mu(G|\mathcal I_\mu)(\omega,x)
\end{align*}
almost everywhere. Hence \eqref{equation in the second case for RDS} holds, so does \eqref{equation in the first case for RDS}.

We will prove the theorem by considering the case of $a_{1}>0$, $a_{2}\geq 0$.
Let $\beta^{(1)}=\alpha^{(1)}\vee \alpha^{(2)}$ and $\beta^{(2)}=\alpha^{(2)}$. Then
$\beta^{(1)}\succeq \beta^{(2)}$, $F_1(\omega,x)=I_{\mu_\omega}(\beta^{(1)}_\omega|\bigvee\limits_{n=1}^\infty(f^n_\omega)^{-1}\beta^{(1)}_{\vartheta^{n}\omega})(x)$ and $F_2(\omega,x)=I_{\mu_\omega}(\beta^{(2)}_\omega|\bigvee\limits_{n=1}^\infty(f^n_\omega)^{-1}\beta^{(2)}_{\vartheta^{n}\omega})(x)$. We note that
\begin{align*}
&\bigvee\limits_{n=0}^{\lceil a_1 N\rceil-1}(f^n_\omega)^{-1}\alpha^{(1)}_{\vartheta^n\omega}\vee \bigvee\limits_{n=0}^{\lceil(a_1+a_2)N\rceil-1}(f^n_\omega)^{-1}\alpha^{(2)}_{\vartheta^n\omega}\\
=&\bigvee\limits_{n=0}^{\lceil a_1 N\rceil-1}(f^n_\omega)^{-1}\beta^{(1)}_{\vartheta^n\omega}
\vee\bigvee\limits_{n=\lceil a_1N\rceil }^{\lceil(a_1+a_{2})N\rceil-1 }(f^n_\omega)^{-1}\beta^{(2)}_{\vartheta^n\omega}.
\end{align*}
It follows from the Shannon-McMillan-Breiman theorem in \cite{Zhu2008On} and \eqref{equation in the first case for RDS} that
\begin{align*}
\lim\limits_{N\to\infty}\dfrac{1}{N}I_{\mu_\omega}(\bigvee\limits_{n=0}^{\lceil a_1 N\rceil-1}(f^n_\omega)^{-1}\beta^{(1)}_{\vartheta^n\omega})(x)= a_1\mathbb E_\mu(F_1|\mathcal I_\mu)(\omega,x)
\end{align*}
and
\begin{align*}
\lim\limits_{N\to\infty}\dfrac{1}{N}I_{\mu_\omega}(\bigvee\limits_{n=\lceil a_1N\rceil }^{\lceil(a_1+a_{2})N\rceil-1 }(f^n_\omega)^{-1}\beta^{(2)}_{\vartheta^n\omega})(x)= a_{2}\mathbb E_\mu(F_{2}|\mathcal I_\mu)(\omega,x)
\end{align*}
almost everywhere.
For $\mu$-a.e. $(\omega,x)\in\Omega\times X$, we define
\begin{align*}
Z_m^\omega(x)=\dfrac{\mu_\omega(\bigvee\limits_{n=0}^{\lceil
a_1 m\rceil-1}(f^n_\omega)^{-1}\beta^{(1)}_{\vartheta^n\omega}(x))
\cdot \mu_\omega(\bigvee\limits_{n=\lceil a_1 m\rceil }^{\lceil
(a_1+a_{2})m\rceil-1}(f^n_\omega)^{-1}\beta^{(2)}_{\vartheta^n\omega}(x))
}{\mu_\omega((\bigvee\limits_{n=0}^{\lceil
a_1 m\rceil-1}(f^n_\omega)^{-1}\beta^{(1)}_{\vartheta^n\omega}\vee\bigvee\limits_{n=\lceil a_1 m\rceil }^{\lceil
(a_1+a_{2})m\rceil-1}(f^n_\omega)^{-1}\beta^{(2)}_{\vartheta^n\omega})(x))
}
\end{align*}
for all $m\in\mathbb N,\omega\in\Omega.$
Since
\begin{align*}
\int_X Z_m^\omega(x)d\mu_\omega(x)
&=\sum\limits_{A\in \bigvee\limits_{n=0}^{\lceil
a_1 m\rceil-1}(f^n_\omega)^{-1}\beta^{(1)}_{\vartheta^n\omega}\atop B\in \bigvee\limits_{n=\lceil a_1 m\rceil)}^{\lceil
(a_1+a_{2})m\rceil-1}(f^n_\omega)^{-1}\beta^{(2)}_{\vartheta^n\omega}} \int_{A\cap B}\dfrac{\mu_\omega(A)\mu_\omega(B)}{\mu_\omega(A\cap B)}d\mu_\omega(x)\\&=\sum\limits_{A\in \bigvee\limits_{n=0}^{\lceil
a_1 m\rceil-1}(f^n_\omega)^{-1}\beta^{(1)}_{\vartheta^n\omega}\atop B\in \bigvee\limits_{n=\lceil a_1 m\rceil)}^{\lceil
(a_1+a_{2})m\rceil-1}(f^n_\omega)^{-1}\beta^{(2)}_{\vartheta^n\omega}} \mu_\omega(A)\mu_\omega(B)\\&=1,
\end{align*}
$\sum_{m=1}^\infty\mu(\{(\omega,x)\in \Omega\times X:Z_m^\omega(x)\geq \exp(\epsilon m)\})$ converges for every $\epsilon>0$ and the Borel-Canteli Lemma implies that for $\mu$-a.e. $(\omega,x)\in \Omega\times X$,
\begin{align*}
\limsup\limits_{m\to\infty}\dfrac{1}{m}\log Z^\omega_m(x)\leq0.
\end{align*}
Then
\begin{align*}
&\limsup\limits_{N\to\infty}\dfrac{1}{N}I_{\mu_\omega}(
\bigvee\limits_{n=0}^{\lceil a_1 N\rceil-1}(f^n_\omega)^{-1}\alpha^{(1)}_{\vartheta^n\omega}\vee \bigvee\limits_{n=0}^{\lceil (a_1+a_{2}) N\rceil-1}(f^n_\omega)^{-1}\alpha^{(2)}_{\vartheta^n\omega})
\\ \leq& a_1\mathbb E_{\mu}(F_1|\mathcal I_\mu)(\omega,x)+a_2\mathbb E_{\mu}(F_2|\mathcal I_\mu)(\omega,x)
\end{align*}
for $\mu$-a.e. $\Omega\times X.$

Conversely, by \eqref{equation in the first case for RDS}, we have
\begin{align*}
\lim\limits_{N\to\infty}\dfrac{1}{N}I_{\mu_\omega}(\bigvee\limits_{n=\lceil a_1 N\rceil}^{\lceil(a_1+a_{2})N\rceil-1}(f^n_\omega)^{-1}\beta^{(1)}_{\vartheta^n\omega})
(x)&=a_{2}\mathbb E_\mu(F_1|\mathcal I_\mu)(\omega,x),\\
\lim\limits_{N\to\infty}\dfrac{1}{N}I_{\mu_\omega}(\bigvee\limits_{n=\lceil a_1 N\rceil}^{\lceil(a_1+a_{2})N\rceil-1}(f^n_\omega)^{-1}\beta^{(2)}_{\vartheta^n\omega})
(x)&=a_{2}\mathbb E_\mu(F_2|\mathcal I_\mu)(\omega,x),
\end{align*}
and
\begin{align*}
&\lim\limits_{N\to\infty}\dfrac{1}{N}I_{\mu_\omega}(
\bigvee\limits_{n=0}^{\lceil(a_1+a_{2})N\rceil-1}(f^n_\omega)^{-1}\beta^{(1)}_{\vartheta^n\omega}
)(x)
=(a_1+a_{2})\mathbb E_\mu(F_1|\mathcal I_\mu)(\omega,x)
\end{align*}
almost everywhere.
Now for $\mu$-a.e. $(\omega,x)\in \Omega\times X$, we define for $m\in\mathbb N,\omega\in\Omega,$
\begin{align*}
R^\omega_m(x)=&\dfrac{\mu_\omega((\bigvee\limits_{n=0}^{\lceil a_1 m\rceil-1}(f^n_\omega)^{-1}\beta^{(1)}_{\vartheta^n\omega}\vee
\bigvee\limits_{n=\lceil a_1 m\rceil}^{\lceil(a_1+a_{2})m\rceil-1}(f^n_\omega)^{-1}\beta^{(2)}_{\vartheta^n\omega})(x))}{\mu_\omega((
\bigvee\limits_{n=0}^{\lceil(a_1+a_{2})m\rceil-1}(f^n_\omega)^{-1}\beta^{(1)}_{\vartheta^n\omega})
(x))}\\
\times&\dfrac{\mu_\omega((\bigvee\limits_{n=\lceil a_1 m\rceil}^{\lceil(a_1+a_{2})m\rceil-1}
(f^n_\omega)^{-1}\beta^{(1)}_{\vartheta^n\omega})(x))}
{\mu_\omega((\bigvee\limits_{n=\lceil a_1 m\rceil}^{\lceil(a_1+a_{2})m\rceil-1}(f^n_\omega)^{-1}\beta^{(2)}_{\vartheta^n\omega})(x))}.
\end{align*}
Then for $\mu$-a.e. $(\omega,x)\in \Omega\times X, R^\omega_m(x)>0.$
Since $\beta_1\succeq \beta_{2},$ we have
\begin{align*}
\int_X R^\omega_m(x)d\mu_\omega(x)
&=\sum_{{A\in\bigvee\limits_{n=0}^{\lceil a_1 m\rceil-1}(f^n_\omega)^{-1}\beta^{(1)}_{\vartheta^n\omega}\atop B\in \bigvee\limits_{n=\lceil a_1 m\rceil}^{\lceil(a_1+a_{2})m\rceil-1}(f^n_\omega)^{-1}\beta^{(2)}_{\vartheta^n\omega}}\atop C\in \bigvee\limits_{n=\lceil a_1 m\rceil}^{\lceil(a_1+a_{2})m\rceil-1}
(f^n_\omega)^{-1}\beta^{(1)}_{\vartheta^n\omega}}\int_{A\cap B\cap C}\dfrac{\mu_\omega(A\cap B)\mu_\omega(B\cap C)}{\mu_\omega(A\cap B\cap C)\mu_\omega(B)}d\mu_\omega(x)\\&=\sum_{{A\in\bigvee\limits_{n=0}^{\lceil a_1 m\rceil-1}(f^n_\omega)^{-1}\beta^{(1)}_{\vartheta^n\omega}\atop B\in \bigvee\limits_{n=\lceil a_1 m\rceil}^{\lceil(a_1+a_{2})m\rceil-1}(f^n_\omega)^{-1}\beta^{(2)}_{\vartheta^n\omega}}\atop C\in \bigvee\limits_{n=\lceil a_1 m\rceil}^{\lceil(a_1+a_{2})m\rceil-1}
(f^n_\omega)^{-1}\beta^{(1)}_{\vartheta^n\omega}}\dfrac{\mu_\omega(A\cap B)\mu_\omega(B\cap C)}{\mu_\omega(B)}
\\&=1.
\end{align*}
Hence, $\sum_{m=1}^\infty\mu(\{(\omega,x)\in \Omega\times X:R^\omega_m(x)\geq \exp(\epsilon m)\})$ converges for every $\epsilon>0$ and the Borel-Canteli Lemma implies that $\limsup_{N\to\infty}\dfrac{1}{N}\log R^\omega_N(x)\leq0$ for $\mu$-a.e. $(\omega,x)\in\Omega\times X.$ Then
\begin{align*}
&\liminf\limits_{N\to\infty}\dfrac{1}{N}I_{\mu_\omega}(
\bigvee\limits_{n=0}^{\lceil a_1 N\rceil-1}(f^n_\omega)^{-1}\alpha^{(1)}_{\vartheta^n\omega}\vee \bigvee\limits_{n=0}^{\lceil (a_1+a_{2}) N\rceil-1}(f^n_\omega)^{-1}\alpha^{(2)}_{\vartheta^n\omega})
\\ \geq&a_1\mathbb E_{\mu}(F_1|\mathcal I_\mu)(\omega,x)+a_2\mathbb E_{\mu}(F_2|\mathcal I_\mu)(\omega,x)
\end{align*}
for $\mu$-a.e. $(\omega,x)\in\Omega\times X.$
\end{proof}

Especially, we have the following result.
\begin{corollary}\label{application for weighted SMB theorem of RDS}
Let $f$ be a continuous bundle RDS over $(\Omega,\mathcal F,\mathbb P,\vartheta)$. Let $\alpha^{(1)},\alpha^{(2)}$ be two finite measurable partitions of $\Omega\times X$ and ${\bf a}=(a_1,a_2)\in\mathbb R^2$ where $a_1>0$ and $a_2\geq0$. if $\mu\in \mathcal{E}_{\mathbb{P}}^{1}(\Omega\times X,f)$ and $\alpha^{(1)}\succeq \alpha^{(2)},$ then for $\mu$-a.e. $(\omega,x)\in \Omega\times X$,
\begin{align}\label{equation for weighted SMB theorem of RDS}
\begin{split}
&\lim\limits_{N\to\infty}\dfrac{1}{N}I_{\mu_\omega}(\bigvee\limits_{n=0}^{\lceil a_1 N\rceil-1}(f^n_\omega)^{-1}\alpha^{(1)}_{\vartheta^n\omega}\vee \bigvee\limits_{n=\lceil a_1 N\rceil}^{\lceil (a_1+a_{2}) N\rceil-1}(f^n_\omega)^{-1}\alpha^{(2)}_{\vartheta^n\omega})(x)\\=& a_1h_{\mu}^{(r)}(f,\alpha^{(1)})+a_2h_{\mu}^{(r)}(f,\alpha^{(2)}).
\end{split}
\end{align}
\end{corollary}
\begin{proof}
Since $\alpha^{(1)}\succeq \alpha^{(2)}$, we have
\begin{align*}
&\bigvee\limits_{n=0}^{\lceil a_1 N\rceil-1}(f^n_\omega)^{-1}\alpha^{(1)}_{\vartheta^n\omega}\vee \bigvee\limits_{n=0}^{\lceil (a_1+a_{2}) N\rceil-1}(f^n_\omega)^{-1}\alpha^{(1)}_{\vartheta^n\omega}\\=&\bigvee\limits_{n=0}^{\lceil a_1 N\rceil-1}(f^n_\omega)^{-1}\alpha^{(1)}_{\vartheta^n\omega}\vee \bigvee\limits_{n=\lceil a_1 N\rceil}^{\lceil (a_1+a_{2}) N\rceil-1}(f^n_\omega)^{-1}\alpha^{(2)}_{\vartheta^n\omega},
\end{align*}
and
\begin{align*}
\alpha^{(1)}\vee\alpha^{(2)}&=\alpha^{(1)}.
\end{align*}
By Themorem \ref{weighted SMB theorem of RDS}, and [\cite{Zhu2008On}, Proposition 2.1], it follows that
\begin{align*}
&\lim\limits_{N\to\infty}\dfrac{1}{N}I_{\mu_\omega}(\bigvee\limits_{n=0}^{\lceil a_1 N\rceil-1}(f^n_\omega)^{-1}\alpha^{(1)}_{\vartheta^n\omega}\vee \bigvee\limits_{n=\lceil a_1 N\rceil}^{\lceil (a_1+a_{2}) N\rceil-1}(f^n_\omega)^{-1}\alpha^{(2)}_{\vartheta^n\omega})(x)\\=& a_1\mathbb E_\mu(F_1| \mathcal I_\mu)(\omega,x)+a_2\mathbb E_\mu(F_2| \mathcal I_\mu)(\omega,x)\\=& a_{1} h_\mu^{(r)}(f,\alpha^{(1)}\vee \alpha^{(2)})+a_{2} h_\mu^{(r)}(f, \alpha^{(2)})
\end{align*}
almost everywhere, where
$F_1(\omega,x):=I_{\mu_\omega}(\alpha^{(1)}_\omega\vee \alpha^{(2)}_\omega|\bigvee\limits_{n=1}^\infty(f^n_\omega)^{-1}
(\alpha^{(1)}_{\vartheta^{n}\omega}\vee \alpha^{(2)}_{\vartheta^{n}\omega}))(x)$ and $F_2(\omega,x):=I_{\mu_\omega}(\alpha^{(2)}_\omega|\bigvee\limits_{n=1}^\infty(f^n_\omega)^{-1}
(\alpha^{(2)}_{\vartheta^{n}\omega}))(x)$
This implies \eqref{equation for weighted SMB theorem of RDS}.
\end{proof}

\section{Appendix B: weighted Brin-Katok local entropy formula}\label{Appendix B: a weighted version of the Brin-Katok theorem}

In this appendix B, we give the proof of Theorem \ref{Brin-Katok local entropy formula}.

\begin{proof}[Proof of Theorem \ref{Brin-Katok local entropy formula}]
We first prove
\begin{align*}
\lim\limits_{\epsilon\to0}\limsup\limits_{n\to\infty}-\dfrac{1}{n}\log\mu_\omega(B^{{\bf a}}_{\omega}(x, n,\epsilon))\leq h_\mu^{{\bf a}}(f_1).
\end{align*}
Given $\epsilon>0.$ Let $\xi^{(i)}$ be a finite measurable partition of $X_i,i=1,2$ with ${\rm diam} (\xi^{(i)})<\epsilon.$ Set $\hat{\xi}^{(i)}=\{\Omega\times A:A\in \xi^{(i)}\}$. Then
\begin{align*}
B^{\bf a}_\omega(x,n,\epsilon)\supseteq
(\bigvee\limits_{j=0}^{\lceil a_1 n\rceil-1}(f^j_{1,\omega})^{-1}\xi^{(1)}\vee \bigvee\limits_{j=0}^{\lceil (a_1+a_{2}) n\rceil-1}(f^j_{1,\omega})^{-1}\pi^{-1}\xi^{(2)})(x),
\end{align*}
for $\omega\in\Omega,x\in X_{1}.$
By Theorem \ref{weighted SMB theorem of RDS}, it follows that
\begin{align*}
&\limsup\limits_{n\to\infty}-\dfrac{1}{n}\log\mu_\omega(B^{\bf a}_\omega(x,n,\epsilon))\\
\leq& \limsup\limits_{n\to\infty}-\dfrac{1}{n}\log\mu_\omega((\bigvee\limits_{j=0}^{\lceil a_1 n\rceil-1}(f^j_{1,\omega})^{-1}\xi^{(1)}\vee \bigvee\limits_{j=0}^{\lceil (a_1+a_{2}) n\rceil-1}(f^j_{1,\omega})^{-1}\pi^{-1}\xi^{(2)})(x))\\
=&\limsup\limits_{n\to\infty}\dfrac{1}{n}I_{\mu_{\omega}}(\bigvee\limits_{j=0}^{\lceil a_1 n\rceil-1}(f^j_{1,\omega})^{-1}\xi^{(1)}\vee \bigvee\limits_{j=0}^{\lceil (a_1+a_{2}) n\rceil-1}(f^j_{1,\omega})^{-1}\pi^{-1}\xi^{(2)})(x)\\
=&a_1 h_\mu^{(r)}(f_1,\hat{\xi}^{(1)}\vee\Pi^{-1}\hat{\xi}^{(2)})+a_2 h_\mu^{(r)}(f_1,\Pi^{-1}\hat{\xi}^{(2)})\\
\leq & a_1h^{(r)}_{\mu}(f_1)+ a_2h^{(r)}_{\mu\circ\Pi^{-1}}(f_2)=h_\mu^{{\bf a}}(f_1).
\end{align*}
Letting $\epsilon\to0,$ we complete the proof of the upper bound.

Next, we show the opposite inequality. It is sufficient to show that for any $\delta>0$, there exists $\epsilon>0$ and a measurable subset $D$ of $\Omega\times X_{1}$ such that $\mu(D)>1-3\delta$ and
\begin{align*}
\liminf\limits_{n\to\infty}\dfrac{-\log\mu_{\omega}(B_\omega^{\bf a}(x,n,\epsilon))}{n}\geq \min\left\{\frac{1}{\delta},h_\mu^{{\bf a}}(f_1)-\delta\right\}-2(1+a_1+a_2)\delta
\end{align*}
for any $(\omega,x)\in D$.

Fix $\delta>0$. We are going to find such $\epsilon$ and $D$. Choose a finite measurable partition $\beta^{(i)}$ of $X_{i}$ for $i=1,2$ such that $\int \mu_{\omega}(\partial \beta^{(1)})d\mathbb{P}=0$, $\int \mu_{\omega}\circ\pi^{-1}(\partial \beta^{(2)})d\mathbb{P}=0$, $\text{diam}(\beta^{(i)})$ are small enough so that
\begin{align*}
h^{(r)}_{\mu}(f_1,\hat{\beta}^{(1)})\geq\left\{
\begin{aligned}
\frac{1}{a_{1}\delta} \ \ \ \ \ \ \  \ \ \ \ \ \ \ \ \ \  \ \ &\text{if} \ h_{\mu}^{(r)}(f_{1})=\infty,\\
h^{(r)}_{\mu}(f_1)-\dfrac{\delta}{a_1} \ \ \ &\text{otherwise},
\end{aligned}
\right.
\end{align*}
and
\begin{align*}
h^{(r)}_{\mu\circ \Pi^{-1}}(f_2,\hat{\beta}^{(2)})\geq\left\{
\begin{aligned}
\frac{1}{a_{1}\delta} \ \ \ \ \ \ \ \ \ \ \ \ \ \ \ \ \ \ \ \ \ \ \ \ \ \ \ \ \ \ \ \ \ &\text{if} \ h_{\mu\circ\Pi^{-1}}^{(r)}(f_{2})=\infty,\\
h^{(r)}_{\mu\circ \Pi^{-1}}(f_2)-\dfrac{\delta}{a_1+a_2} \ \ \ &\text{otherwise},
\end{aligned}
\right.
\end{align*}
where $\hat{\beta}^{(i)}=\{\Omega\times A:A\in \beta^{(i)}\}$ for $i=1,2$.

Next we define the partitions $\alpha^{(i)}$ of $X_i$ for $i=2,1$, by setting $\alpha^{(2)}=\beta^{(2)}$
and $\alpha^{(1)}=\beta^{(1)}\vee\Pi^{-1}(\alpha^{(2)})$. Then for any $\omega\in \Omega$, we have
 \begin{enumerate}
  \item[(1)] $\alpha^{(1)}\succeq\Pi^{-1}\alpha^{(2)}$;
  \item[(2)] $a_1h^{(r)}_{\mu}(f_1,\hat{\alpha}^{(1)})+a_2h^{(r)}_{\mu\circ\Pi^{-1}}(f_2,\hat{\alpha}^{(2)})\geq \min\left\{\frac{1}{\delta},h^{{\bf a}}_{\mu}(f_1)-\delta\right\};$
  \item[(3)] $\mu_\omega(\partial\alpha^{(1)})=0$ and $\mu_\omega\circ\pi^{-1}(\partial\alpha^{(2)})=0$,
\end{enumerate}
where $\hat{\alpha}^{(i)}=\{\Omega\times A:A\in \alpha^{(i)}\}$ for $i=1,2$.

In order to estimate the $\mu_{\omega}$-measure of $B^{\bf a}_\omega(x,n,\epsilon)$, we need the following Hamming metrics. Write $\alpha^{(i)}=\{A_1^i,A_2^i\cdots,A^i_{u_i}\}$ for $ i=1,2.$ Let $M=\max\{u_1,u_{2}\}$ and $\Gamma=\{1,2,\cdots,M\}.$ Given $m\in\mathbb N,$ for ${\bf s}=(s_i)_{i=0}^{m-1},{\bf t}=(t_i)_{i=0}^{m-1}\in\Gamma^{\{0,1,\cdots,m-1\}},$ the Hamming distance between ${\bf s}$ and  ${\bf t}$ is defined to be
\begin{align*}
\dfrac{1}{m}\#\{i\in\{0,1,\cdots,m-1\}:s_i\neq t_i\}.
\end{align*}
For ${\bf s}\in\Gamma^{\{0,1,\cdots,m-1\}}$ and $0<\tau<1,$ set $Q({\bf s},\tau)$ be the total number of those ${\bf t}$ so that the Hamming distance between ${\bf s}$ and ${\bf t}$ does not exceed $\tau.$ Obviously,
\begin{align*}
Q_m(\tau):=\max\{Q({\bf s},\tau):{\bf s}\in\Gamma^{\{0,1,\cdots,m-1\}}\}\leq
\binom{m}{\lceil m\tau\rceil}M^{\lceil m\tau\rceil}.
\end{align*}
By Stirling formula, there exists a $\delta_1>0$ and a positive constant $C>0$ such that for all $m\in\mathbb N,$
\begin{align}\label{estiamtion of total number}
\binom{m}{\lceil m\delta_{1}\rceil}M^{\lceil m\delta_{1}\rceil}\leq e^{\delta m+C}.
\end{align}

For $\eta>0,$ let
\begin{align*}
U^1_\eta(\alpha^{(1)})&=\{x\in X_{1}:B(x,\eta)\setminus \alpha^{(1)}(x)\neq\emptyset\},\\
U^2_\eta(\alpha^{(2)})&=\{x\in X_{1}:B(\pi x,\eta)\setminus \alpha^{(2)}(\pi x)\neq\emptyset\}.
\end{align*}
It is clear that, for $x\in U^1_\eta(\alpha^{(1)})$, the $\epsilon$-ball about $x$ is not contained in the element of $\alpha^{(1)}$ which $x$ belongs, and for $x\in U^2_\eta(\alpha^{(2)})$, the $\epsilon$-ball about $\pi x$ is not contained in the element of $\alpha^{(2)}$ which $\pi x$ belongs. Moreover, set $\hat{U}^1_\eta(\alpha^{(1)})=\Omega\times U^1_\eta(\alpha^{(1)})$ and $\hat{U}^2_\eta(\alpha^{(2)})=\Omega\times U^2_\eta(\alpha^{(2)})$.
Since for $\bigcap_{\eta>0}U^1_\eta(\alpha^{(1)})=\partial\alpha^{(1)}$ and $\bigcap_{\eta>0}U^2_\eta(\alpha^{(2)})=\pi^{-1}(\partial\alpha^{(2)})$, by the construction of $\alpha^{(i)}$, we have $\lim_{\eta\to0}\mu(\hat{U}^1_\eta(\alpha^{(1)}))=\mu(\Omega\times \partial\alpha^{(1)})=0$ and $\lim_{\eta\to0}\mu(\hat{U}^2_\eta(\alpha^{(2)}))=\mu(\Omega\times \pi^{-1}(\partial\alpha^{(2)}))=0$. Therefore, we can take $\epsilon>0$ such that for any $0<\eta\leq\epsilon, i=1,2$,
\begin{align*}
\mu(\hat{U}^i_\eta(\alpha^{(i)}))<\delta_1.
\end{align*}
By the Birkhoff ergodic theorem, for $\mu$-a.e. $(\omega,x)\in \Omega\times X_1,$
\begin{align*}
&\lim\limits_{n\to\infty}\dfrac{1}{\lceil(a_1+a_2)n\rceil}(\sum\limits_{j=0}^{\lceil a_1 n\rceil-1}\chi_{\hat{U}^1_\epsilon(\alpha^{(1)})}(\Theta_1^j(\omega,x))+\sum\limits_{j=\lceil a_{1}n\rceil}^{\lceil (a_1+a_{2}) n\rceil-1}\chi_{\hat{U}^2_\epsilon(\alpha^{(2)})}(\Theta_1^j(\omega,x)))\\=&\dfrac{1}{(a_1+a_2)}(a_1\mu(\hat{U}^1_\epsilon(\alpha^{(1)}))+a_2\mu(\hat{U}^2_\epsilon(\alpha^{(2)})))\\<&\delta_1.
\end{align*}
We can take a number $l_0$ large enough such that $\mu(A_l)>1-\delta$ for any $l\geq l_0,$ where
\begin{align*}
A_l=&\{(\omega,x)\in \Omega\times X_1:\dfrac{1}{\lceil(a_1+a_2)n\rceil}(\sum\limits_{j=0}^{\lceil a_1 n\rceil-1}\chi_{\hat{U}^1_\epsilon(\alpha^{(1)})}(\Theta_1^j(\omega,x))\\&+\sum\limits_{j=\lceil a_{1}n\rceil}^{\lceil (a_1+a_{2}) n\rceil-1}\chi_{\hat{U}^2_\epsilon(\alpha^{(2)})}(\Theta_1^j(\omega,x)))\leq\delta_1 \text{ for all }n\geq l\}.
\end{align*}
Since $\hat{\alpha}^{(1)} \succeq \Pi^{-1}\hat{\alpha}^{(2)}$, by Corollary \ref{application for weighted SMB theorem of RDS} we have
\begin{align*}
&\lim\limits_{n\to+\infty}-\dfrac{1}{n}\log\mu_\omega((\bigvee\limits_{j=0}^{\lceil a_1 n\rceil-1}(f_{1,\omega}^j)^{-1}\alpha^{(1)}\vee \bigvee\limits_{j=\lceil a_{1}n\rceil}^{\lceil (a_1+a_{2}) n\rceil-1}(f_{1,\omega}^j)^{-1}\pi^{-1}\alpha^{(2)})(x))\\
=&a_1 h^{(r)}_\mu(f_1,\alpha^{(1)})+a_2 h^{(r)}_\mu(f_1,\Pi^{-1}\alpha^{(2)})\\=&a_1 h^{(r)}_\mu(f_1,\alpha^{(1)})+a_2 h^{(r)}_{\mu\circ\Pi^{-1}}(f_2,\alpha^{(2)}).
\end{align*}
There exists a large $l_1$ such that $\mu(B_l)>1-\delta$ for any $l\geq l_1,$ where
\begin{align*}
B_l=&\{(\omega,x):-\dfrac{1}{n}\log\mu_\omega((\bigvee\limits_{j=0}^{\lceil a_1 n\rceil-1}(f_{1,\omega}^j)^{-1}\alpha^{(1)}\vee \bigvee\limits_{j=\lceil a_{1}n\rceil}^{\lceil (a_1+a_{2}) n\rceil-1}(f_{1,\omega}^j)^{-1}\pi^{-1}\alpha^{(2)})(x))\\
\\&\geq a_1 h^{(r)}_\mu(f_1,\alpha^{(1)})+a_2 h^{(r)}_{\mu\circ\Pi^{-1}}(f_2,\alpha^{(2)})-\delta,n\geq l\}.
\end{align*}
Fix $l\geq\max\{l_0,l_1\}.$ Set $E=A_l\cap B_l.$ Then $\mu(E)>1-2\delta.$ For $x\in X_1,n\in\mathbb N,\omega\in\Omega,$ the unique element
\begin{align*}
C(n,\omega,x)=(C_j(n,\omega,x))_{j=0}^{\lceil(a_1+a_2)n\rceil-1}\in \Gamma^{\{0,1,\cdots,\lceil(a_1+a_2)n\rceil-1\}}
\end{align*}
satisfying that $f_{1,\omega}^jx\in A^i_{C_j(n,\omega,x)}$ for $0\leq j\leq
\lceil a_1 n\rceil-1$ and $f_{1,\omega}^jx\in\pi^{-1}A^i_{C_j(n,\omega,x)}$ for $\lceil a_{1}n\rceil\leq j\leq
\lceil(a_1+a_{2})n\rceil-1$, is called the $(\omega,\{\alpha^{(i)}\}_{i=1}^2,{\bf a},n)$-name of $x.$ Since each point in one atom $A$ of $\bigvee\limits_{j=0}^{\lceil a_1 n\rceil-1}(f_{1,\omega}^j)^{-1}\alpha^{(1)}\vee \bigvee\limits_{j=\lceil a_{1}n\rceil}^{\lceil (a_1+a_{2}) n\rceil-1}(f_{1,\omega}^j)^{-1}\pi^{-1}\alpha^{(2)}$
has the same $(\omega,\{\alpha^{(i)}\}_{i=1}^2,{\bf a},n)$-name, set $C(n,\omega,A):=C(n,\omega,x)$ for any $x\in A,$ which is called the $(\omega,\{\alpha^{(i)}\}_{i=1}^2,{\bf a},n)$-name of $A.$

By the definition of $U^i_\epsilon(\alpha^{(i)})$, if $y\in B_\omega^{\bf a}(x,n,\epsilon),$ then for $0\leq j\leq
\lceil a_1 n\rceil-1,$ either $f^j_{1,\omega}x$ and $f^j_{1,\omega}y$ belong to  the same element of $\alpha^{(1)}$ or $f^j_{1,\omega}x\in U_\epsilon^1(\alpha^{(1)})$, for $\lceil a_{1}n\rceil\leq j\leq
\lceil (a_1+a_{2}) n\rceil-1,$ either $f^j_{1,\omega}x$ and $f^j_{1,\omega}y$ belong to  the same element of $\pi^{-1}\alpha^{(2)}$ or $f^j_{1,\omega}x\in U_\epsilon^2(\alpha^{(2)})$.
Hence if $(\omega,x)\in E,n\geq l$ and $y\in B_\omega^{\bf a}(x,n,\epsilon),$ then the Hamming distance between $(\omega,\{\alpha^{(i)}\}_{i=1}^2,{\bf a},n)$-name of $x$ and $y$ does not exceed $\delta_1$. Furthermore, $B_\omega^{\bf a}(x,n,\epsilon)$ is contained in the set of points $y$ whose $(\omega,\{\alpha^{(i)}\}_{i=1}^2,{\bf a},n)$-name is $\delta_1$-close to   $(\omega,\{\alpha^{(i)}\}_{i=1}^2,{\bf a},n)$-name of $x.$ Clearly, the total number $L_n(x)$ of such
$(\omega,\{\alpha^{(i)}\}_{i=1}^2,{\bf a},n)$-names admits the following estimate:
\begin{align*}
L_n(x)&\leq \binom{\lceil(a_1+a_2)n\rceil}{\lceil\lceil(a_1+a_2)n\rceil\delta_1\rceil}M^{\lceil\lceil(a_1+a_2)n\rceil\delta_1\rceil}\\
&\leq e^{\delta\lceil (a_1+a_2)n\rceil+C}\\
&\leq e^{\delta (a_1+a_2)n +C+\delta}.
\end{align*}
where the second inequality comes from \eqref{estiamtion of total number}. Moreover, we have shown that for any $(\omega,x)\in E$ and $n\geq l,$
\begin{align*}
B_\omega^{\bf a}(x, n,\epsilon)\subset&\{ y\in X_1:C(n,\omega,y) \text{ is }  \delta_1\text{-close to } C(n,\omega,x)\}\\
=&\bigcup\{A\in\bigvee\limits_{j=0}^{\lceil a_1 n\rceil-1}(f_{1,\omega}^j)^{-1}\alpha^{(1)}\vee \bigvee\limits_{j=\lceil a_{1}n\rceil}^{\lceil (a_1+a_{2}) n\rceil-1}(f_{1,\omega}^j)^{-1}\pi^{-1}\alpha^{(2)}
:\\ &C(n,\omega,A) \text{ is }\delta_1 \text{-close to } C(n,\omega,x)\}.
\end{align*}
And the cardinality of the right item is less than $L_n(x).$

Now for $n\in\mathbb N,$ set $E_{n}=\bigcup_{\omega\in \Omega} \{\omega\}\times E_{n,\omega}$ where $E_{n,\omega}$ be the set of $x\in E_{\omega}$ such that there exists an element $A\in\bigvee\limits_{j=0}^{\lceil a_1 n\rceil-1}(f_{1,\omega}^j)^{-1}\alpha^{(1)}\vee \bigvee\limits_{j=\lceil a_{1}n\rceil}^{\lceil (a_1+a_{2}) n\rceil-1}(f_{1,\omega}^j)^{-1}\pi^{-1}\alpha^{(2)}
$ with
\begin{align*}
\mu_\omega(A)>\exp(-(a_1 h^{(r)}_\mu(f_1,\alpha^{(1)})+a_2 h^{(r)}_{\mu\circ\Pi^{-1}}(f_2,\alpha^{(2)}))n+(2+a_1+a_2)\delta n)
\end{align*}
and the $(\omega,\{\alpha^{(i)}\}_{i=1}^2,{\bf a},n)$-name of $A$ is $\delta_1$-close to the $(\omega,\{\alpha^{(i)}\}_{i=1}^2,{\bf a},n)$-name of $x$.
If $(\omega,x)\in E\setminus E_n,$ then for each $A\in \bigvee\limits_{j=0}^{\lceil a_1 n\rceil-1}(f_{1,\omega}^j)^{-1}\alpha^{(1)}\vee \bigvee\limits_{j=\lceil a_{1}n\rceil}^{\lceil (a_1+a_{2}) n\rceil-1}(f_{1,\omega}^j)^{-1}\pi^{-1}\alpha^{(2)}
$ whose $(\omega,\{\alpha^{(i)}\}_{i=1}^2,{\bf a},n)$-name is $\delta_1$-close to $(\omega,\{\alpha^{(i)}\}_{i=1}^2,{\bf a},n)$-name of $x,$ one has
\begin{align*}
\mu_\omega(A)\leq\exp(-(a_1 h^{(r)}_\mu(f_1,\alpha^{(1)})+a_2 h^{(r)}_{\mu\circ\Pi^{-1}}(f_2,\alpha^{(2)}))n+(2+a_1+a_2)\delta n).
\end{align*}
Now we get ready to estimate $\mu_{\omega}$-measure of the set of points $x\in E_{n,\omega}$ for $n\geq l.$ Let $n\geq l.$ Put
\begin{align*}
\mathcal F_n^\omega=&\{A\in\bigvee\limits_{j=0}^{\lceil a_1 n\rceil-1}(f_{1,\omega}^j)^{-1}\alpha^{(1)}\vee \bigvee\limits_{j=\lceil a_{1}n\rceil}^{\lceil (a_1+a_{2}) n\rceil-1}(f_{1,\omega}^j)^{-1}\pi^{-1}\alpha^{(2)}
:\mu_\omega(A)\\&>\exp(-(a_1 h^{(r)}_\mu(f_1,\alpha^{(1)})+a_2 h^{(r)}_{\mu\circ\Pi^{-1}}(f_2,\alpha^{(2)}))n+(2+a_1+a_2)\delta n)\}.
\end{align*}
For the set $\mathcal F_n^\omega$, we observe that
\begin{align*}
1&\geq\mu_{\omega}(\bigcup_{A\in \mathcal F_n^\omega} A)=\sum_{A\in \mathcal F_n^\omega}\mu_{\omega}(A)\\&\geq\#\mathcal F_n^\omega\cdot\exp(-(a_1 h^{(r)}_\mu(f_1,\alpha^{(1)})+a_2 h^{(r)}_{\mu\circ\Pi^{-1}}(f_2,\alpha^{(2)}))n+(2+a_1+a_2)\delta n)\}.
\end{align*}
Then
\begin{align}\label{third estimation for fiber measure}
\#\mathcal F_n^\omega\leq \exp((a_1 h^{(r)}_\mu(f_1,\alpha^{(1)})+a_2 h^{(r)}_{\mu\circ\Pi^{-1}}(f_2,\alpha^{(2)}))n-(2+a_1+a_2)\delta n).
\end{align}
Let $x\in E_{n,\omega}.$ On the one hand since $x\in B_{l,\omega}$, we have
\begin{align*}
&\mu_\omega((\bigvee\limits_{j=0}^{\lceil a_1 n\rceil-1}(f_{1,\omega}^j)^{-1}\alpha^{(1)}\vee \bigvee\limits_{j=\lceil a_{1}n\rceil}^{\lceil (a_1+a_{2}) n\rceil-1}(f_{1,\omega}^j)^{-1}\pi^{-1}\alpha^{(2)})
(x))\\ \leq &\exp((-(a_1 h^{(r)}_\mu(f_1,\alpha^{(1)})+a_2 h^{(r)}_{\mu\circ\Pi^{-1}}(f_2,\alpha^{(2)}))+\delta)n).
\end{align*}
On the other hand, there exist $A\in\mathcal F_n^\omega $ with the $(\omega,\{\alpha^{(i)}\}_{i=1}^2,{\bf a},n)$-name of $A$ is $\delta_1$-close to the   $(\omega,\{\alpha^{(i)}\}_{i=1}^2,{\bf a},n)$-name of $x$, that is
the $(\omega,\{\alpha^{(i)}\}_{i=1}^2,{\bf a},n)$-name of $A$ is $\delta_1$-close to the   $(\omega,\{\alpha^{(i)}\}_{i=1}^2,{\bf a},n)$-name of
\begin{align*}
(\bigvee\limits_{j=0}^{\lceil a_1 n\rceil-1}(f_{1,\omega}^j)^{-1}\alpha^{(1)}\vee \bigvee\limits_{j=\lceil a_{1}n\rceil}^{\lceil (a_1+a_{2}) n\rceil-1}(f_{1,\omega}^j)^{-1}\pi^{-1}\alpha^{(2)})
(x).
\end{align*}
Hence, we have $E_{n,\omega}\subset\bigcup\{B:B\in \mathcal G^\omega_n\},$ where $\mathcal G_n^\omega$
is the set of all elements
\begin{align*}
B\in \bigvee\limits_{j=0}^{\lceil a_1 n\rceil-1}(f_{1,\omega}^j)^{-1}\alpha^{(1)}\vee \bigvee\limits_{j=\lceil a_{1}n\rceil}^{\lceil (a_1+a_{2}) n\rceil-1}(f_{1,\omega}^j)^{-1}\pi^{-1}\alpha^{(2)}
\end{align*}
and
\begin{align*}
\mu_\omega(B)\leq  \exp((-(a_1 h^{(r)}_\mu(f_1,\alpha^{(1)})+a_2 h^{(r)}_{\mu\circ\Pi^{-1}}(f_2,\alpha^{(2)}))+\delta)n)
\end{align*}
and the $(\omega,\{\alpha^{(i)}\}_{i=1}^2,{\bf a},n)$-name of $B$ is $\delta_1$-close to the $(\omega,\{\alpha^{(i)}\}_{i=1}^2,{\bf a},n)$-name of $A$  for some $A\in\mathcal F_n^\omega.$

Since for each $A\in\mathcal F_n^\omega,$ the total number of $B$ in
\begin{align*}
\bigvee\limits_{j=0}^{\lceil a_1 n\rceil-1}(f_{1,\omega}^j)^{-1}\alpha^{(1)}\vee \bigvee\limits_{j=\lceil a_{1}n\rceil}^{\lceil (a_1+a_{2}) n\rceil-1}(f_{1,\omega}^j)^{-1}\pi^{-1}\alpha^{(2)}
\end{align*}
whose $(\omega,\{\alpha^{(i)}\}_{i=1}^2,{\bf a},n)$-name is $\delta_1$-close to the $(\omega,\{\alpha^{(i)}\}_{i=1}^2,{\bf a},n)$-name of $A,$ is upper bounded by
\begin{align*}
\binom{\lceil(a_1+a_2)n\rceil}{\lceil\lceil(a_1+a_2)n\rceil\delta_1\rceil}M^{\lceil\lceil(a_1+a_2)n\rceil\delta_1\rceil}
\leq\exp((a_1+a_2)\delta n+C+\delta).
\end{align*}
Combining this with \eqref{third estimation for fiber measure}, we have
\begin{align*}
\#\mathcal G_n^\omega&\leq\exp((a_1+a_2)\delta n+C+\delta)\cdot(\#\mathcal F_n^\omega)\\&\leq \exp((a_1 h^{(r)}_\mu(f_1,\alpha^{(1)})+a_2 h^{(r)}_{\mu\circ\Pi^{-1}}(f_2,\alpha^{(2)}))n-2\delta n+C+\delta).
\end{align*}
Moreover, by  and the definition of $\mathcal G_n^\omega$, the $\mu_{\omega}$-measure of the set of points $x\in E_{n,\omega}$ does not exceed
$\exp(-(a_1 h^{(r)}_\mu(f_1,\alpha^{(1)})+a_2 h^{(r)}_{\mu\circ\Pi^{-1}}(f_2,\alpha^{(2)}))n+\delta n)\cdot\#\mathcal G^\omega_n$. It follows that $\mu(E_{n})\leq\exp(-\delta n+C+\delta).$ Take $l_2\geq l$ so that $\sum_{n=l_2}^\infty\exp(-\delta n+C+\delta)<\delta.$
Then $\mu(\bigcup_{n\geq l_2}E_n)<\delta.$ Let $D=E\setminus\bigcup_{n\geq l_2}E_n.$
Then $\mu(D)>1-3\delta.$ For $(\omega,x)\in D$ and $n\geq l_2,$ one has
\begin{align*}
\mu_\omega(B_\omega^{\bf a}(x,n,\epsilon))\leq& \exp((a_1+a_2)\delta n+C+\delta)\\&\cdot\exp(-(a_1 h^{(r)}_\mu(f_1,\alpha^{(1)})+a_2 h^{(r)}_{\mu\circ\Pi^{-1}}(f_2,\alpha^{(2)}))n+(2+a_1+a_2)\delta n).
\end{align*}
Thus for $(\omega,x)\in D,$
\begin{align*}
&\liminf\limits_{n\to\infty}\dfrac{-\log\mu_{\omega}(B_\omega^{\bf a}(x,n,\epsilon))}{n}\\ \geq&  a_1 h^{(r)}_\mu(f_1,\alpha^{(1)})+a_2 h^{(r)}_{\mu\circ\Pi^{-1}}(f_2,\alpha^{(2)})-2(1+a_1+a_2)\delta\\
\geq&\min\left\{\frac{1}{\delta},h_{\mu}^{\bf a}(f_{1})-\delta\right\}-2(1+a_1+a_2)\delta.
\end{align*}
This ends the proof.
\end{proof}

\section{Appendix C: weighted Katok entropy formula}\label{Appendix C: a weighted Katok entropy formula}

By adopting the method of Katok, we are going to give an alternative definition of weighted measure-theoretic entropy $h_{\mu}^{\bf a}(f_{1})$ of RDSs. For $\mu\in \mathcal{M}_{\mathbb{P}}^{1}(\Omega\times X_{1},f_{1})$, $n\in\mathbb{N}$, $\epsilon>0$ and $0<\delta<1$, we call $K\subset X_{1}$ an $(\omega,{\bf a},n,\epsilon,\delta)$ covering set of $X_{1}$, if
\begin{align*}
\mu_{\omega}(\bigcup_{x\in K}B^{\bf a}_{\omega}(x, n,\epsilon))>1-\delta.
\end{align*}
Let $N_{f_{1}}(\omega,{\bf a},n,\epsilon,\delta)$ denote the smallest cardinality of $(\omega,{\bf a},n,\epsilon,\delta)$ covering sets of $X_{1}$.

\begin{definition}\label{an alternative definition of weighted measure-theoretic entropy}
Let $\mu\in \mathcal{E}_{\mathbb{P}}^{1}(\Omega\times X_{1},f_{1})$ and $0<\delta<1$, define
\begin{align*}
{\rm ent}_{\mu}^{\bf a}(\omega,f_{1})=\lim_{\epsilon\rightarrow0}\limsup_{n\rightarrow\infty}\frac{1}{n}\log N_{f_{1}}(\omega,{\bf a},n,\epsilon,\delta).
\end{align*}
\end{definition}

The main result of this Appendix C is the following:

\begin{theorem}\label{weighted Katok entropy formula}
Let $f_{i},i=1,2$ be continuous bundle RDSs over $(\Omega,\mathcal{F},\mathbb{P},\vartheta)$. Assume that ${\bf a}=(a_{1},a_{2})\in\mathbb{R}^{2}$ with $a_{1}>0$ and $a_{2}\geq0$, $f_{2}$ is a factor of $f_{1}$ with a factor map $\Pi:\Omega\times X_{1}\to \Omega\times X_{2}$. Let $\mu\in \mathcal{E}_{\mathbb{P}}^{1}(\Omega\times X_{1},f_{1})$ and $0<\delta<1$, then for $\mathbb{P}$-a.e. $\omega\in\Omega$, we have
\begin{align*}
\lim_{\epsilon\rightarrow0}\liminf_{n\rightarrow\infty}\frac{1}{n}\log N_{f_{1}}(\omega,{\bf a},n,\epsilon,\delta)&=\lim_{\epsilon\rightarrow0}\limsup_{n\rightarrow\infty}\frac{1}{n}\log N_{f_{1}}(\omega,{\bf a},n,\epsilon,\delta)\\&={\rm ent}_{\mu}^{\bf a}(\omega,f_{1})\\&=h_\mu^{{\bf a}}(f_{1}).
\end{align*}
\end{theorem}

\begin{remark}
Noted that Theorem \ref{weighted Katok entropy formula} shows that ${\rm ent}_{\mu}^{\bf a}(\omega,f_{1})$ is
independent of the choice of $\delta\in (0,1)$.
\end{remark}

Next, we prove the Theorem \ref{weighted Katok entropy formula}.

\begin{lemma}\label{power lemma}
Let $\mu\in \mathcal{E}_{\mathbb{P}}^{1}(\Omega\times X_{1},f_{1})$ and $r\in\mathbb{N}$. Then
\begin{align*}
{\rm ent}_{\mu}^{\bf a}(\omega,f_{1}^{r})\leq r{\rm ent}_{\mu}^{\bf a}(\omega,f_{1}),
\end{align*}
where the quantity ${\rm ent}_{\mu}^{\bf a}(\omega,f_{1}^{r})$ is defined
with respect to $\{f_{1}^{r},f_{2}^{r}\}$, and $f_{i}^{r}$ is the continuous bundle RDS over $(\Omega,\mathcal F,\mathbb P,\vartheta^{r})$ defined by $(f_{i}^{r})^{n}_{\omega}:=f_{i,\omega}^{rn}$ for any $i=1,2$.
\end{lemma}

\begin{proof}
For any $n\in\mathbb{N}$ and $a>0$, it is easy to verify that
\begin{align}\label{the inequality for lemma 8.4}
r\lceil a n\rceil-r+1\leq \lceil a n r\rceil\leq r\lceil a n\rceil.
\end{align}
Fix $n\in\mathbb{N}$, $\epsilon> 0$, $0<\delta<1$ and let $m\in\mathbb{N}$ such that $mr\leq n<(m+1)r$. Then
\begin{align*}
N_{f_{1}}(\omega,{\bf a},n,\epsilon,\delta)\leq N_{f_{1}}(\omega,{\bf a},(m+1)r,\epsilon,\delta).
\end{align*}
Then we have
\begin{align*}
\frac{1}{n}N_{f_{1}}(\omega,{\bf a},n,\epsilon,\delta)\leq \frac{1}{mr}N_{f_{1}}(\omega,{\bf a},(m+1)r,\epsilon,\delta).
\end{align*}
It follows that
\begin{align*}
\limsup_{n\rightarrow\infty}\frac{1}{n}N_{f_{1}}(\omega,{\bf a},n,\epsilon,\delta)\leq \limsup_{m\rightarrow\infty}\frac{1}{mr}N_{f_{1}}(\omega,{\bf a},mr,\epsilon,\delta).
\end{align*}
This implies that
\begin{align*}
\limsup_{n\rightarrow\infty}\frac{1}{n}N_{f_{1}}(\omega,{\bf a},n,\epsilon,\delta)= \limsup_{n\rightarrow\infty}\frac{1}{nr}N_{f_{1}}(\omega,{\bf a},nr,\epsilon,\delta).
\end{align*}
By \eqref{the inequality for lemma 8.4}, we have $N_{f_{1}^{r}}(\omega,{\bf a},n,\epsilon,\delta)\leq N_{f_{1}}(\omega,{\bf a},nr,\epsilon,\delta)$. Then we have
\begin{align*}
{\rm ent}_{\mu}^{\bf a}(\omega,f_{1}^{r})\leq r{\rm ent}_{\mu}^{\bf a}(\omega,f_{1}).
\end{align*}
\end{proof}

\begin{proof}[Proof of Theorem \ref{weighted Katok entropy formula}]
In the following, we divide the proof into two steps.

{\bf Step 1.} We prove that for $\mathbb{P}$-a.e. $\omega\in\Omega$, then
\begin{align*}
\lim_{\epsilon\rightarrow0}\limsup_{n\rightarrow\infty}\frac{1}{n}\log N_{f_{1}}(\omega,{\bf a},n,\epsilon,\delta)\leq h_\mu^{{\bf a}}(f_{1}).
\end{align*}

Let $\epsilon>0$. For any $i=1,2$, let $\xi^{(i)}$ be a finite measurable partition of $X_{i}$ and $\text{diam}(\xi^{(i)})<\epsilon$. Set $\hat{\xi}^{(i)}=\{\Omega\times A:A\in \xi^{(i)}\}$. Let $E\subset \Omega$ be a $\mathbb{P}$ full-measure set of elements satisfying \eqref{full measure-theoretic entropy} for $\hat{\xi}^{(i)}$
and $F\subset \Omega\times X_{1}$ a $\mu$ full-measure set of elements satisfying Theorem \ref{weighted SMB theorem of RDS}.  Fix $\omega\in E\cap \pi_{\Omega}F$, and $\gamma>0$. Set
\begin{align*}
X_{\omega,n}=\{x\in X_{1}:&\mu_{\omega}((\bigvee\limits_{j=0}^{\lceil a_1 n\rceil-1}(f^j_{1,\omega})^{-1}\xi^{(1)}\vee \bigvee\limits_{j=0}^{\lceil (a_1+a_{2}) n\rceil-1}(f^j_{1,\omega})^{-1}(\pi^{-1}\xi^{(2)}))(x))
\\&>\exp(-n(a_1 h_{\mu}^{(r)}(f_{1},\hat{\xi}^{(1)}\vee\Pi^{-1}\hat{\xi}^{(2)})+a_2 h_{\mu}^{(r)}(f_{1},\Pi^{-1}\hat{\xi}^{(2)})+\gamma))\}
\end{align*}
By Theorem \ref{weighted SMB theorem of RDS}, $\lim_{n\rightarrow\infty}\mu_{\omega}(X_{\omega,n})=1$. Then there exists $N\in\mathbb{N}$ such that for any $n\geq N$, we have $\mu_{\omega}(X_{\omega,n})>1-\delta$.
Note that $X_{\omega,n}$ contains at most $\lceil\exp(n(a_1 h_{\mu}^{(r)}(f_{1},\hat{\xi}^{(1)}\vee\Pi^{-1}\hat{\xi}^{(2)})+a_2 h_{\mu}^{(r)}(f_{1},\Pi^{-1}\hat{\xi}^{(2)})+\gamma))\rceil$ elements of
\begin{align*}
\bigvee\limits_{j=0}^{\lceil a_1 n\rceil-1}(f^j_{1,\omega})^{-1}\xi^{(1)}\vee \bigvee\limits_{j=0}^{\lceil (a_1+a_{2}) n\rceil-1}(f^j_{1,\omega})^{-1}(\pi^{-1}\xi^{(2)})
\end{align*}
and each element can be covered by some $B^{\bf a}_{\omega}(x, n,\epsilon)$. Therefore,
\begin{align*}
N_{f_{1}}(\omega,{\bf a},n,\epsilon,\delta)\leq \lceil\exp(n(a_1 h_{\mu}^{(r)}(f_{1},\hat{\xi}^{(1)}\vee\Pi^{-1}\hat{\xi}^{(2)})+a_2 h_{\mu}^{(r)}(f_{1},\Pi^{-1}\hat{\xi}^{(2)})+\gamma))\rceil.
\end{align*}
This implies that
\begin{align*}
\lim_{\epsilon\rightarrow0}\limsup_{n\rightarrow\infty}\frac{1}{n}\log N_{f_{1}}(\omega,{\bf a},n,\epsilon,\delta)\leq h_\mu^{{\bf a}}(f_{1}).
\end{align*}

{\bf Step 2.} We prove that for $\mathbb{P}$-a.e. $\omega\in\Omega$, then
\begin{align*}
\lim_{\epsilon\rightarrow0}\liminf_{n\rightarrow\infty}\frac{1}{n}\log N_{f_{1}}(\omega,{\bf a},n,\epsilon,\delta)\geq h_\mu^{{\bf a}}(f_{1}).
\end{align*}
For any $i=1,2$, let $\xi^{(i)}=\{A^{i}_{1},\ldots,A^{i}_{u_{i}}\}$ be a finite measurable partition of $X_{i}$ and $\kappa>0$. Therefore, there exist compact sets $B^{i}_{s}\subset A^{i}_{s}$, $1\leq s\leq u_{i}$, $i=1,2$ with
\begin{align*}
\pi_{X_{1}}\mu(A^{1}_{s}\setminus B^{1}_{s})&=\int\mu_{\omega}(A^{1}_{s}\setminus B^{1}_{s})d\mathbb{P}(\omega)<\frac{\kappa}{u_{1}\log u_{1}},\\
\pi_{X_{2}}\mu\circ\pi^{-1}(A^{2}_{s}\setminus B^{2}_{s})&=\int\mu_{\omega}\circ\pi^{-1}(A^{2}_{s}\setminus B^{2}_{s})d\mathbb{P}(\omega)<\frac{\kappa}{u_{2}\log u_{2}},
\end{align*}
Put $B_{0}^{i}=X_{i}\setminus\bigcup_{s=1}^{u_{i}}B^{i}_{s}$. Notice $\int\mu_{\omega}(B^{1}_{0})d\mathbb{P}(\omega)<\frac{\kappa}{\log u_{1}}$, $\int\mu_{\omega}\circ\pi^{-1}(B^{2}_{0})d\mathbb{P}(\omega)<\frac{\kappa}{\log u_{2}}$. Then for $\eta^{(i)}=\{B_{0}^{i},B_{1}^{i}, \ldots,B_{u_{i}}^{i}\}$, we have, for any $\omega\in \Omega$,
\begin{align*}
H_{\mu_{\omega}}(\xi^{(1)} \ | \ \eta^{(1)})&\leq \mu_{\omega}(B_{0}^{1})\log u_{1}<\kappa,\\
H_{\mu_{\omega}\circ\pi^{-1}}(\xi^{(2)} \ | \ \eta^{(2)})&\leq \mu_{\omega}\circ\pi^{-1}(B_{0}^{2})\log u_{2}<\kappa.
\end{align*}
Then we have
\begin{align*}
\int H_{\mu_{\omega}}(\xi^{(1)} \ | \ \eta^{(1)})d\mathbb{P}(\omega)<\kappa, \ \text{and} \ \int H_{\mu_{\omega}\circ \pi^{-1}}(\xi^{(2)} \ | \ \eta^{(2)})d\mathbb{P}(\omega)<\kappa.
\end{align*}
Since $\mathbb{P}$ is $\vartheta$-invariant, for any $M\in\mathbb{N}$, we have
\begin{align*}
\int H_{\mu_{\omega}}(\bigvee_{s=0}^{M-1}(f_{1,\omega}^{s})^{-1}\xi^{(1)})d\mathbb{P}(\omega)&\leq \int H_{\mu_{\omega}}(\bigvee_{s=0}^{M-1}(f_{1,\omega}^{s})^{-1}\eta^{(1)})d\mathbb{P}(\omega)+M\kappa,\\
\int H_{\mu_{\omega}\circ\pi^{-1}}(\bigvee_{s=0}^{M-1}(f_{2,\omega}^{s})^{-1}\xi^{(2)})d\mathbb{P}(\omega)&\leq \int H_{\mu_{\omega}\circ\pi^{-1}}(\bigvee_{s=0}^{M-1}(f_{2,\omega}^{s})^{-1}\eta^{(2)})d\mathbb{P}(\omega)+M\kappa,
\end{align*}
Set $\hat{\xi}^{(i)}=\{\Omega\times A:A\in \xi^{(i)}\}$ and $\hat{\eta}^{(i)}=\{\Omega\times B:B\in \eta^{(i)}\}$.
Hence we have
\begin{align*}
h_{\mu}^{(r)}(f_{1},\hat{\xi}^{(1)})&\leq h_{\mu}^{(r)}(f_{1},\hat{\eta}^{(1)})+\kappa,\\
h_{\mu\circ\pi^{-1}}^{(r)}(f_{2},\hat{\xi}^{(2)})&\leq h_{\mu\circ\pi^{-1}}^{(r)}(f_{2},\hat{\eta}^{(2)})+\kappa.
\end{align*}
Let $E\subset \Omega$ be a $\mathbb{P}$ full-measure set of elements satisfying \eqref{full measure-theoretic entropy} for both $\hat{\xi}^{(i)}$ and $\hat{\eta}^{(i)}$, and
$F\subset \Omega\times X_{1}$ a $\mu$ full-measure set of elements satisfying Theorem \ref{weighted SMB theorem of RDS}. Fix $\omega\in E\cap \pi_{\Omega}F$, and $\gamma>0$. Let
\begin{align*}
Y_{\omega,n}=&\{x\in X_{1}:\mu_{\omega}((\bigvee\limits_{j=0}^{\lceil a_1 n\rceil-1}(f^j_{1,\omega})^{-1}\eta^{(1)}\vee \bigvee\limits_{j=0}^{\lceil (a_1+a_{2}) n\rceil-1}(f^j_{1,\omega})^{-1}(\pi^{-1}\eta^{(2)}))(x))
\\&\leq\exp(-m(a_1 h_{\mu}^{(r)}(f_{1},\hat{\eta}^{(1)}\vee\Pi^{-1}\hat{\eta}^{(2)})+a_2 h_{\mu}^{(r)}(f_{1},\Pi^{-1}\hat{\eta}^{(2)})-\gamma)),\forall \ m\geq n\}.
\end{align*}
By Theorem \ref{weighted SMB theorem of RDS}, $\lim_{n\rightarrow\infty}\mu_{\omega}(Y_{\omega,n})=1$. Let
$\epsilon<\frac{1}{2}\min_{1\leq i\leq 2}\min_{1\leq j< l\leq u_{i}}\{d(B_{j}^{i},B_{l}^{i})\}$ where $d(A,B)=\inf_{x\in A,y\in B}d(x,y)$.
Then each $B^{\bf a}_{\omega}(x, n,\epsilon)$ can intersect at most $2^{bn+2}$ elements of $\bigvee\limits_{j=0}^{\lceil a_1 n\rceil-1}(f^j_{1,\omega})^{-1}\eta^{(1)}\vee \bigvee\limits_{j=0}^{\lceil (a_1+a_{2}) n\rceil-1}(f^j_{1,\omega})^{-1}(\pi^{-1}\eta^{(2)})$, where $b=a_{1}+a_{1}+a_{2}$. Then we have
\begin{align*}
&\mu_{\omega}(B^{\bf a}_{\omega}(x, n,\epsilon)\cap Y_{\omega,n})\\ \leq &2^{bn+2} \exp(-n(a_1 h_{\mu}^{(r)}(f_{1},\hat{\eta}^{(1)}\vee\Pi^{-1}\hat{\eta}^{(2)})+a_2 h_{\mu}^{(r)}(f_{1},\Pi^{-1}\hat{\eta}^{(2)})-\gamma))\\
\leq&\exp(-n(a_1 h_{\mu}^{(r)}(f_{1},\hat{\eta}^{(1)}\vee\Pi^{-1}\hat{\eta}^{(2)})+a_2 h_{\mu}^{(r)}(f_{1},\Pi^{-1}\hat{\eta}^{(2)})-\gamma-b\log 2)+2\log 2).
\end{align*}
We can assume that $N\in\mathbb{N}$ is sufficiently large that for any $n\geq N$, we have $\mu_{\omega}(Y_{\omega,n})>1-\delta$.

{\bf Claim.} For any $(\omega,{\bf a},n,\epsilon,\delta)$ covering set $K$ of $X_{1}$ and $n\geq N$, we have
\begin{align*}
\text{card}(K)\geq C\exp(n(a_1 h_{\mu}^{(r)}(f_{1},\hat{\eta}^{(1)}\vee\Pi^{-1}\hat{\eta}^{(2)})+a_2 h_{\mu}^{(r)}(f_{1},\Pi^{-1}\hat{\eta}^{(2)})-\gamma-b\log 2))
\end{align*}
for some constant $C>0$.

{\bf Proof of Claim.} Let $K'=\{x\in K:B^{\bf a}_{\omega}(x, n,\epsilon)\cap Y_{\omega,n}\neq\emptyset\}$. For $n\geq N$, we have
\begin{align*}
\text{card}(K)&\exp(-n(a_1 h_{\mu}^{(r)}(f_{1},\hat{\eta}^{(1)}\vee\Pi^{-1}\hat{\eta}^{(2)})\\&+a_2 h_{\mu}^{(r)}(f_{1},\Pi^{-1}\hat{\eta}^{(2)})-\gamma-b\log 2)+k\log 2)\\
\geq& \text{card}(K')\exp(-n(a_1 h_{\mu}^{(r)}(f_{1},\hat{\eta}^{(1)}\vee\Pi^{-1}\hat{\eta}^{(2)})\\&+a_2 h_{\mu}^{(r)}(f_{1},\Pi^{-1}\hat{\eta}^{(2)})-\gamma-b\log 2)+2\log 2)\\
\geq& \sum_{x\in K'}\mu_{\omega}(B^{\bf a}_{\omega}(x, n,\epsilon)\cap Y_{\omega,n})\\
\geq&\mu_{\omega}(\bigcup_{x\in K'}B^{\bf a}_{\omega}(x, n,\epsilon)\cap Y_{\omega,n})\\
=&\mu_{\omega}(\bigcup_{x\in K}B^{\bf a}_{\omega}(x, n,\epsilon)\cap Y_{\omega,n})\geq1-2\delta.
\end{align*}
This implies the Claim.

By the Claim, we have
\begin{align*}
\lim_{\epsilon\rightarrow0}\liminf_{n\rightarrow\infty}\frac{1}{n}\log N_{f_{1}}(\omega,{\bf a},n,\epsilon,\delta)\geq& a_1 h_{\mu}^{(r)}(f_{1},\hat{\eta}^{(1)}\vee\Pi^{-1}\hat{\eta}^{(2)})\\&+a_2 h_{\mu}^{(r)}(f_{1},\Pi^{-1}\hat{\eta}^{(2)})-\gamma-b\log 2.
\end{align*}
Note that
\begin{align*}
&a_1 h_{\mu}^{(r)}(f_{1},\hat{\eta}^{(1)}\vee\Pi^{-1}\hat{\eta}^{(2)})+a_2 h_{\mu}^{(r)}(f_{1},\Pi^{-1}\hat{\eta}^{(2)})\\
\geq& a_1 h_{\mu}^{(r)}(f_{1},\hat{\eta}^{(1)})+a_2 h_{\mu\circ\Pi^{-1}}^{(r)}(f_{2},\hat{\eta}^{(2)})\\
\geq& a_1 h_{\mu}^{(r)}(f_{1},\hat{\xi}^{(1)})+a_2 h_{\mu\circ\Pi^{-1}}^{(r)}(f_{2},\hat{\xi}^{(2)})-\kappa(a_1+a_{2}).
\end{align*}
Taking supremum over $\xi^{(i)}$ of $X_{i}$ for any $i=1,2$, we obtain that
\begin{align*}
\lim_{\epsilon\rightarrow0}\liminf_{n\rightarrow\infty}\frac{1}{n}\log N_{f_{1}}(\omega,{\bf a},n,\epsilon,\delta)\geq h_\mu^{{\bf a}}(f_{1})-\kappa(a_1+a_{2})-\gamma-b\log 2.
\end{align*}
This inequality clearly holds for $f_{1}^{r}$, $r\in\mathbb{N}$. Then by Lemma \ref{power lemma} and [\cite{Bogenschutz1992Entropy}, Theorem 3.6], we have
\begin{align*}
\lim_{\epsilon\rightarrow0}\liminf_{n\rightarrow\infty}\frac{1}{n}\log N_{f_{1}}(\omega,{\bf a},n,\epsilon,\delta)\geq h_\mu^{{\bf a}}(f_{1}).
\end{align*}
Since $\omega$ is chosen from the $\mathbb{P}$ full-measure set $E\cap \pi_{\Omega}F$, the inequality holds for $\mathbb{P}$ almost all $\omega$.

This completes the proof of the theorem.

\end{proof}

\section*{Acknowledgements}

The work was supported by the
National Natural Science Foundation of China (Nos. 12071222 and 11971236), China Postdoctoral Science Foundation (No.2016M591873),
and China Postdoctoral Science Special Foundation (No.2017T100384). The work was also funded by the Priority Academic Program Development of Jiangsu Higher Education Institutions.  We would like to express our gratitude to Tianyuan Mathematical Center in Southwest China (No.11826102), Sichuan University and Southwest Jiaotong University for their support and hospitality. The first author is also supported by Postgraduate Research \& Practice Innovation Program of Jiangsu Province (No. KYCX22\_1529).

\end{document}